\documentclass[12pt]{amsart}
\newtheorem{theo}{{\bfseries Theorem}}[section]
\newtheorem{prop}[theo]{{\bfseries Proposition}}
\newtheorem{lem}[theo]{{\bfseries Lemma}}
\newtheorem{cor}[theo]{{\bfseries Corollary}}
\newtheorem{df}[theo]{{\bfseries Definition}}
\newtheorem{ex}[theo]{{\bfseries Example}}

\input xy

\def \Lra {\Longrightarrow}

\def \ol {\overline}
\def \N {\mathbb N}

\def \Z {\mathbb Z}
\def \R {\mathbb R}
\def \A {\mathcal A}
\def \B {\mathcal B}
\def \E {\mathcal E}
\def \F {\mathcal F}
\def \G {\mathcal G}

\def \L {\mathcal L}
\def \O {\mathcal O}

\def \U {\mathcal U}

\def \f {\hat{f}}
\def \X {\hat{X}}

\usepackage{amssymb,latexsym}
\usepackage[mathcal]{eucal}
\usepackage{amscd}
\usepackage{amsmath}
\usepackage{graphicx}
\usepackage{amsfonts}
\usepackage{amscd}

\numberwithin{equation}{section}

\makeindex

\begin{document}

\title{\bfseries  Compactifications  of Dynamical Systems}

\author{Ethan Akin and Joseph Auslander}

\address{Mathematics Department \\
    The City College \\
    137 Street and Convent Avenue \\
    New York City, NY 10031, USA\\ \\
    Mathematics Department \\
    The University of Maryland \\
    College Park, MD 20742, USA \\}

\date{April, 2006, revised July, 2009}
\vspace{.5cm} \maketitle

{\bfseries Abstract:}  While compactness is an essential assumption for many results in
dynamical systems theory, for many applications the state space is only locally compact.
Here we provide a general theory for compactifying such systems, i.e. embedding them as
invariant open subsets of compact systems.  In the process we don't want to introduce recurrence
which was not there in the original system. For example if a point lies on an orbit which remains
in any compact set for only a finite span of time then the point becomes non-wandering if we use
the one-point compactification. Instead, we develop here the appropriate theory of dynamic
compactification.

\tableofcontents

\section*{Introduction}

In dynamical systems theory the state space is usually a compact metric
space.  Metrizability is usually just a convenience, but compactness is
essential for many arguments.  On the other hand, in many applications
the natural state space is not compact, e.g. a finite-dimensional vector space or an open
subset thereof.  Often local compactness is the best we can get to begin with.
For this reason much attention is given to finding compact invariant
subsets to which the system can be restricted.

Alternatively, we can seek to compactify the system.  That is, we include the original
as a subsystem of one occurring on a larger, compact space.  For example, any
homeomorphism or flow on a locally compact space extends to the one-point compactification.
However, from the dynamic point of view it is undesirable to
thus concatenate of all of the ``points at infinity"
to a single point. Doing so often introduces new recurrence relationships among points of
the original space. What we want are compactifications which introduce no new recurrence.

{\bfseries N. B.} All of the spaces we will consider are Hausdorff completely regular spaces,
i.e. spaces which admit a Hausdorff uniformity.  The state spaces $X, Y$ etc. for our dynamical systems
are locally compact and $\sigma$-compact. In particular,  $X$ is metrizable iff it has a
 countable base. We let $\B (X)$ denote the Banach algebra of bounded, continuous,
real-valued functions on $X$ with the sup norm.  \vspace{.5cm}

{\bfseries 1. Compactifications of a Closed Relation:} \quad While our primary interest
is in the dynamics of a map or a flow, we will follow Akin (1993) in considering
the dynamics of relations.  We regard a map from $X$ to $Y$ as a special case of a relation
$f$ from $X$ to $Y$, that is, a subset of $X \times Y$. We write $f : X \to Y$ for such a relation
letting $f(x)$ denote the -possibly empty- subset $\{ y \in Y : (x,y) \in f \}$. More generally
for $A \subset X$ we let $f(A) = \{ y : (x,y) \in f$ for some $ x \in A \} \subset Y$. We call
$f^{-1} = \{ (y,x) : (x,y) \in f \}$
the \emph{reverse relation} from $Y$ to $X$.

If $X = Y$ we will
call $f$ a relation on $X$. For example, the identity map, $1_X$, is the
relation $ \{ (x,x) : x \in X \}$. Composition generalizes to relations and so we can iterate a relation on $X$.
The closed relations, i.e. the closed subsets of $X \times Y$, are of special interest.
The relation $f$ is called \emph{$+$proper} when it is closed and $A \subset X$ compact implies $f(A) \subset Y$ is
compact and $f$ is \emph{proper} when both $f$ and $f^{-1}$ are proper. A map between locally compact
 spaces is continuous iff it is a  $+$proper relation. The composition of closed relations need
 not be closed but the composition of  $+$proper relations is  $+$proper.

In general,  a \emph{dynamical system} is a pair $(X,f)$ with
 $f$ a  closed relation on $X$. When $f$ is a map, the pair is called a \emph{cascade}.

Relations arise naturally in studying the dynamics of maps.  For example, for a point $x$ the omega
limit point set $\omega f (x) = lim sup_n \ \{ f^n(x) \}$ defines the relation $\omega f$.
The  orbit relation,
$\mathcal{O} f = \bigcup_{n=1}^{\infty} \ f^n $,  is always transitive but is usually not a closed relation.  Its
closure, denoted $\mathcal{N} f$ is closed but usually not transitive.
We define $\G f$, the generalized recurrence relation,
to be the smallest closed and transitive relation which  contains $f$.

We call $|f| = \{ x : (x,x) \in f \}$ the \emph{cyclic set} of $f$.  For $f$ a continuous map,
$|f|$ is the set of fixed points, $|\O f|$ is the set of periodic points, $|\omega f|$ is the
set of recurrent points, $|\mathcal{N} f|$ is the set of
non-wandering points and $|\G f|$ is called the set of generalized recurrent points.

It is easy to check that $\G (f^{-1}) = (\G f)^{-1}$ and so we can omit the parentheses.
$\G f \cap \G f^{-1}$ is a closed equivalence relation on $|\G f|$. Hence, this set is partitioned
by the closed $\G f \cap \G f^{-1}$ equivalence classes.

A \emph{compactification of a space} $X$ is  a continuous map
$j : X \to \X$ such that $j(X)$ is dense in $\X$. Hence, the induced map $j^* : \B(\X) \to \B(X)$ is
an injective isometry of Banach algebras. The compactification is classified by the subalgebra $j^*(\B(\X))$ of
$\B(X)$.  That is, if $\A$ is any closed subalgebra of $\B(X)$ there is a compactification $j : X \to \X$, unique
up to homeomorphism, such that $\A = j^*(\B(\X))$.  A compactification is called \emph{proper} when
$j$ restricts to a homeomorphism of $X$ onto  $j(X) \subset \X$. Since $X$ is locally compact, $j$ is proper
iff $j(X)$ is an open subset of $\X$ and $j : X \to j(X)$ is a proper injective map. For example,
the closed subalgebra $\A_0$ generated by the functions of compact support induces the one-point
compactification $X_*$ of $X$. The compactification $j$ is proper iff $\A$ separates points and closed sets and
$\X$ is metrizable iff $\A$ is countably generated.

When $j$ is proper
we  will usually identify $X$ with $j(X)$ and so regard $j$ as the inclusion of $X$ as an open subset of $\X$.

A \emph{compactification of a dynamical system} $(X,f)$ is a pair $(\X,\f)$ with $j: X \to \X$ a compactification of
$X$ and $\f$ the closure in $\X \times \X$ of $(j \times j)(f)$. If $X \subset \X$ is a proper compactification
then $\f$ is the closure of $f$ and since $f$ is closed, $f = \f \cap (X \times X)$.
 $(\X,\f)$ is a \emph{dynamic compactification} when $\G f = \G \f \cap (X \times X)$.

 \begin{theo} Let $(X,f)$ be a dynamical system with $f$  +proper.
Let  $ (\hat{X},\hat{f})$ be a dynamic
compactification of  $(X,f)$.

If $\hat{E} \subset |\mathcal{G} \hat{f}| $ is a $\mathcal{G}
\hat{f} \cap \mathcal{G} \hat{f}^{-1}$ equivalence class  with $E =
\hat{E} \cap X$, then exactly one of the following three possibilities
holds:
\begin{itemize}
\item[(i)]  $\hat{E}$ is a compact subset of $\hat{X} \setminus
X$ and $E = \emptyset$.
\item[(ii)] $E$ is contained in $ |\mathcal{G}f|$  and is a noncompact
$\mathcal{G} f \cap \mathcal{G}f^{-1}$ equivalence class whose
$\hat{X}$  closure meets $\hat{X}\setminus X$ and is contained in
$\hat{E}  $.
\item[(iii)] $\hat{E} = E$ is contained in $|\mathcal{G}f|$ and is a compact
$\mathcal{G} f \cap \mathcal{G} f^{-1} $ equivalence class.
\end{itemize}

Furthermore, if $x,y \in |\mathcal{G}f|$  lie in distinct $\mathcal{G}f \cap
\mathcal{G}f^{-1}$ equivalence classes then their equivalence
classes have disjoint closures in $\hat{X}$.
\end{theo} \vspace{.5cm}

{\bfseries 2. Lyapunov Function Compactifications:} \quad We construct dynamic compactifications
by using Lyapunov functions.

A  function $L \in \B(X)$
is called a Lyapunov function for a relation $f$ on $X$ if $y \in f(x)$
implies $L(x) \leq L(y)$.  That is, the relation $f$ is contained in the relation
$\leq_L = \{ (x,y) : L(x) \leq L(y) \}$. Since $\leq_L$ is a closed, transitive relation,
it follows that a Lyapunov
function for $L$ is automatically a Lyapunov function for $\G f$.

 A collection $\L$ of Lyapunov functions for $(X,f)$ is {\em sufficient set of Lyapunov functions}
when
\begin{displaymath}
1_X \cup \mathcal{G}f  \quad = \quad \bigcap_{\mathcal{L}} \{\leq_L \}.
\end{displaymath}

\begin{theo}Let $(X,f)$ be a dynamical system with $f$ a +proper relation.

The set of all Lyapunov functions is sufficient.

If $\L$ is a sufficient set of Lyapunov functions and $\A$ is a the closed subalgebra of $\B(X)$ generated
by $\L$ together with the functions of compact support, then the associated compactification $(\X, \f)$ is
a proper dynamic compactification, called the $\L$ compactification.

Let $\L$ be a sufficient set of Lyapunov functions and $(\X,\f)$ be the associated $\L$ compactification.
If $\hat{E} \subset |\mathcal{G} \hat{f}| $ is
 a $\mathcal{G} \hat{f} \cap
\mathcal{G} \hat{f}^{-1}$ equivalence class with $E = \hat{E} \cap
X$, then exactly one of the following four possibilities holds:
\begin{itemize}
\item[(i)]  $\hat{E}$ consists of a single point of $\hat{X} \setminus
X$.
\item[(ii)] $E$ is contained in $ |\mathcal{G}f|$  and is a noncompact
$\mathcal{G} f \cap \mathcal{G}f^{-1}$ equivalence class with
$\hat{E}  $ its one point compactification. That is, there is a
noncompact equivalence class $E \subset |\mathcal{G} f|$ whose
closure in $\hat{X}$ is $\hat{E}$ and $\hat{E} \setminus E$ is a
singleton.
\item[(iii)] $\hat{E} = E$ is contained in $|\mathcal{G}f|$ and is a compact
$\mathcal{G} f \cap \mathcal{G} f^{-1} $ equivalence class.
\end{itemize}

If $X$ is metrizable there exists a countable sufficient set $\L$ of Lyapunov functions and the space
$\X$ of the $\L$ compactification is metrizable.
\end{theo}\vspace{.5cm}

If $(X,f)$ is a cascade, i.e. $f$ is a continuous map on $X$, then $(\X,\f)$ is called a cascade
compactification when the closed relation $\f$ is a map (necessarily continuous) on $\X$. If $f$ is a
proper continuous map and $\L$ is a sufficient set of Lyapunov functions such that $L \in \L$ implies
$L \circ f^n \in \L$ for every $n \in \Z_+$ then the $\L$ compactification is a dynamic cascade
compactification. If $f$ is a homeomorphism on $X$ and $L \circ f^n \in \L$ for every $n \in \Z$ then
$\f$ is a homeomorphism on $\X$.

The set of Lyapunov functions for a closed relation $f$ distinguishes the points of $X$ iff
$\G f \cap \G f^{-1} \subset 1_X$.
If a +proper relation $f$ on $X$ satisfies $\G f \cap \G f^{-1} \subset 1_X$
then the set of Lyapunov functions for
$f$ determines the topology of $X$.  That is, if $\{ x_i \}$ is a net in $X$ and $x \in X$ then
if $\{ L(x_i) \}$ converges to $L(x)$ for every $f$ Lyapunov function $L$ then $\{ x_i \}$ converges
to $x$ in $X$.

\vspace{.5cm}

{\bfseries 3. Compactifications of a Flow:} A compactification for a flow $\phi: \R \times X \to X$ is, to begin with,
a proper compactification $\X \supset X$ whose associate algebra $\A$ is $\phi^*$ invariant.  That is,
$\A = (\phi^t)^*(\A)$ for all $t \in \R$.  This ensures, that each map $\phi^t$ extends to
a homeomorphism $\hat \phi^t$ of $\X$. However, to ensure that $\hat \phi : \R \times \X \to \X$ is continuous
we need that $\A \subset \B_{\phi}(X)$ where $u \in \B_{\phi}(X)$ when  the function
$t \mapsto u \circ \phi^t $ is a uniformly continuous function from $\R$ to $\B(X)$.

For $K$ a compact
subset of $\R$, we let $\phi^K$ denote the closed relation $\bigcup_{t \in K} \phi^t$. We use this especially
with $K = [0,1] = I$ and $K = [1,2] = J$.  Thus, for example,
\begin{displaymath}
\mathcal{O} \phi^I \quad = \quad \phi^I \cup \mathcal{O} \phi^J \quad = \quad \bigcup_{t \in \R_+} \phi^t. \hspace{1cm}
\end{displaymath}
Letting $\G \phi = \G \phi^I$ we show that $\G \phi = \phi^I \cup \G \phi^J$. Observe that $\G \phi^I$ is reflexive
and so we use $|\G \phi^J|$ to define generalized recurrent points.  On $|\G \phi^J|$ the two equivalence
relations $\G \phi^J \cap (\G \phi^J)^{-1}$ and $\G \phi \cap \G \phi^{-1}$ agree. Off  $|\G \phi^J|$ the
latter is $1_X$.

A compactification for the flow $\phi$ is called dynamic when it is dynamic for the proper relation $\phi^J$.

A  function $L \in \B(X)$ is a Lyapunov function for the flow when $L(\phi(t,x))$ is nondecreasing in $t$
for every $x$. A collection $\L$ of Lyapunov functions for the flow
is called a \emph{sufficient set of Lyapunov functions} for
the flow when $\L \subset \B_{\phi}(X)$  and
\begin{displaymath}
\bigcap_{L \in \L} \ \leq_L \quad = \quad \G \phi. \hspace{2cm}
\end{displaymath}
From a $\phi^*$ invariant, sufficient set of Lyapunov functions $\L$ for the flow
we obtain a Lyapunov function compactification of the flow by using the algebra generated
by $\L$ and the functions of compact support.

We will use $|\phi|$ to denote the set of fixed points of the
semiflow $\phi$.  That is,
\begin{displaymath}
|\phi| \quad =_{def} \quad \{ x \in X : \phi^t(x) = x \quad
\mbox{for all} \ t \in {\mathbb R}_+ \}.
\end{displaymath}

\begin{theo} Let $\phi$ be a flow on  $X$.
Lyapunov function compactifications for the flow $\phi$ exist. Let $\hat{\phi}$ on $\hat{X}$ be such
a Lyapunov compactification for $\phi$.

\begin{itemize}
\item[(a)] $\hat{X}$ is a dynamic compactification for $\phi$ with
\begin{displaymath}
\begin{split}
(X \times X) \cap \mathcal{G} \hat{\phi} \quad = \quad
\mathcal{G}\phi \hspace{2cm}\\
(X \times X) \cap \mathcal{G} (\hat{\phi}^J) \quad = \quad
\mathcal{G}(\phi^J).\hspace{1cm}
\end{split}
\end{displaymath}
\item[(b)]  The compact set $\hat{X} \setminus X$ is $\hat \phi $ invariant and
every generalized recurrent point of $\hat{\phi}$ which lies
in $\hat{X} \setminus X$ is a fixed point for $\hat{\phi}$.  That
is,
\begin{displaymath}
(\hat{X} \setminus X) \cap |\mathcal{G}\hat{\phi}^J| \qquad \subset
\qquad |\hat{\phi}|
\end{displaymath}
\item[(c)] If $\hat{E} \subset |\mathcal{G} (\hat{\phi}^J)| $
  is an $\mathcal{G} (\hat{\phi}^J) \cap
\mathcal{G} (\hat{\phi}^J)^{-1}$  equivalence class with $E =
\hat{E} \cap X$ then exactly one of the following three
possibilities holds:
\begin{enumerate}
\item[(i)]  $\hat{E}$ consists of a single point of $\hat{X} \setminus
X$ which is a fixed point for $\hat{\phi}$ and $E = \emptyset$.
\item[(ii)] $\hat{E}$ is the one point compactification of a noncompact
$\mathcal{G}(\phi^J) \cap \mathcal{G}(\phi^J)^{-1}$ equivalence
class $E$. That is, there is a noncompact $\mathcal{G}(\phi^J)
\cap \mathcal{G}(\phi^J)^{-1}$ equivalence class $E \subset X$
whose closure is $\hat{E}$ and $\hat{E} \setminus E$ is a
singleton which is a fixed point of $\hat{\phi}$.
\item[(iii)] $\hat{E}$ is contained in $ X$, i.e. $\hat{E} = E$, and it is a  compact
$\mathcal{G}(\phi^J) \cap \mathcal{G}(\phi^J)^{-1}$ equivalence
class.
\end{enumerate}
\item[(d)]  For $x \in X$
the $\mathcal{G} \hat{\phi} \cap \mathcal{G} \hat{\phi}^{-1}$
equivalence class of $x$ is the closure in $\hat{X}$ of its
$\mathcal{G}\phi \cap \mathcal{G}\phi^{-1}$ equivalence class.

\item[(e)] If $X$ is metrizable  then there exist metrizable Lyapunov function
compactifications for $\phi$.
\end{itemize}
\end{theo}
\vspace{.5cm}

{\bfseries 4. Chain Compactifications:} The chain relation $\mathcal{C} f$ is a uniform notion and
so we require here that our space $X$ is equipped with a uniformity $\U_{X}$. Recall that a compact space
has a unique uniformity consisting of all of the neighborhoods of the diagonal $1_X$. We let $\B_{\U}(X)$
denote the closed subalgebra consisting of those elements of $\B(X)$ which are uniformly continuous.

For a closed +proper
relation $f$ on $X$, we define $\mathcal{C} f = \bigcap_{V \in \U_X} \ \O V \circ f$. It is a closed,
transitive relation which contains $f$ and so contains $\G f$. A proper compactification $(\X, \f)$ of
$(X,f)$ is called \emph{chain dynamic} when the inclusion map $j : X \to \X$ is uniformly continuous and
$(X \times X) \cap \mathcal{C} \f = \mathcal{C} f$.

As usual the dynamic compactifications are constructed using Lyapunov functions.  $L : X \to [0,1]$ is
an \emph{elementary uniform Lyapunov function} for $f$ when $L \in \B_{\U}(X)$ and
\begin{displaymath}
(x,y) \in f \qquad \Longrightarrow \qquad L(x) = 0 \quad \mbox{or} \quad L(y) = 1.
\end{displaymath}
An elementary uniform Lyapunov function for $f$ is a Lyapunov function for $\mathcal{C} f$. If $f$ is
a uniformly continuous map then $L \circ f^n$ is an elementary uniform Lyapunov function for all
$n \in \Z_+$ and if $f$ is a uniform isomorphism then the same is true for all $n \in \Z$.

A set $\L$ of elementary uniform Lyapunov functions is called a \emph{sufficient set} when
\begin{displaymath}
1_X \cup \mathcal{C} f \quad = \quad \bigcap_{L \in \L} \ \leq_L.
\end{displaymath}
For a closed +proper relation $f$ the set of all elementary uniform Lyapunov functions is a sufficient set
and if $X$ is second countable, or, equivalently the topological space $X$ is metrizable, then there
is a countable sufficient set.

\begin{theo} Let  $f$ be a + proper, closed
relation on a uniform space $X$. Let $\mathcal{L} \subset \mathcal{B}_{\U}(X)$ be a sufficient
set of elementary uniform Lyapunov functions for $f$ and
 $\mathcal{A}$ be the closed subalgebra generated by
$\mathcal{L}$ and the continuous functions with compact support.
If $(\X, \f)$ is the  $\L$ compactification of the dynamical system
$(X,f)$ then $(\X, \f)$ is a chain dynamic compactification of $(X,f)$.

Furthermore, if $\hat{E} \subset |\mathcal{C} \hat{f}| $ is
 a $\mathcal{C} \hat{f} \cap
\mathcal{C} \hat{f}^{-1}$ equivalence class with $E = \hat{E} \cap
X$, then exactly one of the following three possibilities holds:
\begin{itemize}
\item[(i)]  $\hat{E}$ consists of a single point of $\hat{X} \setminus
X$.
\item[(ii)] $E$ is contained in $ |\mathcal{C}f|$  and is a noncompact
$\mathcal{C} f \cap \mathcal{C}f^{-1}$ equivalence class with
$\hat{E} \ $its one point compactification. That is, there is a
noncompact equivalence class $E \subset |\mathcal{C} f|$ whose
closure in $\hat{X}$ is $\hat{E}$ and $\hat{E} \setminus E$ is a
singleton.
\item[(iii)] $\hat{E} = E$ is contained in $|\mathcal{C}f|$ and is a compact
$\mathcal{C} f \cap \mathcal{C} f^{-1} $ equivalence class.
\end{itemize}

If $f$ is a uniformly continuous proper map and $\L$ is $f^*$ +invariant then
$(\X, \f)$ is a cascade compactification. If $f$ is a uniform isomorphism and
$\L$ is $f^*$ invariant then $(\X, \f)$ is reversible.
\end{theo}

\vspace{.5cm}

{\bfseries 5. Stopping at Infinity:} By using Beck's ideas for rescaling time for an arbitrary
flow, we are able to construct Lyapunov compactifications $(\X, \hat \phi)$
of a flow $(X, \phi)$ such that every point of $\X \setminus X$ is a fixed point for the flow
$\hat \phi$.  By using the suspension construction we are able to do the same thing for a
homeomorphism, i.e. a reversible cascade. \vspace{.5cm}

{\bfseries 6. Parallelizable Systems:} Following Antosiewicz and
Dugundji as well as Markus, we characterize parallelizable flows.  That is, flows $(X, \phi)$ which
are isomorphic to the product of a constant flow and the translation flow on $\R$.

\begin{theo} Let $\phi$ be a flow on $X$.
\begin{enumerate}
\item[(a)] The following conditions are equivalent:
\begin{itemize}
\item[(i)] The reflexive, transitive relation $\mathcal{O}\phi$ is
closed.
\item[(ii)] The transitive relation $\mathcal{O}(\phi^J)$ is closed.
\item[(iii)] $\mathcal{O}\phi = \mathcal{N}\phi$.
\item[(iv)] $\mathcal{O}\phi = \mathcal{G}\phi$.
\item[(v)] $\mathcal{O}(\phi^J) = \mathcal{G}(\phi^J)$.
\end{itemize}
The above conditions imply that the equivalence relation
$\mathcal{O}(\phi \cup \phi^{-1})$ is closed, or, equivalently,
 $\mathcal{O}(\phi \cup \phi^{-1}) = \mathcal{G}(\phi \cup
\phi^{-1})$.

\item[(b)]  The following are equivalent:
\begin{itemize}
\item[(i)] $\mathcal{O}\phi$ is closed and there are no periodic points, i.e. $|\mathcal{O}(\phi^J)|
= \emptyset$.
\item[(ii)] $\mathcal{O}\phi$ is closed and all points are wandering, i.e. $|\mathcal{N}(\phi^J)|
= \emptyset$.
\item[(iii)] $\mathcal{O}\phi$ is closed and there are no generalized recurrent points, i.e. $|\mathcal{G}(\phi^J)|
= \emptyset$.
\item[(iv)] $\mathcal{O}(\phi \cup \phi^{-1})$ is closed and all points are wandering, i.e. $|\mathcal{N}(\phi^J)|
= \emptyset$.
\item[(v)] $\phi$ is parallelizable.
\end{itemize}
\end{enumerate}
\end{theo}
\vspace{.5cm}

{\bfseries 7. Appendix: Limit Prolongation Relations:} We prove several useful identities connecting
various prolongations of closed relations.
\vspace{.5cm}

{\bfseries 8. Appendix: Paracompactness:}  While we have considered only those locally compact spaces
which are $\sigma$-compact, the results actually apply to locally compact spaces which are paracompact
because such spaces are disjoint unions of clopen subsets each of which is $\sigma$-compact. Furthermore,
if a closed relation $f$ on the space is proper, then members of the partition by $\sigma$-compact
clopen sets can be taken to be $f$ +invariant.
\vspace{.5cm}

\section{Compactifications of a Closed Relation}\label{seccomp}

Except for various Banach algebras, our spaces $X$ are all assumed
to be locally compact, $\sigma$-compact Hausdorff spaces. Thus,
when they are metrizable they are separable and so have a
countable base. A subset $A$ of $X$ is called \emph{bounded}
\index{bounded set} when
it has compact closure.  For $A \subset X$  we write as usual  $\overline{A}$ and $A^{\circ}$
for the closure and interior, respectively. For arbitrary subsets $A, B$ of $X$ we
will use $A \subset \subset B$
\index{$\subset \subset$}
to mean that  $\ol{A} \subset B^{\circ}$.

 ${\mathbb R}_+$ is the subset of
nonnegative elements of ${\mathbb R }$, the set of real numbers.
With ${\mathbb Z}$ the set of integers, ${\mathbb Z}_+ = {\mathbb
Z} \cap {\mathbb R}_+$ and $\N$ is the set of positive integers.

A \emph{relation} \index{relation} $f : X \to Y$ for sets $X$ and $Y$ is just a subset of $X \times Y$. For
$x \in X$ and $A \subset X$  we write $f(x) = \{ y : (x,y) \in f
\}$ and $f(A) = \bigcup \{ f(x) : x \in A \}$. When $f$ is a
 map we use the same notation for the singleton set
$f(x)$ and for the point it contains. The \emph{reverse relation}\index{relation!reverse relation}
$f^{-1} : Y \to X$ is
defined by
$f^{-1} = \{ (y,x) : (x,y) \in f \}$.

For relations $f : X \to Y, g : Y \to Z$
the \emph{composition} \index{relation!composition}$g \circ f : X \to Z$ is defined by $g \circ f = \{ (x,z) :$
there exists $y \in Y$ such that $ (x,y) \in f $ and
$(y,z) \in g   \}$. As with maps the composition operation is
associative.  The identity maps like $1_X = \{ (x,x) : x \in X \}$ act as
identities with respect to composition. Clearly, $(g \circ f)^{-1} = f^{-1} \circ g^{-1}$.

When $X$ and $Y$ are  spaces then $f : X \to Y$ is called a \emph{closed relation} \index{relation!closed relation}
when it is a
closed subset of $X \times Y$.
If $f : X \to Y$ is a  relation between  spaces, we
call $f$ +\emph{proper} \index{relation!+proper relation}if it is closed and if $A \subset X$ compact implies that
$f(A) \subset Y$ is compact. $f$ is called \emph{proper} \index{relation!proper relation}when both
$f$ and $f^{-1}$ are +proper. A continuous map is a +proper relation and it is
proper as a relation iff it is a proper continuous map in the usual sense, i.e. $B$ compact
implies $f^{-1}(B)$ is compact.

If $X_0$ is compact then the projection map $\pi_2 : X_0 \times Y
\to Y$ is a proper continuous map.  Furthermore, if $U \subset X_0 \times Y$
is an open subset then
\begin{equation}\label{eq1.2}
\{ y \in Y :  X_0 \times \{ y \} \ \subset \ U \}
\end{equation}
is an open set by Wallace's Lemma (Kelley (1955) Theorem 5.12).  Hence,
if $A$ is a closed subset of $X_0 \times Y$ then the image
$\pi_2(A)$ is closed, i.e. $\pi_2$ is a closed map.

We collect some elementary results we will need about such relations.

\begin{prop}\label{prop1.1}  Let $f: X \to Y$ be a closed relation  and let $A$
be an arbitrary subset of $X$.
\begin{equation}\label{eq1.1}
f(A) \quad = \quad \bigcap \ \{ f(U) : U \ \mbox{is open and} \ A
\subset U \}.
\end{equation}
\end{prop}

{\bfseries Proof:}  Clearly, $f(A)$ is contained in the
intersection.  If $y \not\in f(A)$ then $f^{-1}(y)$ is a closed
set disjoint from $A$ and hence $U = X \setminus f^{-1}(y)$ is an
open set containing $A$ with $y \not\in f(U)$.

$\Box$ \vspace{.5cm}

\begin{prop}\label{prop1.2}Let $X, Y$ and $Z$ be   spaces.
Assume   $f : X \to Y$  and $g : Y \to Z$ are closed
relations.
  \begin{itemize}
  \item[(a)]  If $A \subset X$ is compact then $f(A) \subset Y$ is
  closed.
  \item[(b)]   Assume that $\mathcal{A}$ is a filterbase of closed subsets
  of $X$
  with intersection $C$. If either (i)
  for each $y \in Y$, $f^{-1}(y)$ is  compact
   ( e.g.  $f^{-1}$ is + proper), or (ii) $A \in \mathcal{A}$ implies $A$ is compact, then
  \begin{equation}\label{eq1.3}
\bigcap_{ A \in \mathcal{A}} \ f(A) \quad = \quad f(C).
\end{equation}
  \item[(c)]  Assume that  $f$ is + proper. If $B \subset Y$ is closed, then $f^{-1}(B) \subset X$ is closed.
  If $U \subset Y$ is open then $\{ x \in X : f(x) \subset U \}$ is open.
  \item[(d)]  If either
  $f$ or $g^{-1}$ is + proper then $g \circ f$ is a closed
  relation.
  \item[(e)]  If both
  $f$ and $g$ are + proper (or proper) relations then $g \circ f$
  is a + proper (resp. proper)  relation.
  \item[(f)] If $f$ is a mapping then it is a continuous map iff it is a $+$proper relation.
  \end{itemize}

   Assume $\mathcal{F}$ is a filterbase of closed
relations from $X$ to $Y$ with intersection $f$ and that
$\mathcal{A}$ is a filterbase of closed subsets of $X$ with
intersection $C$. Assume $U$ is an open containing $f(C)$.
\begin{itemize}
\item[(g)] If either (i)  $g \in \F$ implies $g^{-1}$ is + proper, or
(ii) $A \in \A $ implies $A$ is compact, then
\begin{equation}\label{eq1.4}
\bigcap_{g \in \mathcal{F}, A \in \mathcal{A}} \ g(A) \quad =
\quad f(C).
\end{equation}

\item[(h)] Assume that the complement $X \setminus U$ is compact.
If either (i)  $g \in \F$ implies $g^{-1}$ is + proper, or
(ii) $A \in \A $ implies $A$ is compact, then there exist $g \in
\mathcal{F}$ and $A \in \mathcal{A}$ such that $U$ contains
$g(A)$.

\item[(i)] If $A \in \A$ implies $A$ is compact and
 $g \in \mathcal{F}$ implies $g$ is + proper, then there exist
$g \in \mathcal{F}$ and $A \in \mathcal{A}$ such that $U$ contains
$g(A)$.

\end{itemize}
  \end{prop}

{\bfseries Proof:}(a): $f(A)$ is the image under the closed map $ \pi_2 : A
\times Y \to Y$ of the closed set $(A \times Y) \cap f$.

(b):  For $y$ in $\bigcap_A \ f(A), \ \{f^{-1}(y) \cap A  \}$ is a
filterbase of nonempty compacta and so the intersection,
$f^{-1}(y) \cap C$, is nonempty.

(c): Let $A$ be an arbitrary compact subset of $X$.  By assumption
$f(A)$ is compact in $Y$ and so $f(A) \cap B$ is
compact.  By (a) applied to $f^{-1}$, $A \cap f^{-1}(B)  = f^{-1}(f(A) \cap B) \cap A$
implies that $A \cap f^{-1}(B)$ is closed.  Because $A$ is arbitrary
and $X$ is locally compact, it follows that $f^{-1}(B)$ is closed. If $U$ is
open then $B = Y \setminus U$ is closed and $X \setminus f^{-1}(B) = \{ x : f(x) \subset U \}$.

(d): Let $A \subset X$ and $B \subset Z$ be compact. By (a)
both $f(A)$ and $ g^{-1}(B)$ are closed.  By assumption, at least one of
them is compact.  Hence, $C =_{def} f(A) \cap g^{-1}(B) \subset Y$
is compact.  Furthermore,
\begin{equation}\label{eq1.5}
(A \times B) \cap (g \circ f) \quad = \quad (A \times B) \cap
\pi_{13}((X \times C \times Z) \cap (f \times Z) \cap( X \times g))
\end{equation}
implies that $(A \times B) \cap (g \circ f)$ is closed.  As in (c), this implies that
$g \circ f$ is closed.

(e):  If $f$ and $g$ are +proper and $A$ is a compact subset of
$X$ then $g \circ f (A) = g(f(A))$ is compact. Since $g \circ f$
is closed by (d), it follows that $g \circ f$ is + proper.  For
proper apply this result to $(f \circ g)^{-1} = g^{-1} \circ
f^{-1}$ as well.

(f):  It is clear that a continuous function is a $+$proper relation. The converse follows from (c).

(g): If $y$ is a point of the intersection on the left of (\ref{eq1.4}) then
$\{ g^{-1}(y) \cap A: g \in \mathcal{F}, A \in \mathcal{A}  \}$ is
a filterbase of compacta and so has a nonempty intersection. For
$x $ in this intersection  $x \in C$ and  $(x,y) \in g$ for all $g \in \mathcal{F}$.
Hence, $(x,y) \in f$.

(h): If no $g(A)$ is contained in $U$ and $X \setminus U$ is
compact then  \\ $\{ g(A) \cap (X \setminus U): g \in \mathcal{F}, A
\in \mathcal{A} \}$ is a filterbase of compacta and so has a
nonempty intersection. By (g) $f(C)$ meets $X \setminus U$.

(i): Because the relations in $\mathcal{F}$ are + proper and the sets in $\A$ are compact,
 $\{ g(A): g \in \mathcal{F}, A \in \mathcal{A} \}$ is a filterbase of
compacta with intersection $f(C)$ by (g).  So it is eventually contained in the open set $U$.

$\Box$ \vspace{.5cm}

\begin{ex}\label{ex1.1new} The composition of closed relations need not be closed.\end{ex}
Let ${\mathbb N}_* = {\mathbb N} \cup \{\infty \}$ be the one point compactification of ${\mathbb N}$ and let
$\hat Y = \{ -1, 0, 1 \} \times {\mathbb N}_*  $ and $Y = \hat Y \setminus \{(0,\infty)\}$. Define
the  function $g $ on $ Y$ by
\begin{equation}\label{eq1.30}
\begin{split}
(-1,\infty) \mapsto (1,\infty), (1,\infty) \mapsto (-1,\infty), \qquad \mbox{and} \\
(-1,n) \mapsto (0,n), \quad (0,n) \mapsto (1,n), \quad (1,n)
\mapsto (-1,n)
\end{split}
\end{equation}
for all $ n \in {\mathbb N} $. The bijective map $g$ a closed relation but is not continuous.  Notice that
$g^n = 1_{Y_0}$ for $n \equiv 0$ mod 6 and $g^n = g$ for $n \equiv 1$ mod 6.  If $n \equiv 2,3,4,5$ then
the map $g^n$ is not a closed relation.
%Since $g$ is not continuous and so is not the restriction of
%a continuous map, it cannot have a map compactification. For any compactification
%$(\tilde Y, \tilde g)$ of $(Y,g)$, the compact
%$\mathcal{G}g \cap \mathcal{G}g^{-1} $ class $E = \{(-1,\infty), (1,\infty) \}$ is properly contained in
%its $\mathcal{G}\tilde g \cap \mathcal{G}\tilde g^{-1} $ class $\tilde E$.

$\Box$ \vspace{.5cm}

Thus, the composition of closed relations on a
noncompact space need not be closed, but by (e) the composition of +proper
relations is  a +proper  relation and so is closed. If the space is compact then every
closed relation is proper by (a).

If $X = Y$ then a relation $f : X \to X$ is
called a relation on $X$ and $f^i$ is
the $i$-fold composition when $i$ is a positive integer, the
identity map $1_X$ when $i = 0$ and is $(f^{-1})^{|i|}$ when $i$
is negative. A subset $A \subset X$ is $f \ +$\emph{invariant}
when $f(A) \subset A$ and is $f$ \emph{invariant} when $f(A) = A$.
On the other hand, if $g$ is a map on $X$ then we call the
relation $f \ g +$invariant(or $g$ invariant) if  $f \subset X
\times X$ is $+$invariant(resp. invariant) with respect to the map
$g \times g$ on $ X \times X$, that is, when $(g \times g)(f)
\subset f$ (resp. $(g \times g)(f) = f$).

We call a pair $(X,f)$ a \emph{dynamical system} \index{dynamical system}when $f$ is a closed
relation on a space $X$.  It is a \emph{compact dynamical system} \index{dynamical system!compact dynamical system}
when $X$ is compact.

Our primary interest is in dynamical systems with
$f$  a continuous map. However, it is useful to
consider as far as possible the more general case that $f$ is a
closed relation on $X$.  We follow the notation in Akin (1993). There the
spaces were assumed metrizable but the results which we will quote
do not require metrizability as a glance at their proofs will
show. As in Akin (1993) we define certain extensions of $f$.

The \emph{orbit relation} $\mathcal{O}f $ \index{$\mathcal{O}f$}is defined to be
$\bigcup_{i=1}^{\infty} f^i$ with closure in $X \times X$ denoted
$\mathcal{N}f$. \index{$\mathcal{N}f$} Thus, $\mathcal{O}f$  is transitive but need not be
closed while $\mathcal{N}f$ is closed but need not be transitive.
Define $\mathcal{G}f$ \index{$\mathcal{G}f$} to be the smallest closed, transitive
relation which contains $f$. The relation $f$ is both closed and
transitive exactly when $f = \mathcal{N}f$ in which case $f =
\mathcal{G}f$. Consequently, we can obtain $\mathcal{G}f$ ``from
above" as the intersection of all closed, transitive relations
which (like $X \times X$) contain $f$. Alternatively, we can
proceed ``from below" by using transfinite induction.

With $\mathcal{N}_0 f = f $ we define, inductively for ordinals
$\alpha$, \index{$\mathcal{N}_{\alpha}f$}

\begin{equation}\label{eq1.6}
\begin{split}
\mathcal{N}_{\alpha + 1}f \quad = \quad
\mathcal{N}(\mathcal{N}_{\alpha}f) \qquad \mbox{and} \hspace{2cm} \\
\mathcal{N}_{\beta}f \quad = \quad \overline{ \bigcup_{\alpha <
\beta} \mathcal{N}_{\alpha}f } \qquad \mbox{ for } \ \beta \
\mbox{a limit ordinal}.
\end{split}
\end{equation}
When $X$ is metrizable the subsets of $X \times X$ are all
separable, and this process stabilizes at some countable ordinal
to obtain $\mathcal{G}f$. In any case, it stabilizes at some
ordinal with cardinality at most that of $X \times X$.

\begin{prop}\label{prop1.3}  Let $f$ be a closed relation on a space $X$.
If $f$ is $+$proper then
\begin{equation}\label{eq1.7}
\G f \quad = \quad f \cup \G f \circ f,
\end{equation}
and if $f$ is proper then
\begin{equation}\label{eq1.8}
\G f \quad = \quad f \cup f \circ \G f.
\end{equation}
\end{prop}

{\bfseries Proof:}  Since $\G f$ is transitive and contains $f$ it
contains the right hand side of each equation.  Furthermore, each of
the right hand sides defines a relation on $X$ which contains $f$ and is easily seen to be transitive.
By Proposition \ref{prop1.2}(d) $\G f \circ f$ is closed if $f$ is $+$proper and
$f \circ \G f$ is closed if if $f$ is proper.

$\Box$ \vspace{.5cm}

For each $x \in X$ we define $\mathcal{R}f(x)$ to be the closure
of the orbit $\mathcal{O}f(x)$.  While each $\mathcal{R}f(x)$ \index{$\mathcal{R}f$}is
closed, the \emph{orbit-closure relation} $\mathcal{R}f$ is usually
a proper subset of $\mathcal{N}f$ -and so is
 not a closed relation- even when $X$ is compact and $f$ is a map.

When $f$ is a +proper  relation, e.g. a continuous map, the iterates $f^i$ are
closed for positive integers $i$ and the limit relations are of interest.
The \emph{omega limit point set} of the orbit of $x \in X$ is \index{$\omega f$}
$\omega f(x) = lim sup \{ f^i(x) \}$ which defines the relation
$\omega f$. Recall that for a sequence $\{ A_i \} \ lim sup \{ A_i
\} = \bigcap_i \overline{\bigcup_{j \geq i } \{ A_j \}}$ so that
$\overline{ \bigcup_i \{A_i\}} = \bigcup_i \{A_i \} \cup lim sup
\{ A_i \}$ when all of the sets $A_i$ are closed. Thus, if $f$ is
a +proper  relation then for each
point $x$, the closure of the orbit $\mathcal{O}f(x)$ is $\mathcal{R}f(x) = \mathcal{O}f(x) \cup \omega f(x)$.
If we define $\Omega f = lim sup \{ f^i \}$  \index{$\Omega f$}then
$\mathcal{N} f = \mathcal{O}f \cup \Omega f $ when $f$ is + proper.

For any relation $f$ on $X$ the \emph{cyclic set} \index{cyclic set $|f|$}$|f| =_{def} \{ x : x \in f(x) \}$. If $f$ is a closed
relation then $|f|$ is a closed set. Clearly,
$|f| = |f \cap f^{-1}|$.

Motivated by the case when $f$ is a map, we say that $|f|$ is the
set of \emph{fixed points}, $|\mathcal{O}f|$ is the set of
\emph{periodic points}\index{periodic point}, $|\mathcal{R}f| $ is the set
of \emph{recurrent points} \index{recurrent point} and $|\mathcal{N}f|$ is the set of
\emph{nonwandering points} \index{nonwandering point}. Following Auslander (1964) we call
$|\mathcal{G}f|$ the set of \emph{generalized recurrent points}. \index{recurrent point!generalized recurrent point}

The relation $ \mathcal{G}f \cap \mathcal{G}f^{-1}$ is symmetric
as well as transitive and so restricts to an equivalence relation
on $|\mathcal{G}f | $. We allow ourselves the ambiguous notation
because $\mathcal{G}(f^{-1}) = (\mathcal{G}f)^{-1}$.

Notice that $1_X \cup [ \mathcal{G}f \cap \mathcal{G}f^{-1}]$ is
an equivalence relation on the entire space $X$. Its equivalence
classes are those of $ \mathcal{G}f \cap \mathcal{G}f^{-1}$ in
$|\mathcal{G}f |$ and the singleton sets of points of $X \setminus | \mathcal{G}f |$.

Our goal is to compactify a closed relation on
a locally compact space in such a way that no new recurrence is
introduced. However, we will first prove some compact space results.

What we will need is an analogue of Theorem 4.5 of Akin
(1993), extending the chain result - with some weakening - to get
a theorem for generalized prolongations.

We begin by recalling the chain relation.  This is necessarily
a uniform space construction.  In Section \ref{secchain} we will consider the
chain relation on noncompact spaces, but here we review the well-known
 compact space construction from Conley (1978).

A compact space admits a unique
uniformity consisting of all neighborhoods of the diagonal. Again
the results which we quote from Akin (1993) are stated for compact
metric spaces but they are true with the same proofs for arbitrary
compact spaces.

If $f$ is a closed relation on a compact  space then the
\emph{chain relation} \index{$\mathcal{C}f$} \index{chain relation}
$\mathcal{C}f =_{def} \bigcap_{V}
\mathcal{O}(V \circ f)$ where $V$ varies over $\U = \mathcal{U}_X$, the
set of all neighborhoods in $X \times X$ of the diagonal $1_X$.
That is, $(x,y) \in \mathcal{C}f$ if for every $V \in \mathcal{U}$
there exists a finite sequence $x_0,y_0,....y_{n-1},x_n \in X$
with $n \geq 1$ such that $x = x_0, y = x_n \ (x_i,y_i) \in f$ for
$i = 0,...,n-1$ and $(y_{i-1},x_i) \in V$ for $i = 1,...,n$. Such
a sequence is called a $V$ chain for $f$ from $x$ to $y$. Of
course, when $f$ is a map, $y_i = f(x_i)$ and so we can think of a
$V$ chain for a map $f$ as a sequence $x_0,...,x_n$ with
$(f(x_{i-1}),x_i) \in V$ for $ i = 1,...,n$. When $X$ is
metrizable with metric $d$ and $V = \bar{V}_{\epsilon} =_{def} \{
(x,y) : d(x,y) \leq \epsilon \}$ for some positive $\epsilon$ then
we refer to a $V$ chain as an $\epsilon$ chain.

In general, we will call a sequence $x_0,...,x_n$ with $(x_i,x_i+1) \in f$ for
$i = 0,...,n$ a $0$ chain for $f$.  So a $0$ chain is a piece of an orbit sequence.

The relation $\mathcal{C}f$ is closed, transitive and contains $f$
so that $\mathcal{G}f \subset \mathcal{C}f$.  The inclusion may be
proper. For example, if $f$ is the identity map on a connected,
compact space $X$ then $\mathcal{C}f = X \times X$ but
$\mathcal{G}f = 1_X = f$.

The points of $|\mathcal{C}f|$ are called \emph{chain recurrent
points} \index{recurrent point!chain recurrent point}for $f$.  On $|\mathcal{C}f|$ the relation $\mathcal{C}f
\cap \mathcal{C}f^{-1}$ is a closed equivalence relation.

If $D$ is any closed subset of $X$ and $f$ is any closed relation
on $X$ then $f_D =_{def} f \cap (D \times D)$ is a closed relation
on $D$ called the \emph{restriction} \index{restriction $f_D$} of $f$ to $D$. Notice that no
invariance is assumed. Note too that because $D$ is compact the
neighborhoods $\mathcal{U}_D$ of $1_D$ in $D \times D$ are exactly the
restrictions to $D$ of the elements of $\mathcal{U}_X$.

We will use the notation $A \subset \subset B$ \index{$\subset \subset$} to mean the closure of
$A$ is contained in the interior of $B$.

For a dynamical system $(X,f)$  a subset $C \subset X$ is called
$ f$ \emph{unrevisited} \index{unrevisited set} when $\mathcal{O} f(C) \cap \mathcal{O} f^{-1}(C) \subset C$. Thus,
$A = (1_X \cup \mathcal{O} f)(C)$ is $ f$ +invariant, $B = (1_X \cup \mathcal{O} f^{-1})(C)$ is $ f^{-1}$ +invariant
and $C = A \cap B$.  For example, if $C$ is $ f$ +invariant or $ f^{-1}$ +invariant then it is $ f$ unrevisited.
Intuitively, $C$ is unrevisited when no orbit sequence leaves $C$ and then returns to it. If $x_0,....,x_n$ is
a $0$ chain for $f$ with $x_0, x_n \in C$ then $C$ unrevisited implies $x_1,...,x_{n-1} \in C$ as well.

\begin{theo}\label{th2.1}  Assume  $f$ is a closed relation on a compact
space $X$. Let $C$ be a $\G f$ unrevisited closed subset of $X$ so that
\begin{equation}\label{eq2.1}
\mathcal{G}f (C) \cap \mathcal{G}f^{-1}(C) \quad \subset \quad C.
\end{equation}

If $D$ is any closed neighborhood of $\ C$, i.e. $D$ is closed and

\begin{equation}\label{eq2.2}
C \quad \subset \subset \quad D, \hspace{1cm}
\end{equation}
then:
\begin{equation}\label{eq2.3}
(\mathcal{G}f)_C \quad \subset \quad \mathcal{G}(f_D).
\end{equation}
That is, if $x,y \in C$ with $y \in \mathcal{G}f (x)$ then $(x,y)$
is in the smallest closed transitive relation which contains the restriction
$f_D$.

Furthermore,
\begin{equation}\label{eq2.4}
(\mathcal{G}f)_C \quad \subset \quad \mathcal{C}(f_C).
\end{equation}
That is, if $x,y \in C$ with $y \in \mathcal{G}f (x)$ then $(x,y)$
is in the chain relation for the restriction $f_C$, i.e. for every
$V \in \mathcal{U}_X$ there is a $V$ chain $x_0,y_0,...,x_n$ for
$f$ from $x$ to $y$ with $x_i, y_i \in C$ for all $i$.
\end{theo}

{\bfseries Proof:} We use the transfinite construction (\ref{eq1.6}) for
$\G f$.

We prove, by transfinite induction, the result
with $(\mathcal{G}f)_C$ replaced by $(\mathcal{N}_{\alpha}f)_C$.  To be precise,
we prove by induction on $\alpha$ that if $ C$ is any closed, $\G f$ unrevisited closed set
contained in the interior of $D$ then $(\mathcal{N}_{\alpha}f)_C \quad \subset \quad \mathcal{G}(f_D)$.

For $\alpha = 0$ the result is clear because $(\mathcal{N}_{0}f)_C = f_C \subset f_D$. This is the
initial step of the induction.

For any closed  $V \in \mathcal{U}$ we define
\begin{equation}\label{eq2.5}
\begin{split}
A_1 \quad = \quad (1_X \cup \mathcal{G}f )(V(C)),\hspace{.5cm}
\\
B_1 \quad = \quad (1_X \cup \mathcal{G}f^{-1}
)(V(C)),\\
C_1 \quad = \quad A_1 \cap B_1. \hspace{2cm}
\end{split}
\end{equation}
Clearly, $C \subset \subset C_1$ and by compactness we can choose
$V$ small enough that $C_1 \subset \subset D$, see Proposition \ref{prop1.2}. In proving the
inductive step for  $C$ we may apply the inductive
hypothesis to $C_1$ as it is a closed, $\G f$ unrevisited set which is contained in the interior of $D$.

If $\beta$ is a limit ordinal and $(x,y) \in \mathcal{N}_{\beta}f$
with $x,y \in C$ then, because $C \subset \subset C_1, \ (x,y)$ is
the limit of a net of pairs $(x_i,y_i) \in
(\mathcal{N}_{\alpha_i}f) \cap (C_1 \times C_1)$ with $\alpha_i <
\beta$. By induction hypothesis these points are in the closed
relation $\mathcal{G}(f_D)$ and so their limit is as well.

If $(x,y) \in \mathcal{N}_{\alpha + 1}f $ with $x,y \in C$ then it
is the limit of a net of pairs $(x_i,y_i) \in
(\mathcal{N}_{\alpha}f)^{n_i} \cap (C_1 \times C_1)$. So for each
$i$ in the directed set indexing the net there is a sequence
$a_0,...a_{n_i}$ with $a_0 = x_i, a_{n_i} = y_i$ and
$(a_k,a_{k+1}) \in \mathcal{N}_{\alpha}f $ for $k = 0,...,n_i -
1$.

Since $a_0 = x_i \in C_1 \subset A_1$ it follows from
$\mathcal{G}f$ invariance of $A_1$ that $a_k \in A_1$ for all $k$.
Similarly, $a_{n_i} \in C_1 \subset B_1$ implies that $a_k \in
B_1$ for all $k$.  Hence, for all $k \ (a_k,a_{k+1}) \in
(\mathcal{N}_{\alpha}f)_{C_1} $.  By induction hypothesis this is
contained in the relation $\mathcal{G}(f_D)$ transitivity of which
implies $(x_i,y_i) \in \mathcal{G}(f_D)$.  Again the limit point
$(x,y)$ is in $\mathcal{G}(f_D)$ because the relation is closed.

The induction completed, we choose a filter $\mathcal{F}$ of
closed neighborhoods $D$ of $C$ with intersection $C$.  Applying
(\ref{eq2.3}) we have
\begin{equation}\label{eq2.6}
\bigcap_{D \in \mathcal{F}} \mathcal{C}(f_{D}) \quad \supset \quad
\bigcap_{D \in \mathcal{F}} \mathcal{G}(f_{D}) \quad \supset \quad
(\mathcal{G}f)_C.
\end{equation}
Inclusion (\ref{eq2.4}) now follows because the operation $\mathcal{C}$ is
a monotone, upper semicontinuous map on closed relations, Akin
(1993) Theorem 7.23, which implies that the leftmost intersection
in (\ref{eq2.6}) is
\begin{equation}\label{eq2.7}
  \mathcal{C}(\bigcap_{D \in \mathcal{F}} f_{D}) \quad = \quad \mathcal{C}(f_C).
\end{equation}

$\Box$ \vspace{.5cm}

While we will apply this theorem in a more general context it
has an interesting application in the case when $f$ is a continuous map.
First, we recall one more construction from Akin (1993).

For any closed relation $f$ on a compact space
$X \ \{ (\mathcal{G}f)^n \}$ is a decreasing sequence of closed,
transitive relations so that $\Omega \mathcal{G}f =  \bigcap_n
(\mathcal{G}f)^n $ \index{$\Omega \mathcal{G}f$} is a closed transitive relation. From Akin
(1993) Proposition 2.4 we have
\begin{equation}\label{eq2.8}
\mathcal{G}f \quad = \quad \mathcal{O}f \cup \Omega \mathcal{G}f
\hspace{2cm}
\end{equation}
\begin{equation}\label{eq2.9}
\begin{split}
f \circ \Omega \mathcal{G}f \quad = \quad \mathcal{G}f \circ
\Omega \mathcal{G}f \quad = \quad \Omega \mathcal{G}f  \hspace{2cm}\\ =
\quad \Omega \mathcal{G}f \circ  \mathcal{G}f \quad = \quad \Omega
\mathcal{G}f \circ f. \hspace{1cm}
\end{split}
\end{equation}

Clearly, if $x \in |\mathcal{G}f|$, i.e. $(x,x) \in \mathcal{G}f$
then $(x,x) \in \Omega \mathcal{G}f$ and so $y \in f(x)$ implies
$(x,y) \in f \circ \Omega \mathcal{G}f = \Omega \mathcal{G}f$.
Proceeding by induction on $n$ we get $f^n(x) \subset \Omega
\mathcal{G}f (x)$. On the other hand, if $x \in f(z)$ then $(z,x)
\in (\Omega \mathcal{G}f) \circ f = \Omega \mathcal{G}f$ and so,
inductively, $f^{-n}(x) \subset (\Omega \mathcal{G}f)^{-1}(x)$.
Thus, from (\ref{eq2.8}) we see
\begin{equation}\label{eq2.10}
\begin{split}
|\mathcal{G}f| \quad = \quad |\Omega \mathcal{G}f| \hspace{4cm}\\
x \in |\mathcal{G}f| \qquad \Longrightarrow \hspace{4cm} \\
\mathcal{G}f(x) \  = \ \Omega \mathcal{G}f (x)\qquad \mbox{and}
\qquad  \mathcal{G}f^{-1}(x) \  = \ (\Omega \mathcal{G}f)^{-1}
(x).
\end{split}
\end{equation}
\vspace{.5cm}

\begin{lem}\label{lem2.2} Let $f$ be a continuous map on a compact  space
$X$ and let $x,y \in X$. The compact set $\mathcal{G}f (x)$ is
closed and $\mathcal{G}f \ +$invariant. If $x \in |\mathcal{G}f|$
then $\mathcal{G}f (x)$ is $f$ invariant. The compact set
$\mathcal{G}f^{-1}(y)$ is $\mathcal{G}f^{-1} \ +$invariant. If $y
\in |\mathcal{G}f|$  then $\mathcal{G}f^{-1}(y)$ is $f \
+$invariant and is $f$ invariant when $f$ is a surjective map. If
$x,y \in |\mathcal{G}f|$ then $\mathcal{G}f(x) \cap
\mathcal{G}f^{-1}(y)$ is $f$ invariant.
\end{lem}

{\bfseries Proof:} The $+$invariance results just follow from
transitivity of $\mathcal{G}f$.

If $x,y$  lie in $|\mathcal{G}f|$ then by (\ref{eq2.10}) we can use
$\Omega \mathcal{G}f $ instead of $\mathcal{G}f$. From (\ref{eq2.9}) we
get $f$ invariance of $\Omega \mathcal{G}f(x)$ and from the
equation on the reverse relations we get
\begin{equation}\label{eq2.11}
f^{-1}((\Omega \mathcal{G}f)^{-1}(y)) \quad = \quad (\Omega
\mathcal{G}f)^{-1}(y).
\end{equation}
Apply the map $f$ to both sides and recall that $f \circ f^{-1}
\subset 1_X$ with equality when $f$ is surjective.  In any case,
$f(A \cap f^{-1}(B)) = f(A) \cap B$ which equals $A \cap B$ when
$A$ is invariant.

$\Box$ \vspace{.5cm}

If $\tilde{E}$ is the equivalence class $\mathcal{C}f(x) \cap
\mathcal{C}f^{-1}(x)$ then by Akin (1993) Theorem 4.5 the chain
relation of the restriction $f_{\tilde{E}}$ is exactly the
restriction $(\mathcal{C}f)_{\tilde{E}} = \tilde{E} \times
\tilde{E}$ which says that  on $\tilde{E}$ the map $f$ is
\emph{chain transitive}. The analogous result $\mathcal{G}f$ is
false. \vspace{.5cm}

\begin{ex}\label{ex2.4}For a $\G f \cap \G f^{-1}$ equivalence class $E$ it
need not be true that all points in $E$ are $\G (f|E) \cap \G (f|E)^{-1}$ equivalent.
\end{ex}
 Let $X$ be the unit disc in the plane
${\mathbb R}^2$.  Let $f$ on $X$ be the time-one map of the flow
associated with the differential equations -in polar coordinates-
$\frac{dr}{dt} = r \cdot (1 - r), \ \frac{d \theta}{dt} = r \cdot
(1 - r)$. The set of chain recurrent points equals the set of
fixed points: the unit circle $E$ and the origin. $E$ is a single
$\mathcal{G}f \cap \mathcal{G}f^{-1}$ equivalence class. So the
restriction $(\mathcal{G}f)_E = E \times E$. On the other hand,
$\mathcal{G}(f_E) = 1_E$. Notice that the restriction $f_E = 1_E$
is chain transitive on $E$.  This is the general result.

If we use $\frac{dr}{dt} = r \cdot (1 - r), \ \frac{d \theta}{dt}
= sin^2(\frac{\theta}{2}) + r \cdot (1 - r)$ and remove from $X$
the points of the unit circle $E$ on and above the $x$ axis to obtain the
locally compact space $X_1$ then the new time-one map  restricts
to define a map $f_1$ on $X_1$. $E \cap X_1$ is a noncompact
$\mathcal{G}f_1 \cap \mathcal{G}f_{1}^{-1}$ equivalence class
which is $f_1 \ +$ invariant but not invariant.

$\Box$ \vspace{.5cm}

\begin{cor}\label{cor2.5} Let $f$ be a continuous map on a compact  space
$X$. Let $x$ be a generalized recurrent point, i.e. $x \in
|\mathcal{G}f|$, with $E$ the associated $\mathcal{G}f \cap
\mathcal{G}f^{-1}$ equivalence class, i.e. $E = \mathcal{G}f(x)
\cap \mathcal{G}f^{-1}(x)$. $E$ is a closed $f$ invariant set and
the restriction of $f$ to $E$ is a chain transitive map, that is,
\begin{equation}\label{eq2.12}
\mathcal{C}(f_E) \quad = \quad E \times E.
\end{equation}
\end{cor}

{\bfseries Proof:} By (\ref{eq2.4}) of Theorem \ref{th2.1} again, we have that
$\mathcal{C}(f_E) \supset (\mathcal{G}f)_E $ and the latter set is
$E \times E$ because $E$ is a $\mathcal{G}f \cap
\mathcal{G}f^{-1}$ equivalence class.

$\Box$ \vspace{.5cm}

Now we are ready to compactify.

Let $\mathcal{B}(X)$ \index{$\mathcal{B}(X)$} denote the Banach algebra of bounded,
continuous real-valued functions on $X$ equipped with the sup norm.
If $h : X_1 \to X_2$ is a
continuous map, then $h^* : \mathcal{B}(X_2) \to \mathcal{B}(X_1)$
is the continuous algebra homomorphism defined by $h^*(u) = u
\circ h$.\index{$h^*$}

A \emph{compactification} \index{compactification} of $X$ is a continuous map $j : X \to
\hat{X}$ of $X$ onto a dense subset of compact Hausdorff space
$\hat{X}$. Since $j(X)$ is dense in $\hat{X}$ the algebra
homomorphism $j^* : \mathcal{B}(\hat{X}) \to \mathcal{B}(X)$ is an
isometry with image $j^*(\mathcal{B}(\hat{X}))$ a closed
subalgebra of $\mathcal{B}(X)$. Conversely, for any closed
subalgebra $\mathcal{A}$ of $\mathcal{B}(X)$ the Gelfand
construction yields a compactification, unique up to
homeomorphism, such that $\mathcal{A} =
j^*(\mathcal{B}(\hat{X}))$.

The space $\hat{X}$ is metrizable iff the associated algebra
$\mathcal{A}$, or equivalently $\mathcal{B}(\hat{X})$, is
countably generated. When $\{ j_n : X \to [0,1] \}$ is a sequence
which generates $\mathcal{A}$ then we can define the map $j$ from
$X$ into the product of countably many copies of the unit interval
and the associated compactification is just the closure of $j(X)$
in the product.

The map $j$ is injective iff the subalgebra $\mathcal{A}$
distinguishes the points of $X$. The map $j$ is an embedding, i.e.
a homeomorphism onto its image $j(X)$ equipped with the subspace
topology, iff the subalgebra $\mathcal{A}$ of $\mathcal{B}(X)$
distinguishes points and closed sets. We will call $j : X \to
\hat{X}$ a \emph{proper compactification} \index{compactification!proper compactification}
when $j$ is an
embedding. Since $X$ is locally compact, a compactification $j$ is
proper iff its image $j(X)$ is open as well as dense in  $\hat{X}$
and, in addition, the continuous map $j : X \to j(X)$ is an injective, proper map.

When $j$ is a proper compactification, we will usually regard $j$
as an inclusion, identifying $X$ with its image, and write $X
\subset \hat{X}$. Then $u \in \mathcal{B}(X)$ has a -necessarily
unique- extension $\hat{u} \in \mathcal{B}(\hat{X})$ iff $u$ lies
in the subalgebra $\mathcal{A}$ associated with the
compactification $X \subset \hat{X}$.

 The functions with compact support
generate a closed subalgebra $\mathcal{A}_0$ which distinguishes
points and closed sets. When $X$ is metrizable this algebra is
countably generated. The associated proper compactification is the
\emph{one point compactification} \index{compactification!one point compactification} which we will denote $X_*
\supset X$.

If $(X,f)$ is a dynamical system and $j : X \to \hat{X}$ is a
compactification then we define $(\X,\hat{f})$ to be the compact dynamical system with $\hat f$
the relation on $\hat{X}$ which is the closure of $(j \times j)(f)
\subset \hat{X} \times \hat{X}$. \index{dynamical system!compactification} When  $X \subset \hat{X}$ is a
proper compactification,

\begin{equation}\label{eq1.9}
\hat{f} \cap (X \times X) \quad = \quad f
\end{equation}
because $X \times X$ has the relative topology from $\hat{X}
\times \hat{X}$. In that case, we will call $(\X, \f)$ a \emph{proper compactification}
of $(X,f)$.\index{dynamical system!proper compactification}

If $f$ is a continuous map on $X$ then $\hat{f}$ is a continuous
map on $\hat{X}$ iff the associated algebra $\mathcal{A} \subset
\mathcal{B}(X)$ is $f^* \ +$invariant, i.e.  $f^*(\mathcal{A})
\subset \mathcal{A}$. If $f$ is a homeomorphism on $X$ then
$\hat{f}$ is a homeomorphism on $\hat{X}$ iff $\mathcal{A}$ is
$f^*$ invariant, i.e. $f^*(\mathcal{A}) = \mathcal{A}$.

For an exposition with proofs of these compactification results
see, e.g. Akin (1997) Chapter 5.

\begin{lem}\label{lem1.4} Let  $\hat{X} \supset X$ be a proper compactification
 of $(X,f)$.
If $(x,y) \in (X \times X) \cap \mathcal{G} \hat{f}$ then
either $(x,y) \in \mathcal{G}f$ or there exists $z \in \hat{X} \setminus X$
such that $(x,z), (z,y) \in \mathcal{G} \hat{f}$.
\end{lem}

{\bfseries Proof:}  Let $g_1 = \{ (x,y) \in \mathcal{G}\hat f : y \in \hat X \setminus X \}$ and
$g_2 = \{ (x,y) \in \mathcal{G}\hat f : x \in \hat X \setminus X \}$.
Observe that $g_1 \circ \mathcal{G}\hat f \subset g_1,  \ \mathcal{G}\hat f \circ g_2 \subset g_2$ and
$g_2 \circ g_1 = \mathcal{G}\hat f \circ g_1 = g_2 \circ \mathcal{G}\hat f$. It follows that
$F = \mathcal{G} f \cup g_1 \cup g_2 \cup g_2 \circ g_1$ is a closed
transitive relation which is contained in $\mathcal{G}\hat f$. Since $F$ is closed and contains $f$
 it contains $\hat f$ and so
equals $\mathcal{G}\hat f$.  Clearly, if $(x,y) \in (X \times X) \cap F$ then either $(x,y) \in \mathcal{G} f$
or $(x,y) \in g_2 \circ g_1$.

$\Box$ \vspace{.5cm}

\begin{df} \label{def1.5} Let  $(\hat{X}, \hat f)$  be a  compactification
 of a dynamical system $(X,f)$.We
say that the compactification $X \subset \hat{X}$ is \emph{dynamic
for $f$} or that $(\hat{X},\hat{f})$ is a \emph{dynamic
compactification of $(X,f)$} \index{dynamic compactification} \index{dynamical system!dynamic compactification}
when the compactification is proper and

\begin{equation}\label{eq1.10}
 (\mathcal{G} \hat{f}) \cap (X \times X) \quad =
\quad  \mathcal{G} f.
\end{equation}

We say that  $(\hat{X},\hat{f})$ is an \emph{almost dynamic
compactification of $(X,f)$} \index{almost dynamic compactification}%
 \index{dynamical system!almost dynamic compactification}%
when the compactification is proper and
\begin{equation}\label{eq1.11}
[1_{\hat{X}} \cup \mathcal{G} \hat{f}] \cap (X \times X) \quad =
\quad 1_X \cup \mathcal{G} f.
\end{equation}
\end{df}

$\Box$ \vspace{.5cm}

Thus, a dynamic compactification is  one where the extended relation
$\hat{f}$ on $\hat{X}$ introduces no new generalized recurrence in $X$ or, more
generally, any new $\G \f$ relations between points of $X$.  The ``almost dynamic compactification"
concept is introduced as a tool because it is easier to obtain and because the gap between the two concepts is easy
to investigate.

If $\G f$ is reflexive, i.e. $1_X \subset \G f$ then an almost dynamic compactification is the same as a
dynamic compactification.

Clearly, if $(\X, \f)$ is a proper compactification of $(X,f)$ then $(\X, \f^{-1})$ is a proper compactification
of $(X,f)$. Furthermore, $(\X, \f^{-1})$ is dynamic or almost dynamic for $(X,f^{-1})$ if $(\X, \f)$ satisfies the
corresponding property for $(X,f)$.

\begin{theo}\label{th1.6} Assume that  $ (\hat{X},\hat{f})$ is an almost dynamic
compactification of  $(X,f)$.

\begin{enumerate}
\item[(a)]For $\hat{E} \subset |\mathcal{G} \hat{f}| $ a $\mathcal{G}
\hat{f} \cap \mathcal{G} \hat{f}^{-1}$ equivalence class  let $E =
\hat{E} \cap X$. Exactly one of the following four possibilities
holds:
\begin{itemize}
\item[(i)]  $\hat{E}$ is a compact subset of $\hat{X} \setminus
X$ and $E = \emptyset$.
\item[(ii)] $E$ is contained in $ |\mathcal{G}f|$  and is a noncompact
$\mathcal{G} f \cap \mathcal{G}f^{-1}$ equivalence class whose
$\hat{X}$  closure meets $\hat{X}\setminus X$ and is contained in
$\hat{E} \ $.
\item[(iii)] $\hat{E} = E$ is contained in $|\mathcal{G}f|$ and is a compact
$\mathcal{G} f \cap \mathcal{G} f^{-1} $ equivalence class.
\item[(iv)]  $E$ is a compact $1_X \cup (\mathcal{G}f \cap \mathcal{G}f^{-1})$
equivalence class and the compact set $\hat{E} \setminus E \subset
\hat{X} \setminus X$ is nonempty.  $\hat{E}$ is the disjoint union
of the nonempty compacta $E$ and   $\hat{E} \setminus E$.
\end{itemize}

\item[(b)]  If $\ (x,y) \in  (\mathcal{G} \hat{f}) \cap (X \times X)$ but
$(x,y) \not\in \mathcal{G}f$ then $y = x$ and so $x \in
|\mathcal{G} \hat{f}|$. Furthermore, if $\hat{E}$ is the
$\mathcal{G} \hat{f} \cap \mathcal{G} \hat{f}^{-1}$ equivalence
class of $x$ then $E = \hat{E} \cap X$ is the singleton set $\{ x
\}$ and $\hat{E} \setminus E$ is nonempty. In particular, case
(iv) of (a) applies to $\hat{E}$.

\item[(c)]If $x,y \in |\mathcal{G}f|$  lie in distinct $\mathcal{G}f
\cap \mathcal{G}f^{-1}$ equivalence classes then their equivalence
classes have disjoint closures in $\hat{X}$.
\end{enumerate}
\end{theo}

{\bfseries Proof:} (a): From (\ref{eq1.11}) it follows that  $E = \hat{E}
\cap X$ is a single $1_X \cup (\mathcal{G}f \cap
\mathcal{G}f^{-1})$ equivalence class or else $E = \emptyset$. In
particular, if $E$ contains more than one point then it is a
$\mathcal{G} f \cap \mathcal{G}f^{-1}$ class and so is contained
in $|\mathcal{G}f| $.

The class $\hat{E}$ is a closed subset of $\hat{X}$ and so is
compact.

It follows that if $E$ is empty then case (i) holds.  Now assume
$E$ is nonempty.

If $E$ is noncompact then it certainly contains more than one
point and is a proper subset of its $\X$ closure which is a compact.
Since $\hat{E}$ is closed, it contains the closure of $E$. This is case (ii).

If $\hat{E} \setminus X$ is empty then $E = \hat{E}$ is compact.
If $E$ contains more than one point then it is a $\mathcal{G}f
\cap \mathcal{G}f^{-1}$ class and the case (iii) conditions hold.

In considering case (iii) we have to exclude the odd possibility
that $E = \hat{E}$ is a singleton $\{ x \}$, but $x$ is not in
$|\mathcal{G}f|$ and so the $1_X \cup (\mathcal{G}f \cap
\mathcal{G}f^{-1})$ equivalence class $E$ is not a $\mathcal{G}f
\cap \mathcal{G}f^{-1}$ class.  This possibility is ruled out by Lemma \ref{lem1.4} which says
that if $(x,x) \in \G \f  \setminus \G f$ then there exists $z \in \X \setminus X$ such that
$(x,z), (z,x) \in \G \f$ and so $z \in \hat{E}$.

There remains the possibility that $E$ is compact and that the
compact set $\hat{E} \setminus E = \hat{E} \cap (\hat{X} \setminus
X)$ is nonempty. This is case (iv).

(b): From (\ref{eq1.11}) it follows that $(x,y) \in 1_X$ and so $x = y$.
Since $(x,y) = (x,x) \in \mathcal{G} \hat{f} $ it follows that $x
\in |\mathcal{G} \hat{f}|$. Let $\hat{E}$ be its $\mathcal{G}
\hat{f} \cap \mathcal{G} \hat{f}^{-1}$ equivalence class and $E =
\hat{E} \cap X$. As observed in the first paragraph above, if $E$
contains more than one point then it is a $\mathcal{G} f \cap
\mathcal{G} f^{-1}$ equivalence class.  Since $(x,y) = (x,x)
\not\in \mathcal{G}f$, it follows that $E$ is the singleton set
$\{ x \}$. Thus, neither case (i) nor case (ii) applies. By Lemma \ref{lem1.4} again $\hat{E}
\setminus E \not= \emptyset$ and so case (iii)
too is excluded.  That leaves case (iv).

(c): Observe that the closure of the $\mathcal{G}f \cap
\mathcal{G}f^{-1}$ equivalence class of $x$ is contained in its
$\mathcal{G}\hat{f} \cap \mathcal{G}\hat{f}^{-1}$ equivalence
class. Hence, if the closures of the classes of $x$ and $y$ meet
in $\hat{X}$ then $x$ and $y$ are $\mathcal{G}\hat{f} \cap
\mathcal{G}\hat{f}^{-1}$ equivalent. It follows from (\ref{eq1.11}) that
they are $\mathcal{G}f \cap \mathcal{G}f^{-1}$ equivalent.

$\Box$ \vspace{.5cm}

{\bfseries Remark:} In particular, if   $(\hat{X},\hat{f})$ is an almost
dynamic compactification of $(X,f)$ then the the only way the
extension from $f$ to $\hat{f}$ can introduce  additional
generalized recurrence in $X$ is via the peculiar case (iv)
situation described in part (b) above. In detail, the only way
that $x \in X$ can lie in $|\mathcal{G}\hat{f}| \setminus
|\mathcal{G}f|$ is when there is a $\mathcal{G}\hat{f} \cap
\mathcal{G}f^{-1}$ equivalence class $\hat{E}$ which meets
$\hat{X} \setminus X$ but has $\{ x \} = E = \hat{E} \cap X$ disjoint from $|\G f|$.
Thus, this is the only way an almost dynamic compactification can fail to be
a dynamic compactification.

\vspace{.5cm}

\begin{ex}\label{ex1.25} An anti-symmetric relation need not have any
anti-symmetric compactifications and
other Case (iv) examples.\end{ex}
 Let ${\mathbb N}$ be the set of natural
numbers and let $X_0 = \{ -1, 0, 1 \} \times {\mathbb N}$. Define
the continuous function $f_0 $ on $ X_0$ by

\begin{equation}\label{eq1.29}
(-1,n) \mapsto (0,n), \quad (0,n) \mapsto (1,n), \quad (1,n)
\mapsto (1,n + 1)
\end{equation}
for all $ n \in {\mathbb N} $. The orbit relation
$\mathcal{O}f_0$ is closed and so equals
$\mathcal{G}f_0$.

Let $\hat{X}$ be the two-point compactification of $X$ which
adjoins one new point $z_0$ which is the limit point of the
sequence $\{(0,n)\}$, and one new point $z_{\pm}$ which is the common
limit point of both of the sequences $\{(1,n)\}$ and $\{(-1,n)\}$.
The closure $\hat{f}$ is the relation $f_0 \cup \{
(z_{\pm},z_0),(z_0,z_{\pm}),(z_{\pm},z_{\pm})\}$. It is easy to check that
$(\hat{X},\hat{f})$ is a dynamic compactification of $(X_0,f_0)$.
Notice that  $\hat{f}$ is not a
map. $|\mathcal{G}\hat{f}|$ consists of the single $\mathcal{G}
\hat{f} \cap \mathcal{G} \hat{f}^{-1}$ equivalence class $\hat{E}
= \{z_0,z_{\pm}\}$.

Now let $X_1 = \hat{X} \setminus \{z_0\} = X_0 \cup \{ z_{\pm} \} $ and
$f_1 = f_0 \cup \{(z_{\pm},z_{\pm}) \}$. Notice that $f_1$ is a closed
relation on $X_1$ and is a map but not a continuous map.
$\mathcal{G}f_1 = \mathcal{G}f_0 \cup \{ (x,z_{\pm}) : x \in X_1 \}$.
So $|\mathcal{G}f_1|$ consists of the single $\mathcal{G} f_1 \cap
\mathcal{G} f_1^{-1}$ equivalence class  $E_1 = \{z_{\pm} \} = \hat{E}
\cap X_1$. $(\hat{X},\hat{f})$ is a dynamic compactification of
$(X_1,f_1)$ and $\hat{E}$ is an example of the peculiar case (iv).
Observe that $1_{X_1} \cup \mathcal{G}f_1$ is a closed,
anti-symmetric, transitive relation on $X_1$. Any proper
compactification $\tilde{X}$ of the space $X_1$ will map onto
$\hat{X}$ because the latter is the one-point compactification of
$X_1$. If $\tilde{f}$ is the extension of $f_1$ to $\tilde{X}$
then every point of $\tilde{X} \setminus X_1$ is $\mathcal{G}
\tilde{f} \cap \mathcal{G} \tilde{f}^{-1}$ equivalent to $z_{\pm}$.
Hence, $1_{\tilde{f}} \cup \mathcal{G}\tilde{f}$ is not anti-symmetric for any proper
compactification of $X_1$.

Now, instead, let $X_2 = \hat{X} \setminus \{z_{\pm}\} = X_0 \cup \{z_0\}$
and let $f_2 = f_0$ which is a closed relation on $X_2$. It is not
a map because $f_2(z_0)$ is the empty set. $\mathcal{G}f_2 $ is
$\mathcal{G}f_0 $ which is still closed in $X_2 \times X_2$.
Hence, $|\mathcal{G}f_2| = \emptyset$. That is, the transitive relation
$\G f_2$ is asymmetric.  $(\hat{X},\hat{f})$
is an almost dynamic compactification of $(X_2,f_2)$ and $\hat{E}$ is an
example of case (iv).  In this case the singleton set $E = \{z_0
\} = \hat{E} \cap X_2$ is not a
$\mathcal{G} f_2 \cap \mathcal{G}f_2^{-1}$ equivalence class.
This is a case where new generalized
recurrence is introduced as described in the remark above.
That is, it provides an example of an almost dynamic compactification which is not a
dynamic compactification.

$\Box$ \vspace{.5cm}

When $f$ is +proper, e.g. when $f$ is a continuous map, the anomalous case (iv)
 of Theorem \ref{th1.6} does not occur.

 \begin{prop}\label{prop1.8new} For a dynamical system $(X,f)$ the following
 conditions are equivalent.
 \begin{itemize}
 \item[(a)] The relation $f$ is +proper.
 \item[(b)] For every proper compactification $(\X, \f)$ of $(X,f)$, the open
 set $X \subset \X$ is $\f$ +invariant.  That is, $\f(X) \subset X$.
 \item[(c)] There exists a proper compactification $(\X, \f)$ of $(X,f)$ such that the open
 set $X \subset \X$ is $\f$ +invariant.
 \end{itemize}
 \end{prop}

 {\bfseries Proof:}  (a) $\Rightarrow$ (b):  Assume $(x,z) \in \f$ and $x \in X$. There exists
 a net $\{ (x_i,y_i) \in f$ which converges to $(x,z)$ in $\X \times \X$. Let
 $U$ be a bounded open neighborhood of $x$ in $X$ so that $\ol{U} \subset X$ is compact.
 Since $\{ x_i \}$ is eventually in $U$, $\{ y_i \}$ is eventually in $f(U)$.
 If $f$ is +proper then $f(\ol{U})$ is compact in $X$ and so is closed in $\X$. Hence,
$ f(\ol{U})$ is the closure in $\X$ of $\f(U) = f(U)$. Hence, $z \in f(\ol{U}) \subset X$.
Thus, $(x,z) \in \f \cap (X \times X) = f$.

(b) $\Rightarrow$ (c):  Obvious since proper compactifications exist, e.g. the one point compactification.

(c) $\Rightarrow$ (a):  Let $A \subset X$ be compact and so is closed in
$\X$. If $X$ is +invariant for $\f$ then $f(A) = \f(A)$ By Proposition \ref{prop1.2} (a)
$\f(A)$ is closed in $\X$ and so is compact. Because the compactification is proper,
the topology on $X$ is the subspace topology and so $f(A)$ is compact in $X$.

$\Box$ \vspace{.5cm}

\begin{cor}\label{cor1.8new} Let $(\X ,\f)$ be a proper compactification of
the dynamical system $(X,f)$ with $f$ a +proper relation.

For every positive integer $n$, $(X \times \X) \cap \f^n = f^n$.  Furthermore,
\begin{equation}\label{eq3.1}
\begin{split}
(X \times X) \cap (\mathcal{O}\hat{f})   \quad = \quad \mathcal{O}f,
\\
(X \times X) \cap (\omega \hat{f})\quad = \quad
\omega f,
\\(X \times X) \cap (\mathcal{R}\hat{f})   \quad = \quad \mathcal{R}f,
\\
(X \times X) \cap (\mathcal{N}\hat{f})  \quad = \quad
\mathcal{N}f,
\\
\end{split}
\end{equation}

If $ \mathcal{R}f(x)$ (or $\mathcal{N}f(x)$) is a compact subset of $X$, then
 $\mathcal{R}f(x) = \mathcal{R}\hat{f}(x)$  (resp.
$\mathcal{N}f(x) = \mathcal{N}\hat{f}(x)$).
\end{cor}

{\bfseries Proof:} By Proposition \ref{prop1.8new} $X$ is $\f$ +invariant and so $x \in X$
implies $\f^n(x) = f^n(x) \subset X$ and so
\begin{equation}\label{eq3.2new}
\mathcal{O}\f(x) \quad = \quad \mathcal{O} f(x) \subset X.
\end{equation}
The intersection with $X$ of the closure $\mathcal{R} \f(x)$ is the closure in $X$ which is
$\mathcal{R} f(x)$.

Using the lim sup definition we similarly obtain $X \cap \omega \f(x) = \omega f(x)$.

 If $ \mathcal{R} f(x)$ is compact then it is the $\X$ closure of $\mathcal{O} f(x)$
and so equals $\mathcal{R} \f(x)$ by (\ref{eq3.2new}).

Let $x \in X$. If $y \in  \mathcal{N}\f(x)$ then there is a
net $\{ (x_i,y_i) \}$ in $\mathcal{O}\f$ converging to $(x,y)$ and we can assume $x_i \in X$ for all $i$.
By (\ref{eq3.2new}) $y_i \in \mathcal{O} f(x_i)$ for all $i$.  Thus, if $y \in X$ we have $y \in \mathcal{N} f(x)$,
completing the proof of (\ref{eq3.1}).

Now assume that $C = \{ x \} \cup \mathcal{N}f(x)$ is compact but that $y \not\in X$.
Let $U$ be a bounded open neighborhood of
$C $ in $X$. We can assume that $x_i \in U$ and $y_i \not\in \ol{U}$ for all $i$.
For each $i$ there exists $n_i \in \N$ such that
$y_i \in f^{n_i}(x_i)$. Since $ \{ x_i \} = f^0(x_i) \subset U$ we can define
\begin{equation}\label{3.3new}
m_i \quad =_{def} \quad max \{ m \in \Z_+ : m \leq n_i \qquad \mbox{and} \qquad f^{m}(x_i) \subset U \}.
\end{equation}
Because $y_i \in f^{n_i}(x_i) \setminus \ol{U}$ we have $m_i < n_i$. Choose $z_i \in f^{m_i + 1}(x_i) \setminus U$.
By going to a subnet we can assume that $\{ z_i \}$ converges to $z \in \X$.

Because $f$ is +proper, $f(\ol{U})$ is a compact subset of $X$ and it contains the net $\{ z_i \}$. Hence,
$z \in X \cap \mathcal{N}\f(x) = \mathcal{N}f(x)$. But the net $\{ z_i \}$ is in the closed set $X \setminus U$
and so the limit point $z$ does not lie in $C = \mathcal{N}f(x)$.

This contradiction completes the proof.

$\Box$ \vspace{.5cm}

\begin{theo}\label{th1.9new} Let $(X,f)$ be a dynamical system with $f$ a +proper relation
 and let $(\X, \f)$ be an almost dynamic compactification of $(X,f)$.
\begin{enumerate}
\item[(a)] If $C$ is a compact $\G f$ unrevisited subset of $X$
then $C$ is a $\G \f$ unrevisited subset of $\X$.
That is, $ \G \f (C) \cap \G \f^{-1}(C) \subset C$.
\item[(b)] If $C$ is a compact $\G f$ +invariant subset of $X$ then $C$ is a $\G \f$ +invariant subset of $\X$.
That is, $ \G \f (C)  \subset C$.
\end{enumerate}
\end{theo}

{\bfseries Proof:}  (a): Let $\hat C = C \cup ( \G \f (C) \cap \G \f^{-1}(C)) $ so that $\hat C$ is
a compact $\G \f$ unrevisited subset of $\X$.

First, we prove that
\begin{equation}\label{eq1.13new}
C \quad = \quad \hat C \cap X. \hspace{4cm}
\end{equation}

Let $x \in X \cap (\G \f (C) \cap \G \f^{-1}(C))$. So there exist $y_1,y_2 \in C$ such that
$(y_1,x),(x,y_2) \in \G \f$.  Because the compactification is almost dynamic,
$(y_1,x),(x,y_2) \in 1_X \cup \G f$. If $(y_1,x) \in 1_X$ then $x = y_1 \in C$ and similarly
if $(x,y_2) \in 1_X$.  Otherwise, $(y_1,x),(x,y_2) \in  \G f$ and so $x \in C$ because $C$ is
$\G f$ unrevisited.

Thus, it suffices to show that $\hat C \subset X$. Let $z \in \hat C$. Since $\hat C \subset C \cup \G \f(C)$
we have $z \in C$ or else $(x,z) \in \G \f$ for some $x \in C$. We  show that
$z \in C$ in the latter case as well.

Because $\hat C$ is $\G \f$ unrevisited, Theorem \ref{th2.1} implies that $(x,z) \in \mathcal{C} (\f_{\hat C})$.

Because $f$ is +proper, $C \cup f(C)$ is a compact subset of $X$ which in turn is open in $\X$ and $ \f(C) \subset X$ by
Proposition \ref{prop1.8new}. So $\f(C) = f(C)$ because $f = (X \times X) \cap \f$. Hence,
there exists $V $ a neighborhood of the diagonal in $\X \times \X$ such that
$V(C \cup f(C)) \subset X$.  There exists a $V$ chain in $\hat C$ for $\f$ from $x$ to $z$
That is, there is a sequence  $ x_0,y_0,...x_n,y_n \in \hat C$ with $x_0 = x, y_n = z$, $(x_i,y_i) \in \f$
and $(y_{i-1},x_i) \in V$. Assuming inductively that $x_i \in C$ we show that $y_i, x_{i+1} \in C$.
Since, $x_i \in C$, $y_i \in \f(C) = f(C)$ and so $x_{i+1} \in V(f(C))$.  That is,
\begin{equation}\label{eq1.14newa}
y_i, x_{i+1} \quad \in \quad \hat C \cap V(f(C)) \quad \subset \quad \hat C \cap X \quad = \quad C.
\end{equation}
By induction $z = y_n \in C$ as required.

(b): Proceed as in (a) letting $\hat C = C \cup \G \f(C)$.  If $x \in X \cap \hat C$ then there exists
$y \in C$  such that $(y,x) \in 1_X \cup \G f$ because the compactification is almost dynamic. Hence, $x \in C$
because $C$ is $\G f$ +invariant.  That is, $C = X \cap \hat C$. Again Proposition \ref{prop1.8new} implies
that $\f(C) \subset X$ and so $\f(C) = f(C) \subset C$. Since $\hat C$ is $\G \f$ +invariant, it is
$\G \f$ unrevisited and so Theorem \ref{th2.1} applies. As in (a) one proves that $\hat C \subset X$.

$\Box$ \vspace{.5cm}

\begin{cor}\label{cor1.10new} Let $(X,f)$ be a dynamical system with $f$  +proper.
Let  $ (\hat{X},\hat{f})$ be an almost dynamic
compactification of  $(X,f)$.

The compactification $(\X , \f)$ is dynamic. Furthermore,
if $\hat{E} \subset |\mathcal{G} \hat{f}| $ is a $\mathcal{G}
\hat{f} \cap \mathcal{G} \hat{f}^{-1}$ equivalence class  with $E =
\hat{E} \cap X$, then exactly one of the following three possibilities
holds:
\begin{itemize}
\item[(i)]  $\hat{E}$ is a compact subset of $\hat{X} \setminus
X$ and $E = \emptyset$.
\item[(ii)] $E$ is contained in $ |\mathcal{G}f|$  and is a noncompact
$\mathcal{G} f \cap \mathcal{G}f^{-1}$ equivalence class whose
$\hat{X}$  closure meets $\hat{X}\setminus X$ and is contained in
$\hat{E}  $.
\item[(iii)] $\hat{E} = E$ is contained in $|\mathcal{G}f|$ and is a compact
$\mathcal{G} f \cap \mathcal{G} f^{-1} $ equivalence class.
\end{itemize}
\end{cor}

{\bfseries Proof:}  By Theorem \ref{th1.6} and the remark thereafter the corollary
is proved by showing  that case (iv) of the theorem does not occur.  If $E$ is a compact
$1_X \cup (\G f \cap \G f^{-1})$ equivalence class in $X$ then $E$ is $\G f$ unrevisited.
By Theorem \ref{th1.9new} $E$ is $\G \f$ unrevisited and so
$\hat E = \G \f(E) \cap \G \f^{-1}(E) \subset E$. Thus, $\hat E \cap (\X \setminus X) = \emptyset$.

$\Box$ \vspace{.5cm}

\begin{ex}\label{ex1.7}For any  compactification
$(\hat{X},\hat{f})$ of $(X,f)$, the closure in $\hat{X} \times
\hat{X}$ of $\mathcal{G}f$ is contained in $\mathcal{G}\hat{f}$.
It can happen that this inclusion is strict for every proper
compactification of $(X,f)$. \end{ex}
On  $  X = \{ (x,y) \in [0,1]\times [0,1] : xy = 0 \} \cup (0,1] \times \{-1 \} $.
define the continuous map $ f$ :

\begin{equation}\label{eq1.12new}
(x,y) \quad \mapsto \quad
\begin{cases}
(0,y + \frac{1}{2}y(1-y)) \qquad 0 \leq y \leq 1\\
(x,0) \qquad y = -1,0 .
\end{cases}
\end{equation}
Let $(\X, \f)$ be any proper
compactification of $(X,f)$. If $z$ is any limit point in $\X$ of
$(x,-1)$ as $x \to 0$, then $(z,(0,0)) $ is the limit of $ ((x,-1),(x,0)) \in  f$
and so is in the closure of $\f$. Also, $((0,0),(0,1)) \in \G f$ and so $(z,(0,1)) \in \G \f$.
However, $(z,(0,1))$ is not in the closure of $\G f$.

$\Box$\vspace{.5cm}

When $f$ is a continuous map on $X$ we call $(X,f)$ a \emph{cascade}. \index{cascade} A compactification
$(\X, \f)$ is a \emph{cascade compactification} \index{cascade!cascade compactification}%
\index{compactification!cascade compactification} when $\f$ is a map as well. If $(X,f)$ is a reversible
cascade, \index{cascade!reversible} i.e. $f$ is a homeomorphism on $X$, then we call $(\X, \f)$ a \emph{reversible
cascade compactification} when $\f$ is a homeomorphism on $\X$.  Recall that a  compactification
of a cascade is a cascade compactification iff the associated subalgebra $\A$ of $\B(X)$ is $f^*$
+invariant and a compactification of a reversible cascade is a reversible cascade compactification iff
$\A$ is $f^*$ invariant.

A continuous map is a +proper relation and so Proposition \ref{prop1.8new} and Theorem \ref{th1.9new} and
their corollaries apply to a cascade $(X,f)$.  If $(\X,\f)$ is a proper cascade compactification of $(X,f)$
then by Proposition \ref{prop1.8new}  $X \subset \f^{-1}(X) \subset \X$ and so $\{ \f^{-n}(X) : n \in \Z_+ \}$
is an increasing sequence of open, dense subsets of $\X$.  Define
\begin{equation}\label{eq1.14new}
\begin{split}
\tilde X \quad =_{def} \quad \bigcup_{n \in \Z_+} \ \{ \f^{-n}(X) \} \quad \subset \quad \X, \hspace{2cm}\\
\mbox{and so} \qquad \f^{-1}(\tilde X) \quad = \quad \tilde X.\hspace{2.5cm}
\end{split}
\end{equation}  \hspace{4cm}

\begin{prop}\label{prop1.11new} Let $(\X, \f)$  a proper cascade compactification of $(X,f)$.
\begin{enumerate}
\item[(a)]The restriction of $\f$ to
$\tilde X$ is a proper map and
the following are equivalent.
\begin{itemize}
\item[(i)]  The mapping $f$ is a proper map on $X$.
\item[(ii)]  $X$ is $\f^{-1}$ +invariant.
\item[(iii)]  $\tilde X = X$.
\end{itemize}

\item[(b)]  If $A$ is a compact subset of $\tilde X$ then $\f^n(A) \subset X$ for some positive integer
$n$.  In particular, $\X = \tilde X$ iff there exists a positive integer $n$ such that $f^n(X)$ is a bounded
subset of $X$.

\item[(c)]If $f$ is surjective on $X$ then $\f$ is surjective on $\X$ and $\tilde X$
is $\f$ invariant.
\end{enumerate}
\end{prop}

{\bfseries Proof:}  (a): If $A$ is a compact subset of $\tilde X$, the pre-image $\f^{-1}(A)$
is closed and hence compact in $\X$. By (\ref{eq1.14new}) $\f^{-1}(A) \subset \tilde X$ and so
$\f|\tilde X$ is proper.  In particular, (iii) $\Rightarrow$ (i). Proposition \ref{prop1.8new} applied to
the relation $f^{-1}$ implies  (i) $\Leftrightarrow$ (ii). Finally, (ii) $\Rightarrow$ (iii) is obvious.

(b): If $A \subset \tilde X$ is compact then $\{ f^{-n}(X) \}$ is an open cover of $\X$ and so for some $n \ A \subset
\f^{-n}(X)$ and so $\f^n(A)  \subset X$. In particular, if $\X = \tilde X$ then for some
$n \ f^n(X) \subset \f^n(\X) \subset X$.
On the other hand, if $\f^n(X) = f^n(X) \subset B \subset X$ and
$B$ is compact then, taking the closure in $\X$,  $\f^n(\X) \subset \ol{f^n(X)} \subset B \subset X$ and so
$\X \subset \f^{-n}(X) \subset \tilde X$.

(c): If $f$ is surjective then $X = f(X)$ is a subset of the compact set $\f(\X)$. Since $X$ is dense in $\X$,
$\f(\X) = \X$.  That  $\f(\tilde X) = \tilde X$ then follows from (\ref{eq1.14new}).

$\Box$ \vspace{.5cm}

\begin{ex}\label{ex1.8new} For a cascade $(X,f)$ there is a smallest proper cascade compactification. \end{ex}
If $f$ is a proper map and $u \in \B(X)$ has compact support then $f^*(u) = u \circ f$ has compact support
and so the closed subalgebra  $\A_0$ generated by such functions is $f^*$ +invariant. If $f$ is a homeomorphism
on $X$ then $\A_0$ is $f^*$ invariant. Thus, for the one-point compactification $X_*$ of $X$, the
proper compactification $(X_*, f_*)$ is a cascade compactification when $f$ is proper and is a reversible
cascade compactification when $f$ is a homeomorphism.

If $f$ is not proper, then we let $\A_{0.f}$ denote
the closed subalgebra generated by $\{ u \circ f^n \}$ as $u$ varies over the functions of compact support and
$n$ varies over $\Z_+$. If $\X$ is any proper compactification of $X$, and $u \in \B(X)$ has compact
support then $u$ extends to $\hat u \in \B(\X)$ by $\hat u(z) = 0$ for $z \in \X \setminus X$.
If $(\X, \f)$ is any cascade compactification then $u \circ f^n$ extends to $\hat u \circ \f^n$ for
any $n \in \Z_+$. Furthermore, if
 $(\hat u) \circ \f^n(z) > 0$ for some positive integer $n$ then $\f^n(z)$ is
contained in the support of $\hat u$ which is equal to the support of $u$ since the latter is compact in $X$.
Hence, $z \in \tilde X$. Thus, every function in $\A_{0.f}$ is constant on $\X \setminus \tilde X$.
In particular, if $(X_{*.f},f_*)$ denotes the cascade compactification associated
with $\A_{0.f}$ itself then since $\A_{0.f}$ distinguishes the points of $X_{*.f} \setminus \tilde X$ it follows that
the latter is a single point.  Thus, the restriction of $\f$ to $\tilde X$
defines a proper map cascade $(\tilde X, \tilde f)$ which extends $(X,f)$ and
$X_{*.f}$ is the one-point compactification
of $\tilde X$. For any cascade compactification $(\X,\f)$ the associated subalgebra contains $\A_{0.f}$ and
so we obtain a continuous map $\pi : \X \to X_{*.f}$ which maps $(\X,\f)$ to $(X_{*.f},f_*)$. If $z \in \X$ and
$\pi(z) \in \tilde X \subset X_{*.f}$ then $\hat u \circ \f^n(z) = \hat u \circ \f_*^n(\pi(z)) > 0$ for some $n$ and so
$z \in \tilde X \subset \X$. That is, $\X \setminus \tilde{X}$ is the $\pi$ preimage of the point at infinity
of $X_{*.f}$.

$\Box$ \vspace{.5cm}

In the cascade case we can sharpen Corollary \ref{cor1.10new}.

\begin{prop}\label{prop1.12new} Let  $ (\hat{X},\hat{f})$ be an almost dynamic, cascade
compactification of  a cascade $(X,f)$.

The compactification $(\X , \f)$ is dynamic. Furthermore,
if $\hat{E} \subset |\mathcal{G} \hat{f}| $ is a $\mathcal{G}
\hat{f} \cap \mathcal{G} \hat{f}^{-1}$ equivalence class  with $E =
\hat{E} \cap X$, then exactly one of the following three possibilities
holds:
\begin{itemize}
\item[(i)]  $\hat{E}$ is a compact subset of $\hat{X} \setminus
\tilde X$ and $E = \emptyset$.
\item[(ii)] $E$ is contained in $ |\mathcal{G}f|$  and is a noncompact
$\mathcal{G} f \cap \mathcal{G}f^{-1}$ equivalence class and
$\hat{E}$   meets $\hat{X}\setminus \tilde X$.
\item[(iii)] $\hat{E} = E$ is contained in $|\mathcal{G}f|$ and is a compact
$\mathcal{G} f \cap \mathcal{G} f^{-1} $ equivalence class.
\end{itemize}
\end{prop}

{\bfseries Proof:} By Lemma \ref{lem2.2} $\hat E$ is $\f$ invariant. If $\hat E \subset \tilde X$ then
by Proposition \ref{prop1.11new}(b) there exists a positive integer $n$ such that $\f^n(\hat E) \subset X$.
Since $\hat E$ is $\f$ invariant, $\hat E = \f^n(\hat E) \subset X$ and this is case (iii).

Now assume $\hat E$ meets $\X \setminus \tilde X$. If  also there exists $ x \in \hat E \cap \tilde X$ then
$f^n(x) \in \hat E \cap X = E$ for some $n$ and this is case (ii) of Corollary \ref{cor1.10new} and so
is case (ii) here.

There only remains $\hat E \subset \X \setminus \tilde X$ which is case (i).

$\Box$ \vspace{.5cm}

\section{Lyapunov Function Compactifications}\label{seclyap}

In the previous section we described properties of dynamical compactifications, but we
did not demonstrate their existence.  In this section we construct them by using Lyapunov functions

Given a continuous  $L : X \to {\mathbb R}$ we define the
relations \index{$\leq_L$}

\begin{equation}\label{eq1.13}
\begin{split}
\leq_L \quad = \quad \{ (x,y) : L(x) \leq L(y) \} \\
>_L \quad = \quad \{ (x,y) : L(x) > L(y) \} \\
=_L \quad = \quad \{ (x,y) : L(x) = L(y) \}
\end{split}
\end{equation}
so that $\leq_L$ is a closed, reflexive, transitive relation with
associated equivalence relation  $=_L \ ( \ = \ \leq_L \cap
(\leq_L)^{-1} \ )$. The complement of $\leq_L$, $>_L$, is open,
transitive and asymmetric.

A \emph{Lyapunov function} \index{Lyapunov function} for a closed relation $f$ on a space $X$ (also called a Lyapunov
function for the dynamical system $(X,f)$ ) is a
bounded, continuous real-valued function $L$ on $X$ such that $y
\in f(x) $ implies $L(x) \leq L(y)$. This is clearly equivalent to
the condition $f \ \subset \ \leq_L$. Since $\leq_L$ is a closed,
transitive relation, this implies

\begin{equation}\label{eq1.14}
\begin{split}
\mathcal{G}f  \quad \subset \quad \leq_L \qquad \mbox{and}\\
\mathcal{G}f \cap \mathcal{G}f^{-1} \quad \subset \quad =_L.
\end{split}
\end{equation}
It follows that a Lyapunov function for $f$ is automatically a
Lyapunov function for $\mathcal{G}f$.

By composing with an increasing homeomorphism of ${\mathbb R}$
onto the open unit interval we can replace any  $L: X \to {\mathbb
R}$ by a bounded  function with the same relation $\leq_L$. So the
assumption of boundedness in the definition is just a convenience.
In fact, we need only consider Lyapunov functions which map to
$[0,1]$.

Clearly, the set of Lyapunov functions is closed under finite sums,
multiplication by positive constants, and uniform limits.  Thus, they form
a closed cone in $\B(X)$.

The constant functions are Lyapunov functions but there are
usually many more. From (\ref{eq1.14}) it follows that any Lyapunov
function is constant on each $\mathcal{G}f \cap \mathcal{G}f^{-1}$
equivalence class in $|\mathcal{G}f|$. The theorem on which our
results are based is the observation that otherwise the Lyapunov
functions distinguish points. We review the results from Nachbin
(1965) and Auslander (1964).

The following is essentially a special case of a theorem of
Nachbin (1965) Chapter I.

\begin{theo}\label{th1.8} Let $X_0$ be a closed subset of a compact space $X$.
Let $F$ be a closed, transitive relation on $X$ and let $F_0 = F
\cap (X_0 \times X_0)$ be the closed, transitive relation on $X_0$
which is its restriction.  If $L_0 : X_0 \rightarrow [a,b]$ is a
Lyapunov function for $F_0$ then there exists $L : X \rightarrow
[a,b]$ a Lyapunov function for $F$ which extends $L_0$, ie. $L|X_0
= L_0$. \end{theo}

{\bfseries Proof:} We can assume $a = 0$ and $b = 1$. By replacing $F$ by
$F \cup 1_X$ and $F_0$ by $F_0 \cup 1_{X_0}$ we can assume that $F$ and $F_0$
are reflexive as well as transitive.

We mimic the proof of Urysohn's Lemma. Let $\Lambda = \mathbb{Q}
\cap [0,1]$ counted with $\lambda_0 = 0, \lambda_1 = 1$.  We let
$B_1 = F(L_{0}^{-1}(1))$, $B_0 = X$ and define for all $\lambda
\in \Lambda$ sets $B_{\lambda}$ such that
\begin{itemize}
\item $ F(B_{\lambda}) = B_{\lambda}$.
\item $ F(L_{0}^{-1}((\lambda,1])) \subset
(B_{\lambda})^{\circ}.$
\item $  F(L_{0}^{-1}([\lambda,1])) \subset
B_{\lambda}.$
\item $  (F)^{-1}(L_{0}^{-1}([0,\lambda))) \cap B_{\lambda} = \emptyset.$
\item $ \lambda' > \lambda $ implies $B_{\lambda'} \subset \subset B_{\lambda}$.
\end{itemize}

Notice that for all $\lambda$:

\begin{equation}\label{eq1.15}
\quad (F)(L_{0}^{-1}([\lambda,1])) \cap
(F)^{-1}(L_{0}^{-1}([0,\lambda))) \quad = \quad \emptyset
\end{equation}
because $F$ is transitive with restriction $F_0$ and because $L_0$
is an $F_0$ Lyapunov function.

We observe that if
$C$ is a closed set with $F(C) = C$ and $U$ is an open set containing $C$, then there exists a
closed set $C_1 \subset U$ with $F(C_1) = C_1$ and such that $C
\subset \subset C_1$ (Akin (1993) Proposition 2.7(b)).
This follows because $C = F(C) = \bigcap
F(N)$ as $N$ varies over the closed neighborhoods of $C$ by (\ref{eq1.4}). Let
$C_1 = F(N)$.  By Proposition \ref{prop1.2}(g) or (h) we can choose $N$ small enough that
$C_1 \subset U$.

Proceed inductively assuming that $B_{\lambda}$ has been defined
for all $\lambda$ in $\Lambda_n = \{ \lambda_i : i = 0,...,n \}$
with $n \geq 1$. Let $\lambda = \lambda_{n+1}$ and let $\lambda' <
\lambda < \lambda''$ the  nearest points in $\Lambda_{n}$
below and above $\lambda$.

Choose a sequence $\{ t_{n}^{-} \}$ with $t_{0}^{-} = \lambda'$,
increasing with limit $\lambda$ and $ \{ t_{n}^{+} \}$ with
$t_{0}^{+} = \lambda''$, decreasing with limit $\lambda$.

Define $Q_{0}^{-} = B_{\lambda'}$ and $Q_{0}^{+} = B_{\lambda''}$.
Inductively, apply Akin (1993) Proposition 2.7(b) to choose
$Q_{n}^{+}$ and then $Q_{n}^{-}$ for $n = 1,2,...$ so that
$F(Q_{n}^{\pm}) = Q_{n}^{\pm}$ and

\begin{equation}\label{eq1.16}
\begin{split}
F(L_{0}^{-1}([t_{n}^{+},1]) \cup Q_{n-1}^{+} \quad \subset
\subset \quad Q_{n}^{+} \quad \subset \subset \quad Q_{n-1}^{-}
\setminus
(F)^{-1}(L_{0}^{-1}([0,\lambda]), \\
F(L_{0}^{-1}([\lambda,1]) \cup Q_{n}^{+} \quad \subset \subset
\quad Q_{n}^{-} \quad \subset \subset \quad Q_{n-1}^{-} \setminus
(F)^{-1}(L_{0}^{-1}([0,t_{n}^{-}]).
\end{split}
\end{equation}
Finally, define

\begin{equation}\label{eq1.17}
B_{\lambda} \quad = \quad \bigcap_n Q_{n}^{-}, \hspace{2cm}
\end{equation}
so that

\begin{equation}\label{eq1.18}
 B_{\lambda} \quad \supset  \quad
\bigcup_n Q_{n}^{+}.\hspace{2cm}
\end{equation}

It is easy to check that $B_{\lambda}$ satisfies the required conditions, thus extending the definitions
to $\Lambda_{n+1}$. By induction they can be defined on the entire set $\Lambda$.

Having defined the $B_{\lambda}$'s we proceed as in Urysohn's Lemma to
define $L(x)$ by the Dedekind cut associated with $x$.  That is,

\begin{equation}\label{eq1.19}
L(x) \quad = \quad inf \{ \lambda : x \not\in B_{\lambda} \} \quad
= \quad sup \{ \lambda : x \in B_{\lambda} \}.
\end{equation}
Continuity follows as in  Urysohn's Lemma.
Because each $B_{\lambda}$ is $F  $ invariant, $L$ is a Lyapunov
function.  The additional conditions on these sets imply that if
$x \in X_0$ then $x \in B_{\lambda}$ iff $\lambda \leq L_0(x)$.
Hence, $L$ is an extension of $L_0$.

$\Box$ \vspace{.5cm}

From this we obtain results for a general space, i.e. a locally
compact, $\sigma$ compact space.

\begin{cor}\label{cor1.9} Let $X_0$ be a compact subset of  a space $X$.
Let $F$ be a closed, transitive relation on $X$ and let $F_0 = F
\cap (X_0 \times X_0)$ be the closed, transitive relation on $X_0$
which is its restriction.  If $L_0 : X_0 \rightarrow [a,b]$ is a
Lyapunov function for $F_0$ then there exists $L : X \rightarrow
[a,b]$ a Lyapunov function for $F$ which extends $L_0$, ie. $L|X_0
= L_0$. \end{cor}

{\bfseries Proof:} Let $\{ K_n \}$ be a  sequence of compact sets
with union $X$ and with $X_0 = K_0$ such that $K_{n-1} \subset
\subset K_n$ for $n = 1,2,...$. Let $F_n = F \cap (K_n \times
K_n)$ for $n = 0,1,...$. Apply the theorem inductively to extend
the $F_n$ Lyapunov function $L_n$ on $K_n$ to an $F_{n+1}$
Lyapunov function $L_{n+1}$ on $K_{n+1}$. The union of the
functions $L_n$ is the required extension $L$.  Notice that $L$ is
continuous on $X$ because each $K_n$ is a neighborhood of
$K_{n-1}$.

$\Box$ \vspace{.5cm}

As a corollary we obtain the usual Urysohn Lemma analogue for
Lyapunov functions.

\begin{cor}\label{cor1.10} Let $A, B$ be a pair of disjoint closed subsets of a
space $X$. Let $F$ be a closed, transitive relation on $X$. If
$F(A) \subset A$ and $F^{-1}(B) \subset B$ then there exists $L :
X \rightarrow [0,1]$ a Lyapunov function for $F$ such that $L(x) =
1$ if $x \in A$ and $L(x) = 0$ if $x \in B$.
\end{cor}

{\bfseries Proof:} Let $\{ K_n \}$ be a  sequence of compact sets
with union $X$ and with $K_0 = \emptyset$ such that such that
$K_{n-1} \subset \subset K_n$ for $n = 1,2,...$. Let $F_n = F \cap
(K_n \times K_n)$ for $n = 0,1,...$. For $n = 0,1,...$ let $K_{n
+\frac{1}{2}} = K_n \cup [(A \cup B)\cap K_{n+1}]$ and $F_{n
+\frac{1}{2}} = F \cap (K_{n +\frac{1}{2}} \times K_{n
+\frac{1}{2}})$. If $L_n$ is an $F_{n}$ Lyapunov function on $K_n$
which is $1$ on points of $A$ and $0$ on points of $B$ then we
obtain an $F_{n+\frac{1}{2}}$ Lyapunov function
$L_{n+\frac{1}{2}}$ on $K_{n+\frac{1}{2}}$ by using $1$ on the new
points of $A$ and $0$ on the new points of $B$. Now apply the
extension theorem to $K_{n+\frac{1}{2}}$in order to get an $F_{n+1}$ Lyapunov function
$L_{n+1}$ on $K_{n+1}$. Again $L$ is the union of the functions
$L_{n}$.

$\Box$ \vspace{.5cm}

From this we obtain Theorem 4 of Auslander (1964).

\begin{cor}\label{cor1.11} Let $F$ be a closed, transitive relation on a space $X$.
 If $(x,y) \in (X \times X) \setminus (F \cup 1_X) $ then there exists a
 Lyapunov function $L : X \rightarrow [0,1]$ for $F$ with $L(x) = 1$ and $L(y) = 0$.
\end{cor}

{\bfseries Proof:} With $F_1 = F \cup 1_X$ let $A = F_1(x)$ and $B
= (F_1)^{-1}(y)$. Since $F_1$ is transitive and does not contain
$(x,y)$ these closed sets are disjoint. Apply the previous
corollary.

$\Box$ \vspace{.5cm}

\begin{theo}\label{th1.12}  Let $F$ be a closed, transitive relation on $X$.
Let $A$ be a compact $F \ +$invariant subset of
$X$. Assume that $A$ admits a compact neighborhood $U$
such that the closed set $F(U)$ is compact.  There
exists $L : X \to [0,1]$ a Lyapunov function for $F$ such that $L(x) = 1$ for all $x \in A$ and $L$ has compact
support, i.e. $L(x) = 0$ for all $x$ outside of some compact set.
\end{theo}

{\bfseries Proof:} Let $W$ be an open set with $A \subset W \subset U$.  Let
$G \ = \ \{ x  \in W : F(x) \subset W \} $.
Because $x \in W$ implies $F(x) \subset F(U)$
we have that $G =  W \setminus F^{-1}(F(U) \setminus W)$. Because $F(U) \setminus W$ is compact,
 $G$ is an open set.
 If $y \in F(x)$ with $x \in G$ then by transitivity $F(y) \subset F(x) \subset W$.
 Thus, $G$ is $F \ +$invariant.
 It follows that
 $B = X \setminus G$ is a closed subset which is $F^{-1} \ +$invariant.

By Corollary \ref{cor1.10} there is a Lyapunov function $L$ which is zero on $B$
and one on $A$. Since $W$ is bounded,
$L$ has compact support.

$\Box$ \vspace{.5cm}

Notice that when $A$ is not invariant, it does not suffice
for this result that $A$ have a neighborhood $U$ such that
$F(U)$ is compact.  A neighborhood of $F(A)$ is needed.
This is clear from the conclusion because if $L$ is a
Lyapunov function with compact support such that $L$ is $1$
on $A$ then $L^{-1}[\frac{1}{2},\infty))$ is a compact
+ invariant neighborhood of $A$ and so of $F(A)$

Suppose $x \in X$ such that $F(x)$ is compact but for every
neighborhood $U$ of $F(x)$, the set $F(U)$ is unbounded.
This occurs  in Example \ref{ex1.7} with the point $x = (0,0)$ and $F = \G f^{-1}$.  In any
compactification $(\hat{X},\hat{F})$ of $(X,F)$, if
$z \in \hat{X} \setminus X$ is a limit point of the sets
$F(U)$ then $z \in \hat{F}(x)$.  Now define a new
 space $\tilde{X} = X \cup \{ e \}$ obtained by adjoining
 the isolated point $e$ and define $\tilde{F} =
 F \cup \{ e \} \times (\{ e \} \cup F(x))$.  In this space,
 $e$ is a clopen set with
 $\tilde{F}(e) = \{ e \} \cup F(x)$ which is compact.
 However, any compactification contains points at
 infinity which are related to $e$.

\begin{ex}\label{ex1.13} Corollary \ref{cor1.9} can fail if $X_0$ is merely a closed
subset. \end{ex}
In $\mathbb{R}^2$ define the closed sets
$X_0 \subset X$, and the real-valued map $\pi$  and the relation
$F$ on $X$:

\begin{equation}\label{eq1.20}
\begin{split}
X_0 \quad = \quad \{ (x,y) : |x| \leq 1 \ \mbox{and} \ xy = 1 \},
\hspace{1cm}
\\ X \quad = \quad X_0 \ \ \cup \ \  \{ (x,0) : |x| \leq 1 \}, \hspace{1.5cm}\\
\pi(x,y) \quad = \quad x.\hspace{4cm} \\
F \quad = \quad =_{\pi} \quad = \quad   \hspace{2.5cm} \\
\{ ((x_1,y_1),(x_2,y_2)) \in X \times X : x_1 = x_2 \}.
\end{split}
\end{equation}
Since $F$ is a closed equivalence relation, a Lyapunov function is
a continuous, real-valued function which is constant on
equivalence classes. Since the restriction $F_0 $ to $X_0$ is
 $1_{X_0}$, any continuous, real-valued function on $X_0$ is a
Lyapunov function for $F_0$. Define $L_0(x,y) = 0$ if $y < 0$ and
$ = 1$ if $y > 0$. This does not extend to a Lyapunov function for
$F$ on $X$.  The problem is that while $F_1(L_0^{-1}(1)) $ and
$(F_1)^{-1}(L_0^{-1}(0))$ are disjoint, their closures meet.
Notice too that while each equivalence class is compact, the
closed set $F([-1,1])$ is not compact.

%In this example, we can use, instead of the equivalence relation
%$F$, the  map
%\begin{equation}\label{eq1.21}
%\begin{split}
%f \quad =_{def} \quad \{((0,0),(0,0) \} \quad \cup  \hspace{3cm}\\ \{ ((x,y),(x,1/x)) : (x,y) \in X \quad
%\mbox{and} \quad x \not= 0 \}.\hspace{1cm}
%\end{split}
%\end{equation}
%While $f$ is discontinuous it is a closed, transitive relation and
%so $f = \mathcal{G}f$.  Hence, $|\mathcal{G}f| = |f| = \{ (0,0) \} \cup X_0$. The
%restriction of $f$ to $X_0$ is the identity $1_{X_0}$ and so $L_0$
%is a Lyapunov function for $f$.

$\Box$  \vspace{.5cm}

\begin{lem}\label{lem1.14}  Let $F$ be a closed, transitive relation on $X$. If $x \not\in |F|$ then there
exists a compact neighborhood $U$ of $x$ such that the three sets $U, F(U)$ and $F^{-1}(U)$ are pairwise disjoint.
\end{lem}

{\bfseries Proof:}  Since $(x,x) \not\in F$ we can choose a compact neighborhood
$U_0$ of $x$ such that $U_0 \times U_0$ is disjoint from the closed set $F$.
Hence, $U_0$ is disjoint from $F(U_0)$ and $F^{-1}(U_0)$.

Now assume that for each compact neighborhood $U \subset U_0$ of $x$, there exists $b_U \in
F(U) \cap F^{-1}(U)$ so that there exist $a_U, c_U \in U $ with $(b_U,a_U)$ and $(c_U,b_U) $ in $F$.
By transitivity, $(c_U,a_U) \in F$. As $U$ shrinks to $\{ x \}$ these pairs approach $(x,x)$ and so
$(x,x) \in F$ because $F$ is closed.

$\Box$ \vspace{.5cm}

If $F$ is a closed transitive relation and $U$ is a compact subset of $X$ such that
 $F(U) \cap F^{-1}(U) = \emptyset$ then by Proposition \ref{prop1.2} (a)
 $A = F(U) $ and $B = F^{-1}(U)$ are disjoint closed
 subsets of $X$ with $F(A) \subset A$ and $F^{-1}(B) \subset B$.   By Corollary \ref{cor1.10} there
 exists a Lyapunov function  $L : X \to [0,1]$ for $F$ such that $L(x) = 1$ for $x \in A$ and
 $L(x) = 0$ for $x \in B$.  In general, we will say a Lyapunov function $L$ for $F$ satisfies
  \emph{ splitting for $U$ } \index{splitting} when:
 \begin{equation}\label{eq1.22}
 sup \ L|F^{-1}(U) \quad < \quad inf \ L|F(U).
 \end{equation}

By applying Corollary \ref{cor1.11} to $F = \mathcal{G}f$ we see that for
any closed relation $f$ on $X$

\begin{equation}\label{eq1.23}
1_X \cup \mathcal{G}f  \quad = \quad \bigcap \{ \leq_L \}
\end{equation}
with $L$ varying over all Lyapunov functions for $f$. Thus, the
open sets $ >_{L} $ form an open cover of $(X \times X) \setminus
(1_X \cup \mathcal{G}f) $ as $L$ varies over the Lyapunov
functions.

\begin{df}\label{def1.15} Let $\mathcal{L}$ be a set of bounded Lyapunov functions
 for $f$.  We say that $\mathcal{L}$ a \emph{sufficient set of Lyapunov functions}
 \index{sufficient set of Lyapunov functions}%
 when

\begin{equation}\label{eq1.24}
1_X \cup \mathcal{G}f  \quad = \quad \bigcap_{\mathcal{L}} \{
\leq_L \}
\end{equation}
or, equivalently, when $\{ >_{L} : L \in \mathcal{L} \}$ is an
open cover of $(X \times X) \setminus (1_X \cup \mathcal{G}f) $.

We say that $\L$  {\em satisfies the splitting condition} when
 for every $x \in X \setminus |\G f|$ there is a compact
neighborhood $U$ of $x$ and a Lyapunov function $L \in \L$
which is splitting for $U$.
\index{splitting condition}
\end{df} \vspace{.5cm}

Notice that if the relation $\mathcal{G} f$ is reflexive then the splitting condition
follows vacuously because $X \setminus |\G f| = \emptyset$.

When $X$ is metrizable, and so $X \times X$ is as well, the spaces
are second countable and hence Lindel\"{o}f. By Corollary \ref{cor1.11} we can choose a
sequence $\{ L_n \}$ such that $\{>_{L_{n}} \}$ covers $(X \times
X) \setminus (1_X \cup \mathcal{G}f) $, and so $\mathcal{L} = \{
L_n \}$ is a countable sufficient set of Lyapunov functions.
In addition, by Lemma \ref{lem1.14} we can choose a sequence $\{ U_n \}$ of
compact subsets of $X \setminus |\G f|$ whose
interiors cover $X \setminus |\G f|$ and such that $\G f^{-1}(U_n) \cap \G f(U_n) = \emptyset$.
By Corollary \ref{cor1.10} we can choose for each $U_n$  a  Lyapunov function splitting for $U_n$.
Adjoining these to the previous
sequence we obtain a countable sufficient set of Lyapunov functions for $f$ which also satisfies the
splitting condition.

 \begin{lem}\label{lem1.16} Let $(X.f)$ be a dynamical system with
  $(\hat{X},\f)$  a proper compactification.
 If $U$ is a compact subset
 of $X$ with $\G f(U) \cap \G f^{-1}(U) = \emptyset$ and some
 $\G f$ Lyapunov function $L$ which is splitting
 for $U$ extends
 continuously to $\hat{X}$, then the interior of $U$ is disjoint from $|\mathcal{G}\hat{f}|$.
 \end{lem}

 {\bfseries Proof:}  Let $a = sup \ L| \G f^{-1}(U)$ and $b = inf \ L|\G f(U)$. By (\ref{eq1.22}) $a < b$.

 Because $\f $ is a closed relation on a compact space, it follows from
  Proposition \ref{prop1.3} that $\mathcal{G}\f = \f \cup \mathcal{G}\hat{f} \circ \hat{f}$.

 Suppose $x $ lies in the interior of $U$ and $(x,x) \in \mathcal{G}\hat{f}$.  Since
 $\hat{f} \cap (X \times X) = f$ and $\G f(U) \cap \G f^{-1}(U) = \emptyset$, we have $(x,x) \not\in \f$.
 By Lemma \ref{lem1.4} there exists $z \in \hat{X}$ such that $(x,z) \in \hat{f}$ and $(z,x) \in \mathcal{G}\hat{f}$.
 Since $\hat{f}$ is the closure of $f$ it follows that there is a net of points $(x_i,z_i) \in f$ which converge
 to $(x,z)$. The net $\{ x_i \}$ eventually enters the interior of $U$ and so $z_i$ is eventually in $f(U)$.  Hence,
 $L(z_i) $ is eventually greater than or equal to $b$.
 Letting $\hat{L}$ denote the extension of $L$ to $\hat{X}$ we have that $\hat{L}(z) \geq b$.
Now $\hat{L}$ is a Lyapunov function for $\mathcal{G}\hat{f}$ and so $(z,x) \in \mathcal{G}\hat{F}$ implies
$L(x) \geq b$.

We can apply the same argument to $f^{-1}$.  That is,
$\mathcal{G}\hat{f} = \hat{f} \cup \hat{f} \circ \mathcal{G}\hat{f}$ and so there exists $z' \in \hat{X}$
such that $(z',x) \in \hat{f}, (x,z') \in \mathcal{G}\hat{f}$.  Proceeding as before we obtain $\hat{L}(z') \leq a$
and so $L(x) \leq a$.  Because $a < b$ we obtain a contradiction showing that $(x,x) \not\in \mathcal{G}\hat{f}$.

$\Box$ \vspace{.5cm}

For $\L$  any set of bounded continuous functions on $X$, consider  the closed subalgebra $\A$ of
$\mathcal{B}(X)$ which is generated by the functions in
$\mathcal{L}$ together with the  functions of compact support.
 As $\mathcal{A}$ contains all the functions in $\B(X)$ with compact support,
it distinguishes points and closed sets.  Let $\hat{X}$ be the associated Gelfand space and $j : X \to
\hat{X}$ be the associated proper compactification which we will
regard as an inclusion. We will call this the  $\mathcal{L}$
compactification.
\index{compactification!$\mathcal{L}$ compactification}%
\index{$\mathcal{L}$ compactification} Since $X$ is locally compact, it is an open
dense subset of $\hat{X}$. For any $u \in \mathcal{A}$ we will denote by $ \hat{u}$
 the extension to $\hat{X}$.  In particular, for any $L \in
\mathcal{L}$ we have $\hat{L}$ defined on $ \hat{X} $.  If $X$ is metrizable and $\L$ is a countable set then
the algebra $\A$ is countably generated and so  the compactification $\X$ is metrizable.

Suppose that $\L$ is a set of Lyapunov functions for a closed relation $f$.
As before, the relation $\hat{f}$ is the closure in $\hat{X} \times \hat{X}$
of $f$ and $\mathcal{G} \hat{f} $  is the smallest closed transitive
relation on $\hat{X}$ which contains $\hat{f}$. For each $L \in \L$, the relation
$\leq_{\hat{L}}$ is a closed, transitive relation on $\X$ which contains $\leq_L$ and so
contains $f$.  It follows that

\begin{equation}\label{eq1.25}
\mathcal{G} \hat{f} \quad \subset \quad \bigcap_{\mathcal{L}} \{
\leq_{\hat{L}} \}. \hspace{1cm}
\end{equation}
Thus,  the $\hat{L}$'s are Lyapunov functions for $\mathcal{G} \hat{f}$ and a fortiori for $\f$.

\begin{lem}\label{lem1.17}  Let $\L$ be a set of bounded continuous functions on $X$.
Let $(\hat{X},\hat{f})$ be the $\mathcal{L}$  compactification
for $(X,f)$. If $z_1$ and $z_2$ are distinct points of $\hat{X}
\setminus X$ then $ \ \hat{L}(z_1) \not= \hat{L}(z_2)$ for some $
\ L \in \mathcal{L}$.

If $\L$ is a set of Lyapunov functions for a closed relation $f$ on $X$, then
$z_1$ and $z_2$ do not lie in
the same $\mathcal{G} \hat{f} \cap \mathcal{G} \hat{f}^{-1}$
equivalence class.
\end{lem}

{\bfseries Proof:}  Distinct points of $\hat{X} \setminus X$ are
distinguished by some member of $\mathcal{A}$.  If $u$ has compact support
then $\hat{u} = 0$  on $\hat{X} \setminus X$. So
the points must be distinguished by one of the $\hat{L}$'s.

When $\L$ consists of Lyapunov functions the $\hat{L}$'s are constant on $\mathcal{G} \hat{f}
\cap \mathcal{G} \hat{f}^{-1}$ equivalence classes.  It follows that
distinct points of $\hat{X} \setminus X$ cannot lie in the same
class.

$\Box$ \vspace{.5cm}

\begin{df}\label{df1.19}  A {\em Lyapunov function compactification for a closed relation} $f$ on $X$ is an
$\L$ compactification for $\L$ a sufficient set of Lyapunov functions for $f$.
\index{compactification!Lyapunov function compactification}%
\index{Lyapunov function compactification}
\end{df}
 \vspace{.5cm}

\begin{theo}\label{th1.18} Assume that  $f$ is a closed relation on
 $X$ and that $ \mathcal{L}$ is a sufficient set of Lyapunov
 functions for $f$.
Let $(\hat{X},\hat{f})$ be the Lyapunov function
compactification of $(X,f)$ associated with $\L$.
\begin{enumerate}
\item[(a)]  $(\hat{X},\hat{f})$ is an almost dynamic
compactification of $(X,f)$, i.e.

\begin{equation}\label{eq1.26}
(1_{\hat{X}} \cup \mathcal{G} \hat{f}) \cap (X \times X) \quad =
\quad 1_X \cup \mathcal{G} f.
\end{equation}

\item[(b)] If $\hat{E} \subset |\mathcal{G} \hat{f}| $ is
 a $\mathcal{G} \hat{f} \cap
\mathcal{G} \hat{f}^{-1}$ equivalence class with $E = \hat{E} \cap
X$, then exactly one of the following four possibilities holds:
\begin{itemize}
\item[(i)]  $\hat{E}$ consists of a single point of $\hat{X} \setminus
X$.
\item[(ii)] $E$ is contained in $ |\mathcal{G}f|$  and is a noncompact
$\mathcal{G} f \cap \mathcal{G}f^{-1}$ equivalence class with
$\hat{E} \ $its one point compactification. That is, there is a
noncompact equivalence class $E \subset |\mathcal{G} f|$ whose
closure in $\hat{X}$ is $\hat{E}$ and $\hat{E} \setminus E$ is a
singleton.
\item[(iii)] $\hat{E} = E$ is contained in $|\mathcal{G}f|$ and is a compact
$\mathcal{G} f \cap \mathcal{G} f^{-1} $ equivalence class.
\item[(iv)]  $\hat{E}$ is the union of a  compact $1_X \cup (\mathcal{G}f \cap
\mathcal{G}f^{-1})$ equivalence class $E$ and a single point of
$\hat{X} \setminus X$.
\end{itemize}

\item[(c)] If $x,y \in |\mathcal{G}f|$  lie in distinct $\mathcal{G}f \cap
\mathcal{G}f^{-1}$ equivalence classes then their equivalence
classes have disjoint closures in $\hat{X}$.

\item[(d)] If the set $\L$ satisfies the splitting condition, then $(\X , \f)$ is a
dynamic compactification of $(X,f)$,
i.e. $(\G \f) \cap (X \times X) = \G f$.

\item[(e)] If $f$ is a +proper relation then $(\X, \f)$ is a dynamic compactification and
case (iv) of (b) does not occur.
\end{enumerate}
\end{theo}

{\bfseries Proof:} (a): Clearly, $\mathcal{G}f \subset \mathcal{G}
\hat{f} \cap (X \times X)$. On the other hand, if   $(x_1,x_2) \in
(X \times X) \setminus (1_X \cup \mathcal{G}f) $ then by (\ref{eq1.24})
there exists $L \in \mathcal{L}$ so that $L(x_1) > L(x_2)$ and so
by (\ref{eq1.25}) $(x_1,x_2) \not\in \mathcal{G} \hat{f}$. Hence (\ref{eq1.26})
holds.

(b): Since $(\hat{X},\hat{f})$ is an almost dynamic compactification,
we can apply Theorem \ref{th1.6}. In addition, Lemma \ref{lem1.17} implies that
$\hat{E} \cap (\hat{X} \setminus X)$ is either a singleton or the
empty set.

It follows that in case (i) with $E$ empty $\hat{E}$ is a single
point of $\hat{X} \setminus X$.

In case (ii) and case (iv) as well $\hat{E} \setminus E = \hat{E}
\cap (\hat{X} \setminus X)$ is nonempty and so is a single point.

There remains case (iii) as described in Theorem \ref{th1.6}.

(c): Distinct  $\mathcal{G}f \cap
\mathcal{G}f^{-1}$ classes have distinct closures in $\X$ by Theorem \ref{th1.6} (c).

(d): Assume that $\L$ satisfies the splitting condition.  If $x \in X \setminus |\G f|$ then
by hypothesis there exists a compact neighborhood $U$ of $x$ and a $\G f$  Lyapunov function splitting
for $U$ which is in $\L$.  Since the elements of $\L $ extend to $\X$,
Lemma \ref{lem1.16} implies that $x$, an interior point
of $U$, is not in $| \G \f|$. Together with (\ref{eq1.26}) this implies that the compactification is dynamic.

(e): If $f$ is a +proper relation then the results follow from Corollary \ref{cor1.10new}.

$\Box$ \vspace{.5cm}

From Theorem \ref{th1.18} it follows that a Lyapunov function compactification
 is an almost dynamic compactification which is dynamic when $f$ is +proper or
 $\L$ satisfies the splitting condition. In
fact, a compactification is dynamic exactly when the associated
algebra contains a sufficient set of Lyapunov functions which satisfies the splitting condition.

\begin{theo}\label{th1.20} For $f$ a closed relation on  $X$ let
$ (\hat{X},\hat{f})$ be a proper compactification  of $(X,f)$.
With $j : X \to \hat{X}$  the inclusion map let  $\mathcal{A} =
j^*(\mathcal{B}(\hat{X}))$  so that  $\mathcal{A}$ is a the closed
subalgebra of  $\mathcal{B}(X)$ consisting of those functions
which extend to $\hat{X}$.

$(\hat{X},\hat{f})$ is an almost dynamic compactification iff
$\mathcal{A}$ contains a sufficient set of Lyapunov functions for
$f$.

$(\hat{X},\hat{f})$ is a dynamic compactification iff
$\mathcal{A}$ contains a sufficient set of Lyapunov functions for
$f$ which also satisfies the splitting condition.

The following conditions are equivalent:
\begin{itemize}
\item[(i)] $(\hat{X},\hat{f})$ is a
Lyapunov function compactification of $(X,f)$.
\item[(ii)]$(\hat{X},\hat{f})$ is an almost dynamic compactification such that no two
distinct points of $\hat{X} \setminus X$ are $\mathcal{G} \hat{f}
\cap \mathcal{G} \hat{f}^{-1}$ equivalent.
\item[(iii)] The Lyapunov functions for $f$ in $\mathcal{A}$ is a sufficient set of Lyapunov
functions for $f$  and $\mathcal{A}$ is generated by these Lyapunov functions
together with the functions of compact support.
\end{itemize}
\end{theo}

{\bfseries Proof:} Let $\mathcal{L} \subset \mathcal{A}$ denote
the restrictions to $X$ of the $\hat{f}$ Lyapunov functions
$\hat{L} : \hat{X} \to [0,1]$. Note that $L \in \mathcal{A}$ is a
Lyapunov function for $f$ iff its extension $\hat{L}$ is a
Lyapunov function for $\hat{f}$. This because $\leq_{\hat{L}}$ is
a closed transitive relation on $\hat{X}$ whose intersection with
$X \times X$ is $\leq_L$. Hence, $\leq_{L}$ contains $f$ iff
$\leq_{\hat{L}}$ does, in which case the latter contains the
closure $\hat{f}$. It follows that $\mathcal{L}$ is exactly the
set of $f$ Lyapunov functions which are contained in
$\mathcal{A}$. Thus, $\mathcal{A}$ contains a sufficient set of
Lyapunov functions for $f$  precisely when
$\mathcal{L}$ is sufficient, i.e. (\ref{eq1.24}) holds. Similarly, $\A$
contains a set of Lyapunov functions which satisfy the splitting condition
interiors cover $X \setminus |\G f|$ iff $\L$ is such a set.

Intersecting equation (\ref{eq1.23}) for $\hat{f}$ with $X \times X$ we
obtain

\begin{equation}\label{eq1.27}
(1_{\hat{X}} \cup \mathcal{G}\hat{f})\cap (X \times X)  \quad =
\quad \bigcap_{\mathcal{L}} \{ \leq_L \}.
\end{equation}
Hence, (\ref{eq1.24}) holds iff (\ref{eq1.11}) does.  That is, $\mathcal{L}$ is
sufficient iff $(\hat{X},\hat{f})$ is almost dynamic.

If $\L$ also satisfies the splitting condition then it follows from Lemma \ref{lem1.16}, as in
the proof of Theorem \ref{th1.18}, that $(\X, \f)$ is dynamic.

Conversely, if $(\X, \f)$ is dynamic and $x \in X \setminus |\G f|$ then $x \not\in |\G \f|$ and
so we can apply Lemma \ref{lem1.14} and Corollary \ref{cor1.10}
to $\G \f$ and obtain a compact neighborhood $U$ for $x$ and a Lyapunov
function $\hat{L} : \X \to [0,1]$ with $\hat{L}$ equal to $1$ on $\G \f(U)$ and equal to $0$ on
$\G \f^{-1}(U)$.  By shrinking $U$ if necessary we can assume that $U$ is contained in the open subset
$X$ of $\X$.  Then the restriction $L = \hat{L}| X \ \in \ \L$ is splitting for $U$.  It follows that
$\L$ satisfies the splitting condition.

(i) $\Rightarrow$ (ii) follows from Theorem \ref{th1.18}.

(ii) $\Rightarrow$ (iii)  The set $\mathcal{L}$ is sufficient
because $(\hat{X},\hat{f})$ is almost dynamic.  Assuming (ii), the closed
subalgebra of $\mathcal{B}(\hat{X})$ generated by the Lyapunov
functions for $\hat{f}$ and the functions with support in $X$
distinguish the points of $\hat{X}$ because the $\mathcal{G}
\hat{f} \cap \mathcal{G} \hat{f}^{-1}$ classes intersect $\hat{X}
\setminus X$ in singletons. By the Stone-Weierstrass Theorem it is
all of $\mathcal{B}(\hat{X})$.

(iii) $\Rightarrow$ (i) is just the definition of a Lyapunov
function compactification.

$\Box$ \vspace{.5cm}

Now assume that $(X,f)$ is a cascade so that $f$ is a continuous map on $X$.
 If $L : X \to {\mathbb R}$ is a Lyapunov function for $f$
then so is $L \circ f^m$ for every $m \in {\mathbb Z}_+$ and this
is true for every $m \in {\mathbb Z}$ if $f$ is a homeomorphism.
If $\mathcal{L}_0$ is a sufficient set  of Lyapunov functions for
$f$ then $\mathcal{L} = \{ L \circ f^m : m \in {\mathbb Z}_+, L
\in \mathcal{L}_0 \}$ is an $f^* $ +invariant sufficient set of
Lyapunov functions. When $f$ is a homeomorphism we can let $m$
vary over ${\mathbb Z}$ to get an $f^*$ invariant sufficient set
of Lyapunov functions.  In each case if $\mathcal{L}_0$ is
countable then $\mathcal{L}$ is countable as well. When $f$ is a proper map and
$\mathcal{L}$ is $f^* \ +$invariant (or $f^*$ invariant) then the
algebra generated by $\mathcal{L}$ and the functions of compact
support is $f^* \ +$invariant (resp. $f^*$ invariant) and so the
$\mathcal{L}$ compactification of $(X,f)$ is a cascade
compactification (and is reversible in the $f^*$ invariant case).
If $f$ is not proper then, as in Example \ref{ex1.8new},  we must include all $u \circ f^n$ for $n \in \Z_+$ and
$u$ with compact support.  We call these the \emph{$f$ translations of functions of compact
support}.

We will call $(\hat{X},\hat{f})$ \emph{Lyapunov function
cascade compactification}
\index{compactification!Lyapunov function cascade compactification}%
\index{Lyapunov function cascade compactification}%
\index{cascade!Lyapunov function cascade compactification} for $(X,f)$ when it is the
 compactification for a cascade $(X,f)$ associated with the closed subalgebra $\A$ generated
 by $\mathcal{L}$ an $f^* \ +$invariant sufficient set of Lyapunov
functions for $f$ together with the translations of functions of compact support.

\begin{theo}\label{th1.21new} Assume that $(X,f)$ is a cascade and $\L$ is an
$f^*$ +invariant sufficient set of Lyapunov functions for $f$. Let $(\X, \f)$ be the
associated Lyapunov function cascade compactification. Define $\tilde X = \bigcup_{n \in \Z_+} \{\f^{-n}(X) \}$.
\begin{enumerate}
\item[(a)] $(\X, \f)$ is a dynamic cascade compactification with $|\G \f| \cap (\X \setminus \tilde X) \subset |\f|$.

\item[(b)] If $\hat{E} \subset |\mathcal{G} \hat{f}| $ is a $\mathcal{G}
\hat{f} \cap \mathcal{G} \hat{f}^{-1}$ equivalence class  with $E =
\hat{E} \cap X$, then exactly one of the following three possibilities
holds:
\begin{itemize}
\item[(i)]  $\hat{E}$ consists of a  single fixed point of $\f$ contained in $\hat{X} \setminus
\tilde X$ and $E = \emptyset$.
\item[(ii)] $E$ is contained in $ |\mathcal{G}f|$  and is a noncompact
$\mathcal{G} f \cap \mathcal{G}f^{-1}$ equivalence class and
$\hat{E} \cap \hat{X}\setminus \tilde X$ is a single fixed point of $\f$.
\item[(iii)] $\hat{E} = E$ is contained in $|\mathcal{G}f|$ and is a compact
$\mathcal{G} f \cap \mathcal{G} f^{-1} $ equivalence class.
\end{itemize}

\item[(c)] If $f$ is a proper map on $X$ then $X = \tilde X$ and $(\X, \f)$ is the $\L$ compactification.

\item[(d)] If $f$ is a homeomorphism and $\L$ is $f^*$ invariant then the $\L$ compactification $(\X, \f)$ is a
reversible cascade.

\item[(e)] If $X$ is metrizable and/or $f$ is reversible then there exist Lyapunov function
compactifications for $(X,f)$ with the same properties.
\end{enumerate}
\end{theo}

{\bfseries Proof:} (a), (b): Since $\A$ is $f^*$ +invariant, $(\X, \f)$ is a cascade compactification.
Since $\A$ contains a sufficient set of Lyapunov functions the compactification is almost dynamic
by Theorem \ref{th1.20}.
  A continuous map is a +proper relation and so the compactification is dynamic by  Corollary \ref{cor1.10new}.

 If $u \in \B(X)$ has compact support then $\hat u \circ \f^n$ vanishes on $\X \setminus \tilde X$ (see Exercise
 \ref{ex1.8new}). Hence, as in the proof of Lemma \ref{lem1.17} distinct points of $\X \setminus \tilde X$ are
 distinguished by some $L \in \L$ and so cannot lie in a common $\G \f \cap \G \f^{-1}$ equivalence class $\hat E$.
 If $z \in \hat E \cap (\X \setminus \tilde X) $, then by Lemma \ref{lem2.2}
     $\f(z) \in \hat E \cap (\X \setminus \tilde X)$.  It follows that  $z = \f(z)$, i.e. $z \in |\f|$.

(b) now follows from Proposition \ref{prop1.12new}.

(c),(d): If $f$ is proper then $\tilde X = X$ by Proposition \ref{prop1.11new}
the set of functions with compact support is $f^*$ +invariant and so $\A$ is
generated by $\L$ and the functions with compact support.  If $f$ is a homeomorphism
then the set of functions with compact support is $f^*$ invariant and so $\A$ is as well.

(e): If $X$ is metrizable then we can choose a countable sufficient set $\L_0$ and then extend to a
countable $f^*$ +invariant sufficient set or, in the reversible case, a countable $f^*$ invariant sufficient set.

$\Box$ \vspace{.5cm}

In the +proper case we can sharpen Theorem \ref{th1.12}.

\begin{theo}\label{cor3.3} Let $f$ be a +proper relation on  $X$.

(a) If $A$ is a compact subset of $X$ which is $\mathcal{G}f \
+$invariant and $U$ is a bounded open subset of $X$ with $A
\subset  U$, then there exists a Lyapunov function $L : X \to
[0,1]$  such that  $A \subset L^{-1}(1)$ and with
 support in the closure of $U$. In particular, $L$ has compact support.

(b) Let $(\hat{X},\hat{f})$ be a proper  compactification for
 $(X,f)$ (N.B. we do not assume that the compactification is dynamic).
If for some $x \in X$, the set $\mathcal{G}f(x)$ is
compact then
\begin{equation}\label{eq3.6}
\mathcal{G}\hat{f}(x) \quad = \quad \mathcal{G}f(x).
\end{equation}
\end{theo}

{\bfseries Proof:}  (a): Let $(\hat{X},\hat{f})$ is a
Lyapunov function compactification for $(X,f)$.
By Theorem \ref{th1.9new} (b)  $\mathcal{G}\hat{f}(A) =
\mathcal{G}f(A) \subset A$. Let $B = (1_{\hat{X}} \cup \mathcal{G}
\hat{f}^{-1})(\hat{X} \setminus U)$. By Corollary \ref{cor1.10} applied to
$\mathcal{G} \hat{f}$ there exists a Lyapunov function $\hat{L} :
\hat{X} \to [0,1]$ which is one on $A$ and zero on $B$ and so has
support contained in the closure of $U$. Let $L$ be the
restriction of $\hat{L}$ to $X$.

(b): Assume $y \not\in \mathcal{G} f(x) $. We use (a) to choose a Lyapunov
function $L: X \to [0,1]$ with compact support such that $\{ x \}
\cup \mathcal{G}f(x) \subset  L^{-1}(1)$ and such that $L(y) = 0$. For the
compactification $(\X,\f)$ extend $L$ to $\hat{L}$
on $\hat{X}$ by letting $\hat{L}$ be zero on $\hat{X} \setminus
X$.

If $\hat L(z) > 0$ then $z \in X$ and so by Proposition \ref{prop1.8new}
$\f(z) \subset X$ and so $\f(z) = f(z)$.  Hence, $\hat{L}(y) = L(y) \geq L(z) = \hat L(z)$
for all $y \in \f(z)$.  That is,  $\hat{L}$ is a $\hat{f}$ Lyapunov function. Since
$\hat{L}(x) = 1$, $\hat{L} \geq 1$ on $\mathcal{G}\hat{f}(x)$. But
$\hat{L} = 0$ on $\{y \} \cup \hat{X} \setminus X$. As $y$ was arbitrary (\ref{eq3.6}) follows.

$\Box$ \vspace{.5cm}

{\bfseries Remark:} We can construct the Lyapunov function $L$ so
that  $A = L^{-1}(1)$, iff the compact +invariant set $A$ is
a $G_{\delta}$ subset of $X$.  In particular, such an $L$ always
exists when $X$ is metrizable. \vspace{.5cm}

\begin{ex}\label{ex3.4} The analogue of Theorem \ref{cor3.3}  for $\mathcal{G}f(x) \cap
\mathcal{G}f^{-1}(x)$ does not hold. \end{ex}
With $X = {\mathbb R}$ let $f(x) = e^x
- 1$. $|\mathcal{G}f|$ contains  just the fixed point $0$. If we
let $\hat{X}$ be the one point compactification then
$\mathcal{G}\hat{f} = \hat{X} \times \hat{X}$ and the entire space
is a single $\mathcal{G}\hat{f} \cap \mathcal{G}\hat{f}^{-1}$
equivalence class.

$\Box$ \vspace{.5cm}

\begin{cor}\label{cor3.5} Let $f$ be a +proper relation on  $X$. Let $(X_*,\hat{f})$ be the proper
compactification of $(X,f)$ with $X_*$ the one point
compactification.  If $\mathcal{G}f(x)$ is compact for every $x
\in X$, then $(X_*,\hat{f})$ is a Lyapunov function
compactification.

If $f$ is a proper continuous map on $X$ and  $\mathcal{G}f(x)$ is compact for every $x
\in X$, then $(X_*,\hat{f})$ is a Lyapunov function cascade
compactification.
\end{cor}

{\bfseries Proof:} By Theorem \ref{cor3.3} the Lyapunov functions with
compact support constitute a sufficient set of Lyapunov functions.
The associated algebra is just $\mathcal{A}_0$, generated by the functions
of compact support and associated with the one point
compactification. Recall that if a continuous map $f$ is proper then $\A_0$ is
$f^*$ +invariant and so the one-point compactification is a cascade compactification.

$\Box$ \vspace{.5cm}

If we have two closed relations $f_1 \subset f_2$ on $X$ then any
Lyapunov function for $f_2$ is a Lyapunov function for $f_1$.
Hence, if $\mathcal{L}_i$ is a sufficient set of Lyapunov
functions for $f_i \  (i = 1,2)$ then $\mathcal{L} = \mathcal{L}_1
\cup \mathcal{L}_2$ is a sufficient set of Lyapunov functions for
$f_1$ which contains a sufficient set of Lyapunov functions for
$f_2$. Notice that if $X$ is metrizable we can choose the
$\mathcal{L}_i$'s to be countable and so obtain $\mathcal{L}$
which is countable.

\begin{cor}\label{cor1.21} Let $f_1 \subset f_2$ be closed relations on a space
$X$. Let $\mathcal{L}$ be a sufficient set of $f_1$ Lyapunov
functions which contains a sufficient set of $f_2$ Lyapunov
functions. The $\mathcal{L}$ compactification of $X$  is almost dynamic
for $f_2$ and it is dynamic for $f_2$ if $\L$ satisfies the splitting property for $f_2$.
\end{cor}

{\bfseries Proof:} This is immediate from Theorem \ref{th1.20}.

$\Box$ \vspace{.5cm}

Our main application of this corollary is  the following special
case.

\begin{cor}\label{cor1.22} Let $f$ be a closed relation on a space $X$. There
exist a dynamic Lyapunov compactification $(\hat{X},\hat{f})$ of $(X,f)$
such that $(\hat{X},\hat{f}\cup 1_{\hat{X}} \cup \hat{f}^{-1})$ is
a dynamic compactification of $(X, f \cup 1_X \cup f^{-1})$.
\end{cor}

{\bfseries Proof:} Choose a sufficient set of Lyapunov functions for $f$ which satisfies the splitting
property for $f_1 = f$ and which contains a sufficient set for $f_2 = f \cup 1_X \cup f^{-1}$.
Since the latter is reflexive the splitting property holds vacuously for $f_2$.
Clearly, the closure of $f^{-1}$ is $\hat{f}^{-1}$.

$\Box$ \vspace{.5cm}

Clearly, $\mathcal{G}(f \cup 1_X \cup f^{-1}) $ is the smallest
closed equivalence relation which contains $f$. Observe that if
$f$ is a map then $1_X \subset f^{-1} \circ f$ and so

\begin{equation}\label{eq1.28}
\begin{split}
f \ \mbox{a continuous map} \qquad \Lra \hspace{4cm} \\
\mathcal{O}(f \cup f^{-1}) = \mathcal{O}(f \cup 1_X \cup f^{-1})
; \hspace{2cm}\\
f \ \mbox{a
homeomorphism} \qquad \Lra \hspace{4cm} \\
\mathcal{O}f \cup 1_X \cup \mathcal{O} (f^{-1}) = \mathcal{O}(f
\cup 1_X \cup f^{-1}). \hspace{1cm}
\end{split}
\end{equation}
It  follows by taking closures that the same results hold with
$\mathcal{O}$ replaced by $\mathcal{R}$ or by $\mathcal{N}$
throughout. Furthermore, when $f$ is a continuous map
$\mathcal{G}(f \cup 1_X \cup f^{-1}) \  = \ \mathcal{G}(f \cup
f^{-1})$. However, even when $f$ is a homeomorphism, $\mathcal{G}f
\cup 1_X \cup \mathcal{G}f^{-1}$ is usually not transitive and so
is usually a proper subset of $\mathcal{G}(f \cup 1_X \cup
f^{-1})$. \vspace{.5cm}

\begin{ex}\label{ex1.23}A Lyapunov compactification of $f$ need not be
a dynamic compactification for $f \cup f^{-1}$\end{ex}
Let $Y$ be a compact space and $T :
{\mathbb Z} \to {\mathbb Z}$ be the translation map $T(t) = t+1$.
Let $f = 1_Y \times T$ on $X = Y \times {\mathbb Z}$. Clearly, the
two point compactification $\hat{X} = X \cup \{ \pm \infty \}$
yields a Lyapunov compactification of $f$. But
$\mathcal{G}(\hat{f} \cup \hat{f}^{-1}) = \hat{X} \times \hat{X}$
and so $(\hat{X},\hat{f} \cup \hat{f}^{-1})$ is not dynamic when
$Y$ contains more than one point.

$\Box$ \vspace{.5cm}

The set $\mathcal{L}_f$ of all bounded Lyapunov functions for $f$
is always a sufficient set of Lyapunov functions with the splitting property.
We will denote by $(\beta_f X, \hat{f})$ the
$\mathcal{L}_f$ compactification which is the \emph{maximal
Lyapunov function compactification}.\index{maximal Lyapunov function compactification}%
\index{compactification!maximal Lyapunov function compactification} The space $\beta_f X$ might
be rather large. For example, if $f = 1_X$ then every function in
$\mathcal{B}(X)$ is a Lyapunov function and the compactification
$\beta_{1_X} X$ is the Stone-\v{C}ech compactification.  Often we
wish to remain in the category of metrizable spaces and so
restrict ourselves to countably generated subalgebras
$\mathcal{A}$. However, we have to go to the maximal
compactification to obtain the following:

\begin{theo}\label{th1.24} Let $(\beta_f X, \hat{f})$ be the maximal Lyapunov
function compactification of $(X,f)$. If $A, B$ are disjoint
closed subsets of $X$ such that $\mathcal{G} f(A) \subset A$ and
$\mathcal{G} f^{-1}(B) \subset B$ then $A$ and $B$ have disjoint
closures in $\beta_f X$.
\end{theo}

{\bfseries Proof:} By Corollary \ref{cor1.10} there is a Lyapunov function
$L : X \to [0,1]$ which is zero on $B$ and one on $A$. Since $L$
extends to a continuous $\hat{L}$ on $\beta_f X$ we have that the
closure of $B$ is contained in $\hat{L}^{-1}(0)$ and the closure
of $A$ is in $\hat{L}^{-1}(1)$.

$\Box$ \vspace{.5cm}

{\bfseries Remark:} If $f_1$ is any closed relation on $X$ which
contains $f$ then the bounded Lyapunov functions for $f_1$ are all
included in $\mathcal{L}_f$. From Corollary \ref{cor1.21} it follows that
the compactification $X \subset \beta_f X$ is dynamic for every
closed relation $f_1$ with $f \subset f_1$. \vspace{.5cm}

We now consider when the Lyapunov functions suffice to determine the topology.

For a space $X$ let  $\L \subset \B(X)$. We say that $\L$ \emph{distinguishes points} \index{distinguishes points}
if whenever $x_1,x_2$ are distinct points of $X$ there exists $L \in \L$ such that $L(x_1) \not= L(x_2)$.
We say that $\L$ \emph{distinguishes points and closed sets} \index{distinguishes points and closed sets}
if whenever $A$ is a closed subset of $X$ and
$x \in X \setminus A$ there exists $L \in \L$ such that $L(x)$ is not in the closure of the image $L(A)$.
We say that $\L$  \emph{determines the topology
of} \index{determines the topology} $X$ if whenever $\{x_i \}$ is a net in $X$ and $x \in X$, such that $\{ L(x_i) \}$
converges to $L(x)$ for all $L \in \L$ then $\{ x_i \}$ converges to $x$.  This is equivalent to saying that the
topology on $X$ is the weak topology induced by $\L$, i.e. the coarsest topology with respect to which 
all the functions in $\L$ are continuous. 

For $u \in \B(X)$
we denote by $I_u$ the smallest closed subinterval of $\R$ which contains the image of $u$.
To $\L \subset \B(X)$ we associate the function $j_{\L} : X \to \prod_{L \in \L} \ I_L$ by
$(j_L(x))_L = L(x)$.

We recall some standard results.

\begin{prop}\label{prop1.25new} For $\L \subset \B(X)$ let $\B$ be the closed subalgebra generated
by $\L$.
\begin{enumerate}
\item[(a)] The following are equivalent.
\begin{itemize}
\item[(i)] $\L$ determines the topology of $X$.
\item[(ii)] $\B$ determines the topology of $X$.
\item[(iii)] $j_{\L}$ is an embedding of $X$ into $\prod_{L \in \L} \ I_L$.
\item[(iv)] The compactification $j : X \to \X$ associated with $\B$ is a proper compactification.
\item[(v)] $\B$ distinguishes points and closed sets.
\end{itemize}
\item[(b)] If $X$ is compact then the following are also equivalent.
\begin{itemize}
\item[(vi)] $\L$ distinguishes points.
\item[(vii)] $\B = \B(X)$.
\end{itemize}
\item[(c)] If $Y$ is an open subset of $X$ and $\L$ determines the topology of $X$ then
$\{ L|Y : L \in \L \} \subset \B(Y)$ determines the topology of $Y$.
\end{enumerate}
\end{prop}

{\bfseries Proof:} (a), (b): \quad (i) $\Leftrightarrow$ (ii): The set of $L$ such that $\{L(x_i)\}$ converges to $L(x)$
is closed under algebraic operations and uniform limits. So if it holds for all $L \in \L$ then it
holds for all $L \in \B$. Thus, (ii) implies (i) and the converse is obvious.

(i) $\Rightarrow$ (vi): If $x_1, x_2$ are distinct points of $X$ then the sequence which
is constant at $x_1$ does not converge to $x_2$ and so assuming (i) it cannot be that $L(x_1) = L(x_2)$ for all
$L \in \L$.

(vi) $ \&$ compactness $\Rightarrow$ (vii): This is the Stone-Weierstrass Theorem.

(vii) $\Rightarrow$ (ii): $X$ is assumed to be completely regular.

(i) $\Leftrightarrow$ (iii): To say that $\{ L(x_i) \}$ converges to $L(x)$ for all $L \in \L$
exactly says that $\{ j_{\L}(x_i) \}$ converges to $j_{\L}(x)$.

(iii) $\Rightarrow$ (iv): Let $\bar X $ denote the closure in $\prod_{L \in \L} \ I_L$ of $j_{\L}(X)$.
Thus, $j_{\L} : X \to \bar X$ is a compactification of $X$. I claim it is the compactification
associated with $\B$, i.e. $j_{\L}^*\B(\bar X) = \B$. Let $\bar L $ denote the projection of to the $L$
coordinate. Clearly, $\bar L$ on $\bar X$ extends $L$ on $X$ and $\bar{\L} = \{ \bar L : L \in \L \}$ distinguishes
points of $\bar X$.  By the Stone-Weierstrass Theorem $\bar{\L}$ generates $\B(\bar X)$. By (iii) $j_{\L}$ is
an embedding and so $j_{\L}^*$ an injective isometry on $\B(\bar X)$. It follows that $\L =
j_{\L}^*(\bar{\L})$ generates $j_{\L}^*\B(\bar X)$ and so the latter equals $\B$.

(iv) $\Rightarrow$ (v): If $A \subset X$ is closed and $x \not\in A$ then with $\hat A$ the closure of
$A$ in $\hat X$, we have $x \not\in \hat A$ because the compactification is proper. There exists
$u \in \B(\hat X)$ such that $u(x) = 1$ and $u = 0$ on $\hat A$. The restriction of $u$ to $X$ is an
element of $\B$ which distinguishes $x$ from $A$.

(v) $\Rightarrow$ (ii): Assume $\{ u(x_i) \}$ converges to $u(x)$ for all $u \in \B$. Let $U$ be an open
neighborhood of $x$ and $A = X \setminus U$. Since $\B$ distinguishes points and closed sets there exists
$u \in \B$ such that $u(x)$ is not in the closure of $u(A)$.  It follows that the net $\{ u(x_i) \}$ is
eventually in the complement of $u(A)$ and so $\{ x_i \}$ is eventually in $U$. Since $U$ was arbitrary,
$\{ x_i \}$ converges to $x$.

(c):  Obvious from the definition of the relative topology.

$\Box$ \vspace{.5cm}

Since all Lyapunov functions are constant on the $\G f \cap \G f^{-1}$
equivalence classes, the Lyapunov functions distinguish points exactly when
$\mathcal{G}f \cap \mathcal{G}f^{-1} \subset 1_X$ i.e. when the
reflexive, transitive relation $1_X \cup \mathcal{G}f$ is
anti-symmetric. Because of case (iv) it can happen that a closed,
reflexive, anti-symmetric, transitive relation on a locally
compact space may not have an anti-symmetric compactification, see Exercise \ref{ex1.25}.
In such a case,  the Lyapunov functions,  while they distinguish points of $X$, do
not suffice to yield the correct topology on $X$.

\begin{theo}\label{th1.26} Assume $(\X, \f)$ is a dynamic Lyapunov function compactification of a dynamical
system $(X,f)$.  Let $\hat{\L}$ denote the set of Lyapunov functions for $(\X, \f)$ and let $\L$ denote
restrictions to $X$ of the elements of $\hat{\L}$, i.e. $\L$ consists of the Lyapunov functions for $(\X, \f)$
which extend to $\X$.
\begin{enumerate}
\item[(a)]$\L$ determines the topology of $X$ iff $\G \f \cap \G \f^{-1} \subset 1_{\X}$.

\item[(b)]If  $\G f$ is asymmetric, i.e. $\G f \cap \G f^{-1} = \emptyset$, or $f$ is +proper and
 $\G f \cap \G f^{-1} \subset 1_X$, then $\G \f \cap \G \f^{-1} \subset 1_{\X}$.
\end{enumerate}
\end{theo}

{\bfseries Proof:} (a): If $\G \f \cap \G \f^{-1} \subset 1_{\X}$ then $\hat{\L}$ distinguishes points of $\X$
and so by Proposition \ref{prop1.25new} (b) $\hat{ \L}$ determines the topology of $\X$. By
\ref{prop1.25new} (c) $\L$ determines the topology of $X$ since the compactification is proper.

Conversely, if $\L$ determines the topology of $X$ then  $\L$ distinguishes points of $X$ and so
$\G f \cap \G f^{-1} \subset 1_X$, i.e. the $\G f \cap \G f^{-1}$ equivalence classes in $ |\G f|$
are singletons. This implies that in Theorem \ref{th1.18} (b) case (ii) does not occur. In either case (i)
case (iii) it follows that $\hat E$ is a singleton. Thus, $\G \f \cap \G \f^{-1} \subset 1_{\X}$ fails iff
case (iv) occurs and so there exist $x \in X$ and $z \in \hat X \setminus X$ such that $x, z \in \hat E$ and
so $\hat L(x) = \hat L(z)$ for every $\hat L \in \hat{ \L}$.  Let $\{ x_i \}$ be a net in $X$ which converges to
$z $ and so $\{ L(x_i) = \hat L(x_i)  \}$ converges to $\hat L(z) = L(x)$ for every $L \in \L$.
Since $\{ x_i \}$ does not converge to $x, \ \L$ does not determine the topology of $X$.

(b): Assume that  $\G f \cap \G f^{-1} \subset 1_{X}$.  To show that  $\G \f \cap \G \f^{-1} \subset 1_{\X}$
it suffices, as above, to show that case (iv) does not occur. If $f$ is +proper then case (iv) does not occur by
Theorem \ref{th1.18} (e).

If there exist $x \in X$ and $z \in \hat X \setminus X$ which
are $ \G \f \cap \G \f^{-1}$ equivalent then $x \in |\G \f | \cap X $ which equals $| \G f|$ since the
compactification is dynamic. Hence, $\G f$ is not asymmetric.  Contrapositively, $\G f$ asymmetric
implies that case (iv) does not occur.

$\Box$

{\bfseries Remark:} Even when $\G f$ is asymmetric it is usually not true that $\G \hat{f}$ is asymmetric.
In fact, if $f$ is a
closed relation on a compact space then $|\G f| = \emptyset$  iff there exists a positive integer
$n$ such that the $n$-fold composition $f^n$ is the empty relation, i.e. $f^n = \emptyset$ ( See
Akin (1993) Exercise 2.16).

\vspace{.5cm}

Notice what the theorem does \underline{not} say.

\begin{ex}\label{ex1.27}A sufficient family of Lyapunov functions need not determine the topology.\end{ex}
Let $A_0 = \Z_+$ and $A $ be the one point compactification of $A_0 $ with the
point at infinity identified with $0$. Thus, the ``identity map" $\alpha : A_0 \to A$ is a continuous
bijection which is not a homeomorphism.  Let $B_0 = \Z$ and $B$ be the two point compactification
$B_0 \cup \{ -\infty, +\infty \}$. The inclusion $\beta: B_0 \to B$ is an embedding. Let $T$ be the
homeomorphism on $B_0$ given by $n \mapsto n+1$ with $\hat{T}$ the extension to a homeomorphism of $B$.
Let $X = A_0 \times B_0$ and $\bar X = A \times B$.
The injection $j = \alpha \times \beta : X \to \bar X$ is a compactification which is not proper.
Let $f = 1_{A_0} \times T$ so that $\bar f = 1_A \times \hat{T}$. Let $\bar{\L}$ be the set of all
Lyapunov functions for the homeomorphism $\bar f$ on $\bar X$ and $\L = \{ L = \hat{L} \circ j : \bar{L} \in \bar{\L} \}$.
It is easy to check the $\L$ is a sufficient family for $(X,f)$, but the topology induced by $\L$ is the
one pulled back from $\bar X$ via $j$ rather than the original, discrete one.

The $\L$ compactification
$(\X, \f)$ is the proper compactification associated with the algebra $\A$ generated by $\L$ together
with all the functions of compact support.
By Theorem \ref{th1.26} the set of the restrictions to $X$ of the Lyapunov functions for $(\X, \f)$ does determine
the topology of $X$. This set includes $\L$ but there are additional Lyapunov functions constructed using
functions of compact support which do not factor through $j$.

$\Box$ \vspace{1cm}

We conclude by considering  in detail how case (iv) of Theorem \ref{th1.6} occurs and how it
can be avoided.  As  illustrated by Example \ref{ex1.25}, sometimes it cannot be avoided.

\begin{theo}\label{th1.28}  Let $F$ be a closed, transitive relation on $X$.   Let $E$ be a
compact $F \cap F^{-1}$ equivalence class in $|F|$ of $X$. Assume that for any neighborhood $U$ of $E$,
$F(U) \cap F^{-1}(U)$ is unbounded.  If $(\hat{X},\hat{F})$ is any proper compactification of $(X,F)$, there
exists $z \in \hat{X} \setminus X$  and $a,b \in E$ such that $(a,z), (z,b) \in \hat{F}$. In particular,
$z$ is an element of the $\mathcal{G}\hat{F} \cap \mathcal{G}\hat{F}^{-1}$ class which contains $E$.
\end{theo}

{\bfseries Proof:}   As $U$ varies over compact neighborhoods of $E$ and $W$ varies over all cobounded
closed subsets of $X$ (i.e. $X \setminus W$ is bounded),
$\{ W \cap F(U) \cap F^{-1}(U) \}$ is a filterbase of closed subsets of $X$ with an
empty intersection in $X$.  Let $z$ be any point of $\hat{X}$ in the intersection of the $\hat{X}$ closures of
these sets. As $N$ varies over closed neighborhoods of $z$ let $i$ vary over the sets
$N \cap F(U) \cap F^{-1}(U) \subset X $  For each $i$ there is a $z_i \in N \cap X$ and $a_i, b_i \in U$ such
that $(a_i,z_i) \in F$ and $(b_i,z_i) \in F$.  In the limit, $z_i$ tends to $z$ and limit points $a, b$ of the
$a_i, b_i$ nets lie in $E$.  Since $\hat{F}$ is the closure of $F$, $(a,z), (z,b) \in \hat{F}$.

$\Box$ \vspace{.5cm}

Let $f$ be a closed  relation on $X$. If $A$ is any subset of $X$ then we will denote by $[[A]]_f$
the smallest, closed $f \ +$invariant subset of $X$
which contains $A$.  When the closed relation is understood we will just write $[[A]]$.
If $f$ is reflexive then $A$ is $+$invariant iff it is invariant.\index{$[[A]]_f$}

For example, if $A$ is compact and $f = F$ is transitive then since $F(A)$ is closed,
$[[A]]_F = A \cup F(A)$. If $A$ is a singleton $\{x\}$ and $f$ is a continuous map then
$[[\{x\}]]_f = \{ x \} \cup \mathcal{R} f(x)$.

In general,
$[[A]]$  may require a transfinite construction. Let $f_1 = f \cup 1_X$, i.e. make $f$ reflexive.
\begin{equation}\label{eq1.31}
\begin{split}
K_0 \qquad =_{def} \qquad A, \hspace{2cm}\\
K_{\alpha} \quad =_{def} \qquad \overline{f_1(\bigcup_{\beta < \alpha} K_{\beta})}.
\end{split}
\end{equation}
The increasing transfinite sequence of sets (which are closed once $\alpha > 0$) stabilizes at
$[[A]]$.

\begin{prop}\label{prop1.29}  Let $f$ be a closed  relation on $X$. If $A \subset X$ is closed then
\begin{equation}\label{eq1.32}
[[A]]_f  \quad = \quad  A \cup [[f(A)]]_f
\end{equation}
If the relation $f$ is $+$proper then
\begin{equation}\label{eq1.33}
[[A]]_{\G f} \quad = \quad A \cup [[ \G f(A)]]_{\G f}  \quad = \quad A \cup [[f(A)]]_{\G f}.
\end{equation}
\end{prop}

{\bfseries Proof:}   Since $f(A) \subset [[f(A)]]_f$ and the latter is $f \ +$invariant it follows that
$A \cup [[f(A)]]_f$ is $f \ +$ invariant and so contains $[[A]]_f$. The reverse inclusion is obvious.

Clearly,
\begin{equation}\label{eq1.34}
[[A]]_{\G f} \ \supset \ A \cup [[ \G f(A)]]_{\G f}  \ \supset \ A \cup [[f(A)]]_{\G f}.
\end{equation}
If $f$ is $+$proper then by (\ref{eq1.7}) $\G f(A) =
f(A) \cup \G f(f(A))$ and so is contained in $[[f(A)]]_{\G f}$.  It follows that $A \cup [[f(A)]]_{\G f}$ is
$\G f \ +$invariant and so contains $[[A]]_{\G f}$.

$\Box$ \vspace{.5cm}

\begin{prop}\label{prop1.30} Let $\L$ be the set of all bounded Lyapunov functions for a closed relation $f$.
If $A \subset X$ then
\begin{equation}\label{eq1.35}
[[A]]_{\G f} \quad = \quad \{ x : L(x) \geq inf \ L|A \ \ \mbox{for all}
 \ \ L  \in  \L \  \}.
\end{equation}
\end{prop}

{\bfseries Proof:} The set on the right is closed, $\G f \ +$ invariant and contains $A$.  Hence it
contains $[[A]]_{\G f}$.

On the other hand, if $x \not\in [[A]]_{\G f}$ then $B = \{ x \} \cup
\G f^{-1}(x)$ is a closed, $\G f^{-1} \ +$ invariant subset of $X$ which
is disjoint from $[[A]]_{\G f}$.  By Corollary \ref{cor1.10} there exists a Lyapunov
function $L$ such that $L(x) = 0$ and $L$ on $ [[A]]_{\G f}$ is constant at 1.

$\Box$ \vspace{.5cm}

The remainder of this
section is rather technical and is not used later. The purpose is to prove the following sharpening of
Theorem \ref{th1.26}.

\begin{theo}\label{cor1.35} Assume that a closed, transitive relation $F$ is  reflexive and anti-symmetric,
i.e. $F \cap F^{-1} = 1_X$. The following conditions are equivalent:
\begin{itemize}
\item[(i)] For every $x \in X$ there exists a compact neighborhood $U$ of $x$ such that
$F(U)\cap F^{-1}(U)$ is compact.
\item[(ii)] Every $x \in X$ has a compact unrevisited neighborhood.
\item[(iii)] There exists a dynamic compactification $(\X, \hat{F})$ for $(X, F)$ such that
$\G \hat{F}$ is reflexive and anti-symmetric on $\X$.
\item[(iv)] The Lyapunov functions for $F$ determine the topology of $X$.
\end{itemize}
\end{theo}
\vspace{.5cm}

\begin{lem}\label{lem1.31}  Let $f$ be a closed, reflexive relation on
$X$, $B,G $ be  disjoint subsets of $X$ with $G$ open.
\begin{enumerate}
\item[(a)]   If $G = f^{-1}(G)$ then $G$ is disjoint from $[[B]]_f$.
\item[(b)]  If there exists a closed subset $Q$ of $X$ such that
\begin{equation}\label{eq1.36}
\begin{split}
f(Q) \quad = \quad Q \hspace{2.5cm}\\
B \quad \subset \quad Q \hspace{3cm}\\
f^{-1}(G) \cap Q \quad \subset \quad G \hspace{2cm} \\
\end{split}
\end{equation}
then $f^{-1}(G)$ contains $G$ and  is disjoint from $[[B]]_f$.
\end{enumerate}
\end{lem}

{\bfseries Proof:}  (a)  Since $G$ is open and $f^{-1}(G)$ invariant, $X \setminus G$ is
closed and $f \ +$invariant.  Because  $B \subset X \setminus G, [[B]]_f \subset X \setminus G. $

(b) Because $f$ is reflexive, $f^{-1}(G)$ contains $G$.

 We use the transfinite construction of (\ref{eq1.31}), $\{ K_{\alpha} \}$, for $[[B]]$ beginning with $K_0 = B$ and
we show inductively that $f^{-1}(G) \cap K_{\alpha}$ is empty for all $\alpha$.

Notice first that since $Q$ is closed and $f$ invariant,  $B \subset Q$ implies $K_{\alpha} \subset [[B]] \subset Q$ .

For $\alpha = 0$, $B \subset Q$ implies
$f^{-1}(G) \cap B \subset f^{-1}(G) \cap Q \cap B \subset G \cap B = \emptyset$.

For the inductive step, recall that for any subsets $E, F $ of $X$
\begin{equation}\label{eq1.37}
E \cap f^{-1}(F) = \emptyset \qquad \Longleftrightarrow \qquad f(E) \cap F = \emptyset
\end{equation}
because each is equivalent to $f \cap E \times F = \emptyset$.

By inductive hypothesis, $f^{-1}(G) \cap K_{\beta} = \emptyset$ for all $\beta < \alpha$ and
so $G \cap f(K_{\beta}) = \emptyset$.  Because $G$ is open,
$G \cap K_{\alpha} = \overline{\bigcup_{\beta < \alpha} f(K_{\beta})} = \emptyset$.

Finally, as in the $\alpha = 0$ case
$f^{-1}(G) \cap K_{\alpha} \subset f^{-1}(G) \cap Q \cap K_{\alpha} \subset G \cap K_{\alpha} = \emptyset$.

$\Box$ \vspace{.5cm}

{\bfseries Remark:} If $\{ A_n \}$ is a sequence of subsets such that for all $n$
\begin{equation}\label{eq1.38}
\begin{split}
A_{n} \quad \subset \quad (A_{n+1})^{\circ} \hspace{3cm}\\
f^{-1}(A_{n}) \cap Q \quad \subset \quad A_n \hspace{3cm} \\
\end{split}
\end{equation}
then $G = \bigcup_n A_n = \bigcup_n (A_n)^{\circ}$ is an open set such that $f^{-1}(G) \cap Q  \subset  G$.
\vspace{.5cm}

Recall that for a relation $f$ on $X$ we call  $A \subset X$ an $f$
\emph{unrevisited subset}  when
\begin{equation}\label{eq1.39}
\mathcal{O}f(A) \cap \mathcal{O}f^{-1}(A) \quad \subset \quad  A.
\end{equation}
That is, if there exist $n, m > 0$ such that $f^n(x)$ and $f^{-m}(x)$ both meet $A$, then $x \in A$.

\begin{lem}\label{lem1.32} Assume $F$ is a reflexive, transitive relation on $X$ and $A \subset X$.  The set
$B = F(A) \cap F^{-1}(A)$ is the smallest unrevisited set which contains $A$.   In particular,
$A$ is unrevisited iff $A = F(A) \cap F^{-1}(A)$.
\end{lem}

{\bfseries Proof:}  Since $F$ is transitive, $\mathcal{O}F = F$ and since $1_X \subset F$, i.e. $F$
is reflexive, $A \subset F(A)\cap F^{-1}(A) = B$. Clearly,  $A$ is unrevisited iff $A = B$.

In any case, $A \subset B \subset F(A), F^{-1}(A)$ and so
$B = F(A) \cap F^{-1}(A) \subset F(B) \cap F^{-1}(B) \subset FF(A) \cap F^{-1}F^{-1}(A) \subset F(A) \cap F^{-1}(A) = B$.
So $B$ is always unrevisited and if $A \subset C$ then $B \subset F(C) \cap F^{-1}(C)$.

$\Box$ \vspace{.5cm}

Notice that for a point $x \in X$ the unrevisited set $F(x) \cap F^{-1}(x)$ is the $F \cap F^{-1}$ equivalence
class of $x$.  Any $F$ unrevisited set is saturated by the $F \cap F^{-1}$ relation and $F(A) \cap F^{-1}(A)$
consists of all points which lie ``between" points of $A$.

\begin{lem}\label{lem1.33}  Let $F$ be a closed, transitive, reflexive relation on $X$.   Let $A$ be a
compact unrevisited subset of $X$, so that $A = F(A) \cap F^{-1}(A)$. Assume that $A$ admits a
compact neighborhood $U$ such that the unrevisited set $F(U) \cap F^{-1}(U)$ is compact as well.
Let $W$ be an open  set which contains $A$ with $W \subset \subset U$.
\begin{enumerate}
\item[(i)]There exists a compact set $C$ such that
\begin{equation}\label{eq1.40}
\begin{split}
F(A) \cap \overline{W} \quad \subset \subset \quad C \quad  \subset \subset U, \\
 F^{-1}(A) \cap F(C)  \quad \subset \subset \quad W. \hspace{1cm}
 \end{split}
 \end{equation}
 \item[(ii)] If $C$ is any compact set which satisfies the conditions of (i) then
 the closed $F$ invariant set $B = [[(X \setminus W) \cap F(C)]]_F$
is contained in $F(C)$ and is disjoint from $F^{-1}(A)$.
 There exists $L : X \to [0,1]$
a Lyapunov function $L : X \to [0,1]$ for $F$ such that $L(x) = 0 $ for $x \in F^{-1}(A)$ and $L(x) = 1$ for
$x \in F(C) \setminus W$.
\item[(iii)] Let $(\hat{X},\hat{F})$ be a dynamic extension of $(X,F)$ such that the function $L$ extends to a
continuous function $\hat{L}$ on $\hat{X}$.  The function $\hat{L}$ is a Lyapunov function for
$\mathcal{G}\hat{F}$ and if $z \in (\hat{X} \setminus X) \cap \mathcal{G}\hat{F}(x)$ for some $x \in  A$
then $\hat{L}(z) = 1$.
\end{enumerate}
\end{lem}

{\bfseries Proof:} (i) Since $F(A) \cap \overline{W}  \subset \subset U$
there exist compact sets $C$ such that
$F(A) \cap \overline{W}  \subset \subset  C  \subset \subset U$.

Observe that $A \subset F(A) \cap \overline{W} \subset F(A)$ and
$F(A) = F(F(A))$ by transitivity. Hence, $F(A) = F(F(A) \cap \overline{W})$.  For every compact $C \subset U$,
the closed set $F^{-1}(A) \cap F(C)$ is contained in
the compact set $F(U) \cap F^{-1}(U)$ and so is compact as well. Letting $C$ vary, the intersection of the
compacta
\begin{equation}\label{eq1.41}
\{F^{-1}(A) \cap F(C): F(A) \cap \overline{W} \subset \subset C = \overline{C} \subset \subset U \}
\end{equation}
is $ F^{-1}(A) \cap F(F(A) \cap \overline{W}) = F^{-1}(A) \cap F(A)$ which is contained  in the open
set $W$. Hence,  we can choose $C$ small enough to satisfy the second condition of (\ref{eq1.40}) as well.

 (ii) Because $C$ is compact, $F(C)$ is closed.  Because $F$ is reflexive and transitive $F(C)$ contains $C$ and
 is $F$ invariant.  Hence, $F(C) = [[F(C)]] $. Since $F(C)$ contains
 $F(C) \cap (X \setminus W)$, $ F(C) = [[F(C)]]$ contains $B $.

For the disjointness result we will apply Lemma \ref{lem1.31} (b) with $Q = F(C)$.

Notice that $F$ reflexive and $A \subset W$ implies $A \subset F(A) \cap W$ and so $C$ is a neighborhood of $A$.

Assume that $K$ is compact and $K' = K \cup (F^{-1}(K) \cap F(C)) \subset W$.  Then $K \subset K' \subset F^{-1}(K)$
implies $F^{-1}(K) \subset F^{-1}(K') \subset F^{-1}(F^{-1}(K)) = F^{-1}(K)$.  Hence, $F^{-1}(K) = F^{-1}(K')$ and
$K' = K' \cup (F^{-1}(K')\cap F(C))$.  If $N \subset W$ is a compact neighborhood of $K'$ then $F^{-1}(N) \cap F(C)$
is a closed subset of the compact set $F^{-1}(U) \cap F(U)$ and so is compact. As $N$ decreases to $K'$ the
sets $N' = N \cup (F^{-1}(N) \cap F(C))$ have intersection $K'$.  Hence, there exists a compact neighborhood $N$ of
$K'$ such that $ N' \subset W$.

Let $A_0 = A \cup (F^{-1}(A) \cap F(C))$. Choose $N $ a compact neighborhood of $A_0$  in $W$
so that the compact set $A_1 = N'$ is contained in $W$ and $A_1 = A_1'$.  Equivalently,
$F^{-1}(A_1) \cap F(C) \subset A_1$.

Continue inductively, choosing a
sequence $\{ A_n \}$ of compact subsets of $W$ such that for $n= 0,1,...$
\begin{equation}\label{eq1.42}
\begin{split}
A_n \quad \subset \subset \quad A_{n+1} \quad \subset   W, \\
F^{-1}(A_n) \cap F(C) \quad \subset \quad A_n. \hspace{1cm}
\end{split}
\end{equation}

By the remark following Lemma \ref{lem1.31}, $G = \bigcup_n A_n$ is an open set with $F^{-1}(G) \cap F(C) \subset G$.
The Lemma implies that $F^{-1}(G)$ is disjoint from $[[B]]$.  Since $G$ contains $A$, $F^{-1}(A)$ is
disjoint from $[[B]]$.

The existence of the Lyapunov function $L$ then follows from Corollary \ref{cor1.10} applied to the disjoint sets
$F^{-1}(A)$ and $B$.

(iii) If $L$ extends to the -unique- continuous function $\hat{L}$ on $\hat{X}$ then
$F \subset \leq_{\hat{L}}$ and so its closure $\hat{F}$ and transitive extension $\mathcal{G}\hat{F}$ are
contained in $\leq_{\hat{L}}$ as well. Thus, as usual, $\hat{L}$ is a Lyapunov function for $\mathcal{G}\hat{F}$.

For the final result we use the transfinite construction of (\ref{eq1.6}) for $\mathcal{G}\hat{F}$:
\begin{equation}\label{eq1.43}
\begin{split}
R_0 \quad ={def} \quad \hat{F} \quad = \quad \overline{F} \\
R_{\alpha + 1} \quad =_{def} \quad \overline{ \bigcup_{n=1,...} R_{\alpha}^n } \\
R_{\alpha} \quad =_{def} \quad \overline{\bigcup_{\beta < \alpha} R_{\beta}}
\end{split}
\end{equation}
with the third applying when $\alpha$ is a limit ordinal.

Recall that $C$ is a neighborhood of $F(A) \cap \overline{W}$ which contains $A$.
 Consider the compact relation
$\tilde{F} = F \cap \overline{W} \times \overline{W}$ which is clearly a reflexive and
transitive relation on $\overline{W}$.

With $C^{\circ}$ the interior of $C$, define
\begin{equation}\label{eq1.44}
\tilde{G} \quad = \quad \{ x \in  C^{\circ} : \tilde{F}(x) \subset C^{\circ} \}
\quad = \quad C^{\circ} \setminus \tilde{F}^{-1}(\overline{W} \setminus C^{\circ}).
\end{equation}
Because $\tilde{F}$ is a compact relation $\tilde{G}$ is an open subset of $X$ and transitivity implies that
$\tilde{G}$ is $\tilde{F}$ invariant.  Clearly, $A \subset \tilde{G} \subset C$.

It suffices to show that if $(x,z) \in R_{\alpha} \cap \tilde{G} \times (\hat{X} \setminus X)$ then $\hat{L}(z) = 1$.

Since the compactification is dynamic and $F$ is reflexive we have that $F = (X \times X) \cap \mathcal{G}\hat{F}$ and
hence $F = (X \times X) \cap R_{\alpha}$ for all $\alpha$.

{\bfseries Case} $ 0 $:   If $(x,z) \in R_0$ with $x \in  G$ and $z \in \hat{X} \setminus X$ then there is a net
$(x_i,y_i) \in F$ converging to $(x,z)$ we can assume that $x_i \in G$ and $y_i \in X \setminus W$ for
all $i$ (Recall that $W$ is bounded).  Hence, $y_i \in (X \setminus W) \cap F(C)$. Thus, $L(y_i) = 1$
for all $i$ and in the limit $L(z) = 1$.

{\bfseries Case} $\alpha + 1$: If $(x,z) \in R_{\alpha + 1}$ then there exists a positive integer $n$ and
a sequence $x = a_0, a_1,..., a_n = z$ in $\hat{X}$ such that $(a_i,a_{i+1}) \in R_{\alpha}$ for $0 \leq i < n$.
Let $j$ be the smallest index such that $a_j \in \hat{X} \setminus X$.  It suffices to show that $\hat{L}(a_j) = 1$
for then $(a_j,z) \in \mathcal{G}\hat{F}$ implies $L(z) = 1$. Now for all $i < j$
$a_{i} \in X$ and since the compactification
is dynamic  and $x \in C$, $a_{i} \in F(C)$. If for some $i < j$ $a_{i} \not\in W $ then
$a_{i} \in F(C) \setminus W  \subset B$
and so $\hat{L}(a_{i}) = L(a_{i}) = 1$.  So $\hat{L}(a_{j}) = 1$ because $\hat{L}$ is a Lyapunov function.
Assume now that $a_i \in W$ for $i = 0,..,j-1$.  We have assumed that $x \in \tilde{G}$.  By definition of
$\tilde{F}$ $a_i \in \tilde{F}(x)$ for $i = 0,...,j-1$ and it follows that $\tilde{F}(a_i) \subset \tilde{F}(x)$
by transitivity of $F$. Hence, $a_i \in \tilde{G}$ for $i = 0,...,j-1$. Thus,
$(a_{j-1},a_j) \in R_{\alpha} \cap \tilde{G} \times (\hat{X} \setminus X)$.
By induction hypothesis, $\hat{L}(a_j) = 1$
as required.

{\bfseries Case} limit $\alpha$: As in the first case, there is a net $(x_i,y_i) \in R_{\beta_i}$ converging
to $(x,z)$ with $\beta_i < \alpha$. We can assume that $x_i \in \tilde{G}$ and $y_i \not\in \overline{W}$ for all $i$.
If $y_i \in X$ then, again because the compactification is dynamic,
$(x_i,y_i) \in F$. Thus, $y_i \in F(C) \setminus W$ and so $\hat{L}(y_i) =
L(y_i) = 1$. On the other hand, if $y_i \not\in X$ then
$(x_i,y_i) \in R_{\beta_i} \cap \tilde{G} \times (\hat{X} \setminus X)$ and so
$\hat{L}(y_i) = 1$ by induction hypothesis.  We then obtain $\hat{L}(z) = 1$ by taking the limit.

$\Box$

The question arises  about the case when $B$ is empty, i.e. $F(C) \subset W$.  The above proof does not
need that $B$ is nonempty, but when $B$ is empty we can make a stronger statement.

If $F(C)$ is a subset of $W$ then
$A \subset C$  implies $F(A) \subset W$ and so
\begin{equation}\label{eq1.45}
F(A) \ = \ F(A) \cap \overline{W} \ \subset \ C.
\end{equation}
That is,  $F(A)$ is a compact $F$ invariant set which admits a
neighborhood $C$ such that  $F(C)$ is compact. Hence, by Theorem \ref{th1.12} there is a Lyapunov function function
$L$ with compact support which is $1$ on $A$. Consequently, for any proper compactification $\hat{X}$ of $X$,
$L$ extends continuously to $\hat{L}$ which is $0$ on $\hat{X} \setminus X$.

Now we apply these results.

\begin{theo}\label{th1.34} Assume that $f$ is a closed relation on $X$ such that whenever $E$ is a compact
$\G f \cap \G f^{-1}$ equivalence class in $|\G f|$ there exists a compact neighborhood $U$ of $E$
such that $\G f(U) \cap \G f^{-1}(U)$ is compact. There then exist Lyapunov compactifications $(\X, \f)$
of $(X,f)$ such that if $E$ is a compact $\G f \cap \G f^{-1}$ equivalence class in $|\G f|$
then $E = \G \f(E) \cap \G \f^{-1}(E)$, i.e. $E$ is a  $\G \f \cap \G \f^{-1}$ equivalence class in $\X$.
The maximal Lyapunov compactification is such compactification and if $X$ is metrizable $(\X, \f)$
can be chosen to be a Lyapunov compactification with $\X$ metrizable.
\end{theo}

{\bfseries Proof:}  Let $\L_0$ be any sufficient set of $f$ Lyapunov functions and so of $F = 1_X \cup \G f$
Lyapunov functions.  If $E$ is a compact equivalence class for $F \cap F^{-1}$ then either $F$ is
a $\G f \cap \G f^{-1}$ equivalence class in $|\G f|$ or $E = \{ x \}$ with $x \not\in |\G f|$.  In the
latter case, by Lemma \ref{lem1.14} we can choose $U$ a compact neighborhood of $x$ such that
$\G f(U) \cap \G f^{-1}(U) = \emptyset$ and in particular with $U$ disjoint from $|\G f|$.
Hence, $U = F(U) \cap F^{-1}(U)$.

Let $\E =  \{ E : E $ is a compact $F \cap F^{-1}$ equivalence class  $\}$. For each $ E \in \E$
we can choose a compact neighborhood $U_E$ of $E$ such that $ F(U_E) \cap F^{-1}(U_E)$ is a compact
unrevisited set.  As the compact neighborhoods $N$ of $E$ in $U_E$ decrease to $E$ the
intersections $F(N) \cap F^{-1}(N)$  decrease to $ F(E) \cap F^{-1}(E) = E$.
So we can choose an unrevisited compact set $A_E$ with $E  \subset \subset A_E \subset \subset U_E$.
Apply Lemma \ref{lem1.33} to obtain  a Lyapunov function $L_E$ as in the statement of the
lemma. Let $\E_0 $ be any subset of $ \E$ such that $\{ (A_E)^{\circ} : E \in \E_0 \}$ is an open
cover of $ \bigcup \ \E$.  Let $\L $ be any collection of Lyapunov functions
which contains $ \L_0 \cup \{ L_E : E \in \E_0 \}$.

Let $(\X, \f)$ be the $\L$ compactification of $(X,f)$. Since $\L$ contains $\L_0$ it is a sufficient
set of Lyapunov functions and so $(\X, \f)$ is a dynamic compactification by Theorem \ref{th1.18}. If $E$ is
any compact $\G f \cap \G f^{-1}$ equivalence class then $E \in \E$ and so there exists $E' \in \E_0$
such that $E$ meets the interior of $A_{E'}$.  Because $A_{E'}$ is unrevisited and $E$ is an equivalence
class, it follows that $E \subset A_{E'}$.  Hence, if $x \in E$ and $ z \in (\X \setminus X) \cap \G \f(x)$
Lemma 1.24 implies that $\hat{L}_{E'}(z) = 1$ and $\hat{L}_{E'}(x) = 0$.  Hence,   $z$ is not
$\G \f \cap \G \f^{-1}$ equivalent to $x$.

If we choose $\L$ to be the set of all Lyapunov functions then we obtain the maximal Lyapunov compactification.
On the other hand, if $X$ is metrizable then we can choose $\L_0$ countable and by using the Lindel\"{o}f
property of $ \bigcup \E$, we can choose $\E_0$ countable as well.  Hence, $ \L_0 \cup \{ L_E : E \in \E_0 \}$
is countable and the associated compactification is metrizable.

$\Box$ \vspace{.5cm}

Now we prove Theorem \ref{cor1.35} which says that for  a closed, transitive relation $F$
which is  reflexive and anti-symmetric,
the following conditions are equivalent \begin{itemize}
\item[(i)] For every $x \in X$ there exists a compact neighborhood $U$ of $x$ such that
$F(U)\cap F^{-1}(U)$ is compact.
\item[(ii)] Every $x \in X$ has a compact unrevisited neighborhood.
\item[(iii)] There exists a proper compactification $(\X, \hat{F})$ for $(X, F)$ such that
$\G \hat{F}$ is reflexive and anti-symmetric on $\X$.
\item[(iv)] The Lyapunov functions for $F$ determine the topology of $X$.
\end{itemize}

{\bfseries Proof of Theorem \ref{cor1.35}:}  (i) $\Leftrightarrow$ (ii):
Since $F$ is reflexive, $U \subset U' = F(U) \cap F^{-1}(U)$.
So if $U$ is a neighborhood of $x$ then $U'$ is an unrevisited neighborhood of $x$.  If $U$ is unrevisited
then $U = U'$.

(i) $\Rightarrow$ (iii):     Choose a Lyapunov compactification $(\X, \hat{F})$ of $(X, F)$ which satisfies the
conditions of Theorem \ref{th1.34}.  Since $1_X \subset F, \ 1_{\X} \subset \hat{F} \subset \G \hat{F}$ and so the
latter is reflexive. If $x \in X$ then $\{ x \}$ is the $\F \cap \F^{-1}$ equivalence class by Theorem \ref{th1.34}.
In particular,  if $z \in \X \setminus X$ then $z$ is not $F \cap F^{-1}$ equivalent to any point of $X$.
By Lemma \ref{lem1.17} it is not equivalent to any other point of $\X \setminus X$.
Hence, $\G \hat{F}$ is anti-symmetric
as well as reflexive.

(iii) $\Rightarrow$ (iv):  Let $(\X, \hat{F})$ be a proper compactification such that $\G \hat{F}$ is
reflexive and anti-symmetric. The set of Lyapunov functions for $\hat{F}$ distinguishes points of $\X$ and so
generates the topology of $\X$ by Proposition \ref{prop1.25new} (b). Since the
compactification is proper, the restrictions
of these to $X$ form a set of Lyapunov functions for $F$ which determines the topology of $X$ by Proposition
\ref{prop1.25new} (c).

(iv) $\Rightarrow$ (iii): Apply Theorem \ref{th1.26} (a) to the maximal Lyapunov compactification $(\beta_F X, \hat{F})$.

$\Box$ \vspace{.5cm}

\section{Compactifications of a Flow}\label{secflow}

 With ${\mathbb R}_+ = [0,\infty)$ and $\phi :
{\mathbb R}_+ \times X \to X$ a continuous map, we write
\begin{equation}\label{eq4.1}
 \phi(t,x) \quad = \quad \phi^t(x) \quad = \quad
\phi_x(t)
\end{equation}
\index{$\phi^t$}%
\index{$\phi_x$}
so that $\phi^t$ is a continuous map on $X$. For $K$ any compact
subset of ${\mathbb R}_+$ we define the relation $\phi^K$ on $X$
by
\begin{equation}\label{eq4.2}
\phi^K \quad = \quad \bigcup_{t \in K} \phi^t \quad = \quad \{
(x,y) : y = \phi(t,x) \quad \mbox{for some} \ t \in K \}.
\end{equation}
Regarding $\phi$ as a closed subset of $\R_+ \times X \times X$, we have that
$\phi^K = \pi_{23}(\phi \cap (K \times X \times X))$ \index{$\phi^K$}which is closed because the
restriction $\pi_{23}|(K \times X \times X)$ is proper. That is,
$\phi^K$ is a closed relation. Furthermore,  for $A$ a compact subset of $X$
\begin{equation}\label{eq4.3}
\phi^K(A) \quad = \quad \phi(K \times A)
\end{equation}
which is compact. Thus, the relation $\phi^K$ is + proper.
%
%We will call $\phi$ \emph{proper} when for every compact subset
%$K$ of ${\mathbb R}_+$ the restriction $\phi : K \times X \to X$
%is a proper map, i.e. the preimages of compacta are compact. This
%is equivalent to the condition that $\pi_1 \times \phi :{\mathbb
%R}_+ \times X \to {\mathbb R}_+ \times X$ is proper. Of course,
%$\phi$ proper implies that $\phi^t : X \to X$ is a proper map for
%every $t \in {\mathbb R}$. Clearly, if $X$ is compact then $\phi$
%is proper.

When $X$ is compact the composition of closed relations is closed.
Given merely local compactness this need not be true. However, we
do have:

\begin{lem}\label{lem4.1} Let $F$ be a closed relation on a $X$ and
let $\phi : {\mathbb R}_+ \times X \to X$ be a continuous map. If
$K$ is a compact subset of ${\mathbb R}_+$ then the relation $F
\circ \phi^K$ is closed. If  the restriction $\phi|K \times X $ is
a proper continuous map, then $\phi^K \circ F$ is closed.
\end{lem}

{\bfseries Proof:} The relation $\phi^K$ is $+$proper relation and it a
proper relation if $\phi|K \times X $ is
a proper continuous map.  So the results follow from Proposition \ref{prop1.2} (d).
%
%Assume that a net $\{ (x_k,y_k) \} $ converges
%to $(x,y)$ in $X \times X$. With $A$ a compact neighborhood of $y$
%eventually $y_k \in A$
%
%If  for all $ k \ (x_k,y_k) \in F \circ \phi^K$ then there exists
%$t_k \in K$ such that $(\phi(t_k,x_k),y_k) \in F$. By going to a
%subnet we can assume that $\{ t_k \}$ converges to $t \in K$ and
%so $\{ \phi(t_k,x_k) \}$ converges to $\phi(t,x)$.  Since $F$ is
%closed, $(\phi(t,x),y) \in F$. Hence, $(x,y) \in F \circ \phi^K$.
%
%If for all $k \ (x_k,y_k) \in \phi^K \circ F$ then there exist
%$t_k \in K$ and $ z_k \in X$ such that $(x_k,z_k) \in F$ and
%$\phi(t_k,z_k) = y_k$.   The net $\{ (t_k,z_k) \}$ is eventually
%in the set $(\phi|(K \times X))^{-1}(A)$ which is compact since
%$\phi|(K \times X)$ is assumed to be proper. Again, by going to a
%subnet, we can assume
% $\{ (t_k,z_k) \}$ converges to a point $(t,z) \in K \times
%X$. Then $(x,z) \in F$ and $\phi(t,z) = y$. Hence, $(x,y) \in
%\phi^K \circ F$.

$\Box$ \vspace{.5cm}

\begin{ex}\label{ex4.2}The composition of a closed relation with a continuous
map need not be closed.\end{ex}
For $x \in {\mathbb R}$ let $F(x) = 1/x$
for $x > 0$ and $F(x) = 0$ otherwise. Let $g(x) = arctan(x)$. The
map $F$ is a closed relation but the map $g \circ F$ is not
closed.

$\Box$ \vspace{.5cm}

The map $\phi : {\mathbb R}_+ \times X \to X$ is called a
\emph{semiflow} \index{semiflow} when it is an action of ${\mathbb R}_+$ on $X$.
That is,
\begin{equation}\label{eq4.4}
\phi^0 \quad = \quad 1_X \qquad \mbox{and} \qquad \phi^t \circ
\phi^s \quad = \quad \phi^{t+s}
\end{equation}
for all $t,s \in {\mathbb R}_+$. For a semiflow $\phi$  we will
call $f = \phi^1$ the \emph{time one} map of the semiflow.

If $f$ is a homeomorphism then so is $\phi^t$ for every $t \geq
0$. For such a \emph{reversible semiflow} \index{semiflow!reversible semiflow} we can extend $\phi$ to
an action of ${\mathbb R}$ on $X$, so that $\phi^{-t} =
(\phi^t)^{-1}$. The map from ${\mathbb R} \times X$ to $X$ is
continuous and the semigroup identity (\ref{eq4.4}) holds for all $t,s \in
{\mathbb R}$.  Such an ${\mathbb R}$  action is called a
\emph{flow}. For a reversible semiflow $\phi$ we define the
\emph{reverse semiflow} \index{semiflow!reverse semiflow} $\phi^{-1}: {\mathbb R}_+ \times X \to X$
by $(\phi^{-1})^t = (\phi^t)^{-1}$ for $t \in {\mathbb R}_+$.

For example, the \emph{constant flow on $X$} \index{semiflow!constant flow}is the projection
$\pi : {\mathbb R} \times X \to X$ with $\pi^t = 1_X$ for all $t$.
The \emph{translation flow on ${\mathbb R}\times X$} \index{semiflow!translation flow}
$ \ \tau:{\mathbb R} \times ({\mathbb R} \times X) \to {\mathbb R}
\times X $ is defined by $\tau^t(s,x) = (t+s,x)$. It restricts to
a semiflow on ${\mathbb R}_+ \times X$.

\begin{prop}\label{prop4.3} Let $\phi$ be a semiflow on $X$. The following
conditions are equivalent and when they hold we call $\phi$ a
\emph{proper semiflow} \index{semiflow!proper semiflow}
\begin{itemize}
\item[(i)] $\pi_1 \times \phi : {\mathbb R}_+ \times X  \to {\mathbb R}_+ \times
X$ is a proper continuous map, i.e. the preimage of every compact
set is compact.
\item[(ii)] For every compact $K \subset  {\mathbb R}_+ $ the
restriction $\phi|K \times X$ is a proper continuous map.
\item[(iii)] There exists $\epsilon > 0$ such that the restriction
$\phi|[0,\epsilon] \times X$ is a proper continuous map.
\end{itemize}
If $\phi$ is proper then $\phi^K$ is a proper relation on
$X$ for every compact $K \subset  {\mathbb R}_+ $ and, in particular, $\phi^t$ is a proper continuous map on
$X$ for every  $t \in  {\mathbb R}_+ $.

A reversible semiflow is proper.
\end{prop}

{\bfseries Proof:} (i) $\Rightarrow$ (ii): If
$A \subset X$ then
\begin{equation}\label{eq4.5x}
  (\phi|(K \times X))^{-1}(A) = (\pi_1 \times \phi)^{-1}(K \times A).
\end{equation}

(ii) $\Rightarrow$ (iii): Obvious.

(iii) $\Rightarrow$ (i): Clearly, (iii) implies that $\phi^s$ is
proper for $s \in [0,\epsilon]$. Now express any compact $A
\subset {\mathbb R}_+ \times X$ as a finite union of pieces $A_n
\subset [n\epsilon,(n+1)\epsilon]$ and observe that for $t \in
[n\epsilon,(n+1)\epsilon] \  \ \ (t,\phi(t,x)) =
(t,\phi(t-n\epsilon,\phi^{n\epsilon}(x)))$. This shows that the
restriction of $\pi \time \phi$ to each $[n\epsilon,(n+1)\epsilon]
\times X$ is proper and so the entire map is proper.

The set $(\phi^K)^{-1}(A)$ is the projection of $(\pi_1 \times \phi)^{-1}(K \times A)$ to the second coordinate.
Hence, $\phi^K$ is proper when $K$ is compact and $\phi$ is proper.

When $\phi$ is a reversible semiflow the map $\pi_1 \times \phi$
is a homeomorphism of ${\mathbb R}_+ \times X$. Hence, a
reversible semiflow is proper.

$\Box$ \vspace{.5cm}

In what follows we will let $I = [0,1]$ and $J = [1,2]$.
\index{$\phi^I$}%
\index{$\phi^J$} Observe
that any  real number $t \geq 0$ can be written $t = n + s$ with
$n = 0,1,2...$ and $s \in I$ and if $t \geq 1$ we can instead use
$s \in J$. Thus, any $\phi^t$ is a composition of a finite number
of functions $\phi^s$ with $s \in I$ and if $t \geq 1$ we can use
functions with $s \in J$. From this we obtain:
\begin{equation}\label{eq4.5}
\begin{split}
\mathcal{O} (\phi^I) \quad = \quad \phi^I \cup \mathcal{O}
(\phi^J) \hspace{1.5cm} \\ \mathcal{O}(\phi^J) \quad = \quad
\phi^I
\circ \mathcal{O}f \quad = \quad \mathcal{O}f \circ \phi^I \\
= \quad  \mathcal{O}(\phi^J) \circ \mathcal{O}(\phi^I) \quad =
\quad \mathcal{O}(\phi^I) \circ \mathcal{O}(\phi^J).
\end{split}
\end{equation}

For a semiflow $\phi : {\mathbb R}_+ \times X \to X$  we use the
closed relation $\phi^I$ to define
\begin{equation}\label{eq4.6}
\begin{split}
\mathcal{O} \phi  \quad =_{def}  \quad \mathcal{O} (\phi^I) \quad
= \quad
\bigcup \{ \phi^t : t \geq 0 \}.  \\
\mathcal{R}\phi   \quad =_{def}  \quad \mathcal{R} (\phi^I). \hspace{3cm} \\
\mathcal{N} \phi  \quad =_{def}  \quad \mathcal{N} (\phi^I) \quad
= \quad \overline{\mathcal{O} \phi} \quad = \quad \mathcal{O} \phi \ \cup \ \Omega \phi, \\
\mbox{where} \quad \Omega \phi \quad =_{def} \quad limsup_{t \to
\infty} \{ \phi^t \}.  \\
 \mathcal{G} \phi \quad =_{def}
\quad \mathcal{G} (\phi^I).\hspace{3cm} \\
 \mathcal{G} \phi^{-1} \quad =_{def}
\quad \mathcal{G} (\phi^I)^{-1}.\hspace{2.5cm} \\
 \mathcal{G} (\phi \cup \phi^{-1}) \quad =_{def}
\quad \mathcal{G} (\phi^I \cup (\phi^I)^{-1}).\hspace{1cm}
\end{split}
\end{equation}
Clearly, $\mathcal{G}\phi$ is the smallest closed transitive
relation on $X$ which contains  the maps $\phi^t$ for all $t \geq
0$ and $\mathcal{G} \phi^{-1}$ is its reverse relation.  Since
$\phi^0$ is the identity, $\mathcal{G} \phi$ is reflexive.
$\mathcal{G}(\phi \cup \phi^{-1})$ is the smallest closed
equivalence relation which contains the maps $\phi^t$.

$\mathcal{G}(\phi^J)$ is the smallest closed, transitive relation
on $X$ which contains the maps $\phi^t$ for all $t \geq 1$. We
call a point a \emph{generalized recurrent point}
\index{semiflow!generalized recurrent point}%
\index{recurrent point!generalized recurrent point for $\phi$}%
for $\phi$ when
it lies in $|\mathcal{G}( \phi^J)|$.

We will call $A \subset X \ \phi \ +$invariant (or $\phi$
invariant) for a semiflow $\phi$ on $X$ when it is $\phi^t \
+$invariant (resp. $\phi^t$ invariant ) for every $t \in {\mathbb
R}_+$. Clearly, $A$ is $\phi \ +$invariant iff it is $+$invariant
for the relation $\phi^I $.

%
%So we use $\phi^J$ to define
%%4.6
%\begin{equation}
%\Omega \mathcal{G} \phi \quad =_{def} \quad \Omega \mathcal{G}
%(\phi^J) \quad = \quad \bigcap_n (\mathcal{G} (\phi^J))^n.
%\end{equation}
We collect some useful identities for these relations.  Recall that for a closed relation
$f$ on $X$ and $A \subset X, \ [[A]]_f$ denotes the smallest closed $f \ +$invariant subset which contains
$A$.

\begin{prop}\label{prop4.4} Let $\phi$ be a semiflow on  $X$ with $f$  the time one map, i.e.
$f = \phi^1$. \begin{enumerate}
\item[(a)]For every $t \in
{\mathbb R}_+$ and for $K$ any compact subset of ${\mathbb R}_+$
\begin{equation}\label{eq4.7}
%\begin{split}
%\phi^t \circ \mathcal{G}f \quad \subset \quad \mathcal{G}f \circ
%\phi^t \\
%\phi^t \circ \Omega \mathcal{G}f \quad \subset \quad
%\Omega \mathcal{G}f \circ \phi^t \\
 \phi^t \circ \mathcal{G}(\phi^K) \quad \subset \quad
\mathcal{G}(\phi^K) \circ \phi^t
%\end{split}
\end{equation}
with equality when $\phi$ is reversible.

\item[(b)] For $K$ any compact subset of ${\mathbb R}_+$
\begin{equation}\label{eq4.8}
\phi^K \cup \phi^K \circ \mathcal{G}(\phi^K) \quad \subset \quad
\phi^K \cup \mathcal{G}(\phi^K) \circ \phi^K \quad = \quad
\mathcal{G}(\phi^K)
\end{equation}
with equality when  $\phi$ is proper.

\item[(c)]
\begin{equation}\label{eq4.9}
\mathcal{G} \phi \quad = \quad \phi^I \cup \mathcal{G}(\phi^J)
\hspace{2cm}
%\quad = \quad \mathcal{O}\phi \cup \Omega \mathcal{G}\phi.
\end{equation}
\vspace{.2cm}
\begin{equation}\label{eq4.10}
\begin{split}
\mathcal{G}(\phi^J) \quad = \quad  \mathcal{G}\phi \circ  \mathcal{G}(\phi^J)
\quad = \quad \mathcal{G}(\phi^J) \circ \mathcal{G}\phi.\\
\mathcal{G}(\phi^J) \quad  = \quad  \mathcal{G}\phi \circ  f
\quad \supset \quad f \circ \mathcal{G}\phi. \hspace{1cm}\\
\mathcal{G}(\phi^J) \quad = \quad \mathcal{G}f \circ \phi^I \quad
\supset
\quad \phi^I \circ \mathcal{G}f. \hspace{.5cm}% \\
%\Omega \mathcal{G}\phi \quad = \quad \Omega \mathcal{G}f \circ
%\phi^I \quad = \quad \phi^I \circ \Omega \mathcal{G}f.
\end{split}
\end{equation}
The first inclusion is an equality when $\phi$ is proper. Both
inclusions are equalities when $\phi$ is reversible.
\item[(d)]
If $\phi$ is reversible then
\begin{equation}\label{eq4.11}
\mathcal{O}(\phi^{-1}) \quad = \quad (\mathcal{O}\phi)^{-1} \qquad
\mbox{and} \qquad \mathcal{G}(\phi^{-1}) \quad = \quad
(\mathcal{G}\phi)^{-1}
\end{equation}
and so we can -without ambiguity - omit the parentheses in these
expressions.
\item[(e)] If $A$ is a closed subset of $X$,  then
\begin{equation}\label{eq4.12}
\begin{split}
[[A]]_{\G \phi} \quad \supset \quad \phi^I(A) \cup [[A]]_{\G (\phi^J)} \quad = \hspace{2cm}
\\ \phi^I(A) \cup [[\G (\phi^J)(A)]]_{\G (\phi^J)} \quad = \quad \phi^I(A) \cup [[\phi^J(A)]]_{\G (\phi^J)},
\end{split}
\end{equation}
with equality if $\phi$ is reversible.
\end{enumerate}
\end{prop}

{\bfseries Proof:} If $g: X_1 \to X_2 $ is a map and $A_1,A_2$ are
relations on $X_1$ and $X_2$ respectively then the following six
inclusions are all equivalent:
%4.12
\begin{center}\begin{equation}\label{eq4.13}
\begin{split}
 (g \times g)(A_1)  \subset  A_2, \qquad A_1
\subset (g \times g)^{-1}(A_2), \\
g \circ A_1 \circ g^{-1} \subset A_2, \qquad A_1 \subset g^{-1} \circ A_2 \circ g, \\
g \circ A_1 \subset A_2 \circ g, \qquad A_1 \circ g^{-1} \subset
g^{-1} \circ A_2,
\end{split}
\end{equation}\end{center}
because each says that $(x,y) \in A_1$ implies $(g(x),g(y))
\subset A_2$. When   $X_1 = X_2$ and $A_1 = A_2$, these say that
$A_1$ is a $g \ +$invariant relation.

It is clear that $\phi^K$ is $\phi^t \ +$invariant and hence
$(\phi^t \times \phi^t)^{-1}(\mathcal{G}(\phi^K))$ is a closed
transitive relation which contains $\phi^K$. Hence, it contains
$\mathcal{G}(\phi^K)$. So (\ref{eq4.7}) follows from (\ref{eq4.13}).  When $\phi$
is reversible we can apply the same argument to the entire
associated flow and so (\ref{eq4.7}) holds with $t$ replaced by $-t$. This
implies equality in (\ref{eq4.7}).

The three inclusions in (\ref{eq4.8}) and (\ref{eq4.10}) follow from (\ref{eq4.7}) as does
equality when $\phi$ is reversible.

In (\ref{eq4.8}) each of the three relations contains $\phi^K$ and is
contained in $\mathcal{G}(\phi^K)$. Since $\mathcal{G}(\phi^K)$ is
transitive and composition distributes over union, it is easy to
check that each relation is transitive. Since \\
$\mathcal{G}(\phi^K)\circ \phi^K \cup \phi^K$ is closed by Lemma 4.1, it equals
$\mathcal{G}(\phi^K)$. The same is true for the first relation
when it is closed, e.g. when $\phi$ is proper.

Next we will use (\ref{eq4.8}) to show that
\begin{equation}\label{eq4.14}
\mathcal{G}(\phi^J)\quad \subset \quad \phi^I \circ \mathcal{G}(\phi^J) \quad \subset \quad \mathcal{G}(\phi^J) \circ
\phi^I \quad \subset \quad \mathcal{G}(\phi^J), 
\end{equation}
and so they are all equal. 

Because $\phi^I$ is reflexive, $\phi^I \circ \mathcal{G}(\phi^J)$
contains $\mathcal{G}(\phi^J)$. The second inclusion follows from (\ref{eq4.7}).  By (\ref{eq4.8}) 
$\mathcal{G}(\phi^J) \circ \phi^I $ equals
$[\mathcal{G}(\phi^J) \circ \phi^J \cup \phi^J] \circ \phi^I$.
Since $\phi^J \circ \phi^I \subset \mathcal{O}(\phi^J) \subset
\mathcal{G}(\phi^J)$, it is contained in $\mathcal{G}(\phi^J)$ as
well.

The closed relation $\mathcal{G}(\phi^J) \cup \phi^I$ contains
$\phi^I$ and is contained in $\mathcal{G}(\phi^I)$.   From (\ref{eq4.14})
it follows that $\mathcal{G}(\phi^J) \cup \phi^I$ is transitive and so contains
$\mathcal{G}(\phi^I)$, proving (\ref{eq4.9}).

Each of the relations in the first line of (\ref{eq4.10}) contains
$\mathcal{G}(\phi^J)$ because $\mathcal{G}\phi$ is reflexive. For
the reverse inclusion substitute from (\ref{eq4.9}) and observe that
$\phi^I \circ \mathcal{G}(\phi^J) = \mathcal{G}(\phi^J)
\circ \phi^I = \mathcal{G}(\phi^J)$ by  (\ref{eq4.14}).

Next observe that $\mathcal{G}(\phi^I) \circ f $ and $\mathcal{G}f
\circ \phi^I$ contain $\phi^J = \phi^I \circ f = f \circ \phi^I$.
Both of these relations are closed by Lemma \ref{lem4.1}.

 From (\ref{eq4.9})
\begin{equation}\label{eq4.15}
\mathcal{G}(\phi^I) \circ f = [\mathcal{G}(\phi^I) \circ \phi^I
\cup \phi^I] \circ f = \mathcal{G}(\phi^I) \circ \phi^J \cup
\phi^J
\end{equation}
which is contained in $\mathcal{G}(\phi^J)$ by the already proved
first line of (\ref{eq4.10}). It is easy to see that this relation is
transitive and so it contains $\mathcal{G}(\phi^J)$.

Also, $f \circ \mathcal{G}(\phi^I)$ is a transitive relation which
contains $\phi^J$.  If $\phi$ is proper then the relation is
closed as well and so contains $\mathcal{G}(\phi^J)$.

Similarly, by (\ref{eq4.8}) with $K = \{ 1 \}$
\begin{equation}\label{eq4.16}
\mathcal{G}f \circ \phi^I = [\mathcal{G}f \circ f \cup f] \circ
\phi^I \subset \mathcal{G}(\phi^J).
\end{equation}
On the other hand, by (\ref{eq4.7}) and transitivity of $\mathcal{G}f$ the
composition of this relation with itself is contained in
\begin{equation}\label{eq4.17}
\mathcal{G}f \circ \phi^I \circ \phi^I \subset \mathcal{G}f \circ
[\mathcal{O}f \cup 1_X] \circ \phi^I \subset \mathcal{G}f \circ
\phi^I.
\end{equation}
That is, $\mathcal{G}f \circ \phi^I$ is transitive and so contains
$\mathcal{G}(\phi^J)$.  This completes the proof of (\ref{eq4.10}).

Finally, the reversibility results of part (d) follow because for
the reverse semiflow $\phi^{-1}$ the closed relation
$(\phi^{-1})^I$ is the reverse relation of $\phi^I$. That is, $y =
\phi^t(x)$ iff $x = \phi^{-t}(y)$.

(e)  Clearly,
\begin{equation}\label{eq4.18}
\begin{split}
[[A]]_{\G \phi} \ \supset \ \phi^I(A) \cup [[A]]_{\G (\phi^J)}  \ \supset \hspace{3cm}\\
\phi^I(A) \cup [[\G (\phi^J)(A)]]_{\G (\phi^J)}\ \supset \  \phi^I(A) \cup [[\phi^J(A)]]_{\G (\phi^J)}.\hspace{1cm}
\end{split}\end{equation}
Because $\phi^J$ is always $+$proper, (\ref{eq1.32}) implies that
\begin{equation}\label{eq4.19}
\begin{split}
\phi^I(A) \cup [[\phi^J(A)]]_{\G (\phi^J)} \quad = \hspace{4cm}\\
\phi^I(A) \cup (A \cup [[\phi^J(A)]]_{\G (\phi^J)}) \ = \ \phi^I(A) \cup [[A]]_{\G (\phi^J)}.\hspace{1cm}
\end{split}\end{equation}

Since $Q =_{def} \phi^I(A) \cup [[A]]_{\G \phi^J}$ contains $A$ it suffices to show that it is $+$invariant
with respect to $\G \phi = \phi^I \cup \G \phi^J$ when $\phi$ is reversible.

First we observe that (\ref{eq4.14}) implies
\begin{equation}\label{eq4.20}
\G (\phi^J)(\phi^I(A)) \ = \ \G (\phi^J)(A) \ \subset \ [[A]]_{\G (\phi^J)},
\end{equation}
and $\phi^I(\phi^I(A)) = \phi^I(A) \cup \phi^J(A)$. Thus, $\G \phi (\phi^I(A)) \subset Q$.

Obviously, $\G (\phi^J)([[A]]_{\G (\phi^J)})$ is contained in $Q$.  We are left with showing that
$\phi^t([[A]]_{\G (\phi^J)}$ is contained in $Q$ for all $t \in I$.

Since $\phi$ is reversible, each $\phi^t$ is a homeomorphism preserving the flow and so
$\phi^t([[A]]_{\G (\phi^J)}) = [[\phi^t(A)]]_{\G (\phi^J)}$ which equals
$\phi^t(A) \cup [[\phi^J(\phi^t(A)]]_{\G (\phi^J)}$ by (\ref{eq1.32}) again.   $\phi^t(A) \subset \phi^I(A) $ and
$\phi^J(\phi^t(A)) \subset \mathcal{O}  \phi^J(A)) $. Because
$\phi^J(A) \subset \mathcal{O} (\phi^J)(A) \subset \G (\phi^J)(A)$
it follows from (\ref{eq1.32}) that \\ $[[\phi^J(\phi^t(A))]]_{\G (\phi^J)} \subset Q$ as required.

$\Box$ \vspace{.5cm}

Recall that whether $\phi$ is reversible or not we  write
$\mathcal{G} \phi^{-1}$ for the reverse relation $(\mathcal{G}
\phi)^{-1} = \mathcal{G}(\phi^I)^{-1}$.

\begin{cor}\label{cor4.5} Let $\phi$ be a semiflow on $X$ with time one map $f$. Let $x \in X$.
\begin{itemize}
\item[(a)] The point $x$ is generalized recurrent,
i.e. $x \in |\mathcal{G}(\phi^J)|$, iff $(f(x),x) \in
\mathcal{G}\phi $.
\item[(b)]If $x \in |\mathcal{G}(\phi^J)|$ then
%4.17
\begin{equation}\label{eq4.21}
\begin{split}
\mathcal{G}(\phi^J)(x)  \quad = \quad \mathcal{G}\phi(x),\hspace{2.5cm} \\
\mathcal{G}(\phi^J)^{-1}(x)  \quad = \quad \mathcal{G}\phi^{-1}(x)
\hspace{2cm}
\end{split}
\end{equation}
and each of these sets is $\phi \ +$invariant.  If $\phi$ is
reversible  then each is $\phi$ invariant.
\item[(c)] If the equivalence class $\mathcal{G}\phi(x) \cap
\mathcal{G}\phi^{-1}(x)$ contains more than one point then $x$ is
generalized recurrent and $ \mathcal{G}\phi(x) \cap
\mathcal{G}\phi^{-1}(x) = \mathcal{G}(\phi^J)(x) \cap
\mathcal{G}(\phi^J)^{-1}(x)$.Thus,
\begin{equation}\label{eq4.22}
\mathcal{G}\phi \cap \mathcal{G}\phi^{-1} \quad = \quad 1_X \cup
[\mathcal{G}(\phi^J) \cap \mathcal{G}(\phi^J)^{-1}].
\end{equation}
\item[(d)] If $X$ is compact and $x \in |\mathcal{G}(\phi^J)|$
then $\mathcal{G}(\phi^J)(x)$ and $\mathcal{G}(\phi^J)(x) \cap
\mathcal{G}(\phi^J)^{-1}(x)$ are $\phi$ invariant.
\end{itemize}
\end{cor}

{\bfseries Proof:} (a): By (\ref{eq4.10}) $\mathcal{G}\phi^J =
\mathcal{G}\phi \circ f$ and so
\begin{equation}\label{eq4.23}
(y,x) \in \mathcal{G}(\phi^J) \qquad \Longleftrightarrow \qquad
(f(y),x) \in \mathcal{G}\phi.
\end{equation}
Apply this first with $y = x$ to prove (a). Next note that
$(f(x),x) \in \G \phi$ and so $(x,x) \in \G (\phi^J)$ imply $(f(x),x) \in
\G (\phi^J) \circ \G \phi = \G (\phi^J)$.

Next, observe that $(y,x) \in
\mathcal{G}(\phi^K)$ implies $(f(y),f(x)) \in \mathcal{G}(\phi^K)$
by (\ref{eq4.7}). By induction and (\ref{eq4.23}) with $y = f^n(x)$ for $n =
1,2,...$ we see that $x \in |\mathcal{G}(\phi^J)|$ iff $(f(x),x)
\in \mathcal{G}\phi$ iff $(f^n(x),x) \in \mathcal{G}(\phi^J)$ for
any positive integer $n$.

(b):  Now assume that $x$ is  generalized recurrent.

For any $t \in {\mathbb R}_+$ choose $n > t + 1$ a positive integer.
Because $n \geq t \geq 0$, we have $ (x,\phi^t(x)) \in
\mathcal{G}\phi \circ \G (\phi^J) = \G (\phi^J),$ and $ (\phi^t(x),f^n(x)) \in \mathcal{G}(\phi^J) $.
By the argument in (a) $ \ (f^n(x),x) \in
\mathcal{G}(\phi^J)$.  It follows that $x$ is $\mathcal{G}(\phi^J)\cap
\mathcal{G}(\phi^J)^{-1}$ equivalent to $\phi^t(x)$ for every $t
\in {\mathbb R}_+$. In particular, $\mathcal{G}(\phi^J)(x)$
contains $\phi^I(x)$ and so equals $\mathcal{G}\phi(x)$ by (\ref{eq4.9}).

By (\ref{eq4.7}) $y \in \mathcal{G}(\phi^J)(x)$ implies $\phi^t(y) \in
\mathcal{G}(\phi^J)(\phi^t(x)) = \mathcal{G}(\phi^J)(x)$ for all
positive $t$. Similarly,  $y \in \mathcal{G}(\phi^J)^{-1}(x)$
implies $(\phi^t(y),\phi^t(x)) \in \mathcal{G}(\phi^J)$ and so
$(\phi^t(y),x) \in \mathcal{G}(\phi^J)$ for all positive $t$.
Hence, $\mathcal{G}(\phi^J)(x)$ and $\mathcal{G}(\phi^J)^{-1}(x)$
are $\phi \ +$invariant.

If $\phi$ is reversible then $(\phi^{-t}(x),\phi^{-t}(x)) \in
\mathcal{G}(\phi^J)$ for every positive $t$ by (\ref{eq4.7}). Hence,
$\phi^{-t}(x)$ is $\mathcal{G}(\phi^J)\cap
\mathcal{G}(\phi^J)^{-1}$ equivalent to every point in its forward
orbit which includes $x$. Thus, $x$ is $\mathcal{G}(\phi^J)\cap
\mathcal{G}(\phi^J)^{-1}$ equivalent to every point in its
backward orbit as well as its forward orbit. Furthermore, $y \in
\mathcal{G}(\phi^J)(x)$ implies $\phi^{-t}(y) \in
\mathcal{G}(\phi^J)(\phi^{-t}(x)) = \mathcal{G}(\phi^J)(x)$ for
all positive $t$. Hence, $\mathcal{G}(\phi^J)(x)$ is $\phi$
invariant. Similarly, $(y,x) \in \mathcal{G}(\phi^J)$ implies
$(\phi^t(y),\phi^t(x)) \in \mathcal{G}(\phi^J)$ and so
$\mathcal{G}(\phi^J)^{-1}(x)$ is $\phi$ invariant.

(c): Suppose that $x$ and $y$ are two distinct points with $(x,y),
(y,x) \in \mathcal{G}\phi$. If neither $x \in \phi^I(y)$ nor $y
\in \phi^I(x)$ then by (\ref{eq4.9}) $(x,y), (y,x) \in
\mathcal{G}(\phi^J)$ and by transitivity $(x,x), (y,y) \in
\mathcal{G}(\phi^J)$. Thus, $x$ and $y$ are generalized recurrent
and in the same $\mathcal{G}(\phi^J)\cap \mathcal{G}(\phi^J)^{-1}$
equivalence class.

Suppose instead that $y = \phi^t(x)$ for some $t \in I$. Since $y
\not= x, t > 0$. By (\ref{eq4.7}) $(\phi^{nt}(y),\phi^{nt}(x)) \in
\mathcal{G}\phi \cap \mathcal{G}\phi^{-1}$ for every positive
integer $n$ and so by induction $x \in
\mathcal{G}\phi(\phi^{nt}(x))$ for all such $n$. For $n$ large
enough that $nt > 1 \ (f(x), \phi^{nt}(x)) \in \mathcal{G}\phi$
and so by (a), $x$ is generalized recurrent. Equality of the
equivalence classes then follows from (b).

(d): When $X$ is compact we can  apply (\ref{eq2.8}) and (\ref{eq2.9})  to the
closed relation $\phi^J$ to get
\begin{equation}\label{eq4.24}
\begin{split}
\mathcal{G}(\phi^J) \quad = \quad \mathcal{O}(\phi^J) \cup \Omega
\mathcal{G}(\phi^J), \qquad \mbox{and} \hspace{1cm}\\
\Omega \mathcal{G}(\phi^J) \circ \phi^J \quad = \quad \Omega
\mathcal{G}(\phi^J) \circ \mathcal{G}(\phi^J) \quad = \quad \Omega
\mathcal{G}(\phi^J) \\ = \quad \mathcal{G}(\phi^J) \circ \Omega
\mathcal{G}(\phi^J)\quad = \quad \phi^J \circ \Omega
\mathcal{G}(\phi^J).
\end{split}
\end{equation}

From (\ref{eq2.10}) we have
\begin{equation}\label{eq4.25}
\begin{split}
|\mathcal{G}(\phi^J)| \quad = \quad |\Omega \mathcal{G}(\phi^J)| \hspace{4cm}\\
x \in |\mathcal{G}(\phi^J)| \qquad \Longrightarrow \hspace{4cm} \\
\mathcal{G}(\phi^J)(x) \  = \ \Omega \mathcal{G}(\phi^J) (x)\qquad
\mbox{and} \qquad  \mathcal{G}(\phi^J)^{-1}(x) \ = \ (\Omega
\mathcal{G}(\phi^J))^{-1} (x).
\end{split}
\end{equation}

Now fix  $ t \in (0,n+1]$. Clearly,
\begin{equation}\label{eq4.26}
\begin{split}
\mathcal{O}(\phi^J)\circ (\phi^J)^n \quad = \quad (\phi^J)^n \circ
\mathcal{O}(\phi^J) \quad \subset \\  \phi^t \circ
\mathcal{O}(\phi^J) \quad = \quad \mathcal{O}(\phi^J) \circ \phi^t
\quad \subset \quad \mathcal{O}(\phi^J)
\end{split}
\end{equation}

Composing these on the left and right as necessary with $\Omega
\mathcal{G}(\phi^J)$ we obtain:
\begin{equation}\label{eq4.27}
\Omega \mathcal{G}(\phi^J) \circ \phi^t \quad = \quad \Omega
\mathcal{G}(\phi^J) \quad = \quad \phi^t \circ \Omega
\mathcal{G}(\phi^J).
\end{equation}

From (\ref{eq4.27}) and (\ref{eq4.25}) it follows that when $x \in
|\mathcal{G}(\phi^J)|, \mathcal{G}(\phi^J)(x) = \Omega
\mathcal{G}(\phi^J)(x) $ is $\phi$ invariant. The proof that the
$\mathcal{G}(\phi^J) \cap \mathcal{G}(\phi^J)^{-1}$ equivalence
class of $x$ is $\phi$ invariant follows by the same argument as
was used in Lemma \ref{lem2.2}.

 $\Box$ \vspace{.5cm}

At least in the reversible case, the properties of generalized
recurrence for the semiflow and for the time one map agree. This
requires the following:

\begin{lem}\label{lem4.6} Let $T$ be a closed subset of ${\mathbb R}$ such that
\begin{itemize}
\item $t,s \in T$ implies $t + s \in T$.
\item $t \in T$ implies $t - 1 \in T$.
\item There exists $t \in T$ with $t > 0$.
\end{itemize}
Either $T = {\mathbb R}$ or $ T = \frac{1}{N} {\mathbb Z}$ for
some positive integer $N$. In particular, $T$ is an additive
subgroup of ${\mathbb R}$.
\end{lem}

{\bfseries Proof:}  Let $\xi = inf \{ t > 0 : t \in T \}.$

Case (1) $\xi = 0$: For every $\epsilon > 0$ there exists $t \in
T$ with $0 < t < \epsilon$.  Then $t{\mathbb Z}_+ \subset T$ and
is $\epsilon$ dense in ${\mathbb R}_+$.  Hence, ${\mathbb R}_+
\subset T$ because $T$ is closed. Repeated translation by $-1$
shows that ${\mathbb R} \subset T$.

Case (2) $\xi > 0$:  Since $T$ is closed $\xi \in T$. Let $N$ be
the smallest positive integer so that $N \cdot \xi \geq 1$. Hence,
$r = N \cdot \xi - 1 \in T$ and $ \xi > r \geq 0$. Minimality of
$\xi$ implies that $r = 0$ and so $\xi = \frac{1}{N}$. Hence,
$\frac{1}{N} {\mathbb Z}_+ \subset T$.  Translating repeatedly by
$-1$ shows that $\xi {\mathbb Z} = \frac{1}{N} {\mathbb Z} \subset
T$.  Finally, if $t \in T$ then let $m$ be the largest integer
such that $m \xi \leq t$ and so with $s = t - m \xi, \ 0 \leq s <
\xi$. Choose $k$ a positive integer large enough that $m + k$ is a
positive integer divisible by $N$ and so that $(m + k)\xi$ is a
positive integer. Then $s = (t + k\xi) -((m + k)\xi) \in T$. Since
$0 \leq s < \xi$, minimality again implies $s = 0$ and so $t = m
\xi$.  Thus, $\xi {\mathbb Z} = T$.

$\Box$ \vspace{.5cm}

\begin{theo}\label{th4.7} Let $\phi$ be a reversible semiflow on  $X$ with time one map $f$.
\begin{equation}\label{eq4.28}
|\mathcal{G}(\phi^J)| \quad = \quad |\mathcal{G}f|. \hspace{2cm}
\end{equation}
For $x \in |\mathcal{G}(\phi^J)| \ y$ is $ \mathcal{G}(\phi^J)
\cap \mathcal{G}(\phi^J)^{-1}$ equivalent to $x$ iff there exists
$t \in I$ such that $\phi^t(y)$ is $\mathcal{G}f \cap
\mathcal{G}f^{-1}$ equivalent to $x$. Furthermore,
\begin{equation}\label{eq4.29}
\mathcal{G}(\phi^J)(x) \cap \mathcal{G}(\phi^J)^{-1}(x) \quad =
\quad \phi^I(\mathcal{G}f(x) \cap \mathcal{G}f^{-1}(x)).
\end{equation}
\end{theo}

{\bfseries Proof:}  Write $\phi $ for the associated flow as well
as for the original semiflow. Let $x \in |\mathcal{G}(\phi^J)|$.
Define
\begin{equation}\label{eq4.30}
T \quad = \quad \{ t \in {\mathbb R} : (\phi^t(x),x) \in
\mathcal{G}f \}.
\end{equation}
If $t \in T$ then by (\ref{eq4.7}) $(\phi^{t + s}(x),\phi^{s}(x)) \in
\mathcal{G}f$ and so if $s \in T$ as well then transitivity
implies that $t + s \in T$. Similarly, $(\phi^{t -1}(x),f^{-1}(x))
\in \mathcal{G}f$ and since $(f^{-1}(x),x)$ is always in
$\mathcal{G}f$ we have that $t - 1 \in T$.

By (\ref{eq4.10}) $\mathcal{G}(\phi^J) = \mathcal{G}f \circ \phi^I $ and
so $x \in |\mathcal{G}(\phi^J)|$ implies there exists $s \in I$
such that $(\phi^s(x),x) \in \mathcal{G}f$.  If $s = 0$ then $x
\in |\mathcal{G}f|$ and so $(f^n(x),x) \in \mathcal{G}f$ for every
integer $n$. Otherwise, $s > 0$.  Either case shows there exist
positive $t$ in $T$.

Finally, $T$ is a closed set because $\mathcal{G}f$ is a closed
relation and $\phi$ is continuous.

By Lemma \ref{lem4.6}, $T$ is an additive subgroup of ${\mathbb R}$. Since
$0 \in T, \ x \in |\mathcal{G}f|$ proving (\ref{eq4.28}).

Since $t \in T$ iff $-t \in
T,$
\begin{equation}\label{eq4.31}
T \quad = \quad \{ t \in {\mathbb R} : \phi^t(x) \in
\mathcal{G}f(x) \cap \mathcal{G}f^{-1}(x) \}.
\end{equation}

Now if $y \in \mathcal{G}(\phi^J)(x) \cap
\mathcal{G}(\phi^J)^{-1}(x)$ then by (\ref{eq4.10}) again there exist $s,t
\in T$ such that $(\phi^s(x),y),(\phi^t(y),x) \in \mathcal{G}f$.
By (\ref{eq4.7}) $(\phi^{s+t}(x),\phi^{t}(y)) \in \G f$ and so by transitivity
$(\phi^{s+t}(x),x) \in \mathcal{G}f$.
In particular, $s + t \in T$. By (\ref{eq4.31}) $(x,\phi^{s+t}(x)) \in
\mathcal{G}f$ and so by transitivity $(x,\phi^t(y)) \in
\mathcal{G}f$. Hence, $\phi^t(y) \in \mathcal{G}f(x) \cap
\mathcal{G}f^{-1}(x)$.

Furthermore, $\phi^{-s}(y) \in \mathcal{G}f(\phi^{-s-t}(x)) \cap
\mathcal{G}f^{-1}(\phi^{-s-t}(x)) = \mathcal{G}f(x) \cap
\mathcal{G}f^{-1}(x)$ and $y = \phi^s(\phi^{-s}(y))$. From this
(\ref{eq4.29}) follows.

$\Box$ \vspace{.5cm}

We will need some results about the chain relation for a semiflow
$\phi$  on a compact space $X$. In Akin (1993) Chapter 6 the
notation $\mathcal{C}\phi$ was used for $\mathcal{C}(\phi^J)$ and
it follows from Proposition 2.4 (c) of Akin (1993) that $\mathcal{C}(\phi^J)
= \mathcal{O}(\phi^J) \cup \Omega \mathcal{C}(\phi^J)$.  We define
\index{$\tilde{\mathcal{C}}\phi$}
\begin{equation}\label{eq4.32}
\tilde{\mathcal{C}}\phi \quad =_{def} \quad \phi^I \cup
\mathcal{C}(\phi^J) \quad = \quad \mathcal{O}\phi \cup \Omega
\mathcal{C}(\phi^J).
\end{equation}

Usually the inclusion $\tilde{\mathcal{C}}\phi \subset
\mathcal{C}(\phi^I)$ is strict (contrast this with equation
(\ref{eq4.9}) \ ). For example, whenever $X$ is connected $X \times X =
\mathcal{C}1_X $ and so $ X \times X = \mathcal{C}(\phi^I)$. On
the other hand, if $f$ is surjective then the condition $X \times
X = \tilde{\mathcal{C}}\phi$ is equivalent to chain transitivity
of the semiflow which says there are no proper attractors.

\begin{prop}\label{prop4.8}  If $\phi$ is a semiflow on a compact  space $X$
then
\begin{equation}\label{eq4.33}
\Omega \mathcal{C}(\phi^J) \quad = \quad \bigcap_{n=1}^{\infty} \
\mathcal{C}((\phi^J)^n).
\end{equation}
\end{prop}

{\bfseries Proof:}  For a positive integer $n$ and a relation $g$
on $X$ we define the transitive relation $\mathcal{O}_n g =
\bigcup_{i \geq n} g^i$. Adjusting the result to eliminate the
metrizability hypothesis, Proposition 2.15 of Akin (1993) says
that $\Omega \mathcal{C}(\phi^J) = \bigcap_{n,V} \mathcal{O}_n (V
\circ \phi^J)$ intersecting over positive integers $n$ and
neighborhoods of the diagonal $V \in \mathcal{U}_X$. Clearly, $V
\circ (\phi^J)^n \subset \mathcal{O}_n ( V \circ \phi^J)$. Since
the latter is transitive: $\mathcal{C}((\phi^J)^n) \subset
\mathcal{O}(V \circ (\phi^J)^n) \subset \mathcal{O}_n( V \circ
\phi^J)$. Intersecting over $n$ and $V$ we get $\bigcap_n
\mathcal{C}((\phi^J)^n) \subset \Omega \mathcal{C}(\phi^J)$.

For the other direction we need uniform continuity of the semiflow
$\phi$ on $J \times X$. It then follows that for every $V \in
\mathcal{U}_X$ and positive integer $n$ there exists $W \in
\mathcal{U}_X$ such that for all integers $k$ with  $n \leq k \leq
2n \ \ \ (W \circ \phi^J)^k \subset V \circ (\phi^J)^k \subset
\mathcal{O}(V \circ (\phi^J)^n)$. Writing any $p \geq n$ as a sum of $k$'s between
$n$ and $2n$ we see that for any  $p \geq n,  (W \circ \phi^J)^p \subset \mathcal{O}(V \circ (\phi^J)^n)$.  
Now take the union over the
$p$'s, and use Akin (1993) Proposition 2.15
again to get $\Omega \mathcal{C}(\phi^J) \subset \mathcal{O}_n(W
\circ \phi^J) \subset \mathcal{O}(V \circ (\phi^J)^n)$. Compare
first and third of these and intersect first over $V$ and then
over $n$ to get the reverse inclusion in (\ref{eq4.33}).

$\Box$ \vspace{.5cm}

This result says that $(x,y) \in \Omega \mathcal{C}(\phi^J)$ iff
for every $V \in \mathcal{U}_X$ and every positive integer $n$,
there exists a chain $(x_0,t_0,x_1,....,t_k,x_k)$ with $x_0 = x,
x_k = y, t_i \geq n$ for $i = 0,...,k$ and
$(\phi(t_i,x_i),x_{i+1}) \in V$ for $i = 0,...,k-1$. This is the
relation $P_f$ introduced by Conley, e.g.  Conley (1988).

From (\ref{eq4.32}) it follows as with (\ref{eq4.22}) that
\begin{equation}\label{eq4.34}
\tilde{\mathcal{C}}\phi \cap \tilde{\mathcal{C}}\phi^{-1} \quad =
\quad 1_X \cup [\mathcal{C}(\phi^J) \cap
\mathcal{C}(\phi^J)^{-1}].
\end{equation}

Using Proposition 2.4 of Akin(1993) it is easy to check  the
following analogue of (\ref{eq4.25})
\begin{equation}\label{eq4.35}
\begin{split}
|\mathcal{C}(\phi^J)| \quad = \quad |\Omega \mathcal{C}(\phi^J)|.
\hspace{4cm}\\  \\
x \in |\mathcal{C}(\phi^J)| \qquad \Longrightarrow \hspace{4cm} \\
\tilde{\mathcal{C}}\phi(x)\quad = \quad \mathcal{C}(\phi^J)(x) \quad
= \quad \Omega \mathcal{C}(\phi^J) (x)\hspace{1cm}\\
\mbox{and}\hspace{4cm}
\\
(\tilde{\mathcal{C}}\phi)^{-1} (x) \quad = \quad
\mathcal{C}(\phi^J)^{-1}(x) \quad = \quad (\Omega
\mathcal{C}(\phi^J))^{-1} (x).
\end{split}
\end{equation}

In order to compactify a  semiflow $\phi$ on $X$ we begin by
compactifying the closed relation $\phi^I$, but a bit more is
needed to obtain a semiflow on the compactification.

For a bounded, continuous real valued function $u : X \to {\mathbb
R}$ and $t \in [0, \infty)$ define $\Delta_t u : X \to {\mathbb
R}$ by
\begin{equation}\label{eq4.36}
\Delta_t u \quad =_{def} \quad u \circ\phi^t - u.
\end{equation}
Since $\phi$ is a semiflow, we have for all $s,t \in {\mathbb
R}_+$
\begin{equation}\label{eq4.37}
\Delta_t u(\phi(s,x)) \quad = \quad u(\phi(s+t,x)) - u(\phi(s,x)).
%\quad = \quad \Delta_t (u \circ \phi^s)(x).
\end{equation}

We will call $u$ \emph{$\phi$ uniform} \index{$\phi$ uniform function}when for every $\epsilon >
0$ there exists a $\delta > 0$ such that  $t \leq \delta$ implies
$|\Delta_t u| \leq \epsilon $, i.e. $ |u(\phi(t,x)) - u(x)| \leq
\epsilon$ for all $(t,x) \in  [0,\delta] \times X$. We call $u$
\emph{$\phi$ Lipschitz} \index{$\phi$ Lipschitz function} with constant $M \in {\mathbb R}_+$ if for
all $t \in {\mathbb R}_+ \ |\Delta_t u| \leq M t$. Clearly, $\phi$
Lipschitz implies $\phi$ uniformity.

Clearly, $\phi$ uniformity says exactly that the function from
${\mathbb R}_+$ to $\mathcal{B}(X)$ given by $t \mapsto u \circ
\phi^t$ is continuous at 0 and from (\ref{eq4.37}) we see that this
implies uniform continuity on ${\mathbb R}_+$.

If $X$ is compact then every continuous real valued function $u$
is $\phi$ uniform. In the general, locally compact case, the
$\phi$ uniform functions comprise a closed subalgebra of
$\mathcal{B}(X)$ which we will denote  by $\mathcal{B}_{\phi}(X)$.\index{$\mathcal{B}_{\phi}(X)$}
It is easy to check that if $\phi$ is proper then
$\mathcal{B}_{\phi}(X)$ contains  $\mathcal{A}_0$ the closed
subalgebra generated by the  functions of compact support.

We say that a subset $\mathcal{K}$ of $\mathcal{B}(X)$ $\phi^* \
+$invariant (or $\phi^*$ invariant) when it is $(\phi^t)^* \
+$invariant (resp. when it is $(\phi^t)^*$ invariant) for all $ t
\in {\mathbb R}_+$. If $j: X \to \hat{X}$ is the compactification
associated with a $\phi^* \ +$ invariant closed subalgebra
$\mathcal{A}$ then the continuous maps $\phi^t$ on $X$ all extend
to maps $\hat{\phi}^t$ of $\hat{X}$. We obtain the semiflow
composition rule $\hat{\phi}^t \circ \hat{\phi}^s =
\hat{\phi}^{t+s}$ from the corresponding equation (\ref{eq4.4}) on $X$.
However, continuity of the semiflow $\hat{\phi}$, that is, joint
continuity of the map $\hat{\phi} : {\mathbb R}_+ \times \hat{X}
\to \hat{X}$ requires exactly that $\mathcal{A} \subset
\mathcal{B}_{\phi}(X)$. Necessity easily follows from the remarks
above concerning the compact case. For sufficiency:

\begin{lem}\label{lem4.10} Let $\phi : {\mathbb R}_+ \times X \to X$ be a
semiflow and $D$ be a dense subset of ${\mathbb R}_+$. Let
$\mathcal{A}$ be a closed subalgebra of $\mathcal{B}_{\phi}(X)$
with  $j : X \to \hat{X}$ be the associated compactification.

If for every $t \in D$
\begin{equation}\label{eq4.38}
(\phi^t)^*(\mathcal{A}) \quad \subset \quad \mathcal{A}.
\end{equation}
then $\mathcal{A}$ is $\phi^* \ +$ invariant and there is a unique
semiflow $\hat{\phi} : {\mathbb R}_+ \times \hat{X} \to \hat{X}$
such that $\hat{\phi}^t \circ j = j \circ \phi^t$ for all $t \in
{\mathbb R}_+$.  If $\phi$ is reversible and for all $t \in D$
\begin{equation}\label{eq4.39}
(\phi^t)^*(\mathcal{A}) \quad = \quad \mathcal{A}
\end{equation}
 then $\hat{\phi}$ is reversible.
\end{lem}

{\bfseries Proof:} Recall that $j^* : \mathcal{B}(\hat{X}) \to
\mathcal{A}$ is an isometric algebra isomorphism and for $u \in
\mathcal{A}$ we let $\hat{u}$ denote the unique function on
$\hat{X}$ such that $u = \hat{u} \circ j$.

By definition $t \to (\phi^t)^*(u)$ is continuous for $u \in
\mathcal{B}_{\phi}$ and so if $u \in \mathcal{A}$. If
$(\phi^t)^*(u) \in \mathcal{A}$ is true for $t \in D$ then it is
true for all $t \in {\mathbb R}_+$ because $\mathcal{A}$ is closed
and $D$ is dense. Thus, $\mathcal{A}$ is $\phi \ +$invariant. It
follows that for each $t$ the continuous map $\hat{\phi}^t$ is
uniquely defined on $\hat{X}$ extending $\phi^t$.

To prove joint continuity, let $t \in {\mathbb R}_+$ and $a,b \in
\hat{X}$ with $\hat{\phi}(t,a) = b$. Let $\hat{u} \in
\mathcal{B}(\hat{X})$ with $\hat{u}(b) = 1$ and with support an
arbitrarily small neighborhood of $b$. Let $u = \hat{u} \circ j$.
Since $|\hat{u} \circ \hat{\phi}^t - \hat{u}| = |u \circ \phi^t -
u|$ it follows that $t \mapsto \hat{u} \circ \hat{\phi}^t $ is a
continuous map from ${\mathbb R}_+$ to $\mathcal{B}(\hat{X})$.
From this it follows that $\hat{u} \circ \hat{\phi} : {\mathbb
R}_+ \times \hat{X} \to {\mathbb R}$ is continuous. Hence, if
$(s,c)$ is close enough to $(t,a)$ then $\phi(s,c)$ is in the
support of $\hat{u}$.

If $\phi$ is reversible then (\ref{eq4.39}) for $t$  implies (\ref{eq4.38}) for $-t$
and so the reverse semiflow extends to $\hat{X}$. Clearly, the
extension of the reverse flow is the reverse of the extension.

$\Box$ \vspace{.5cm}

\begin{df}\label{df4.11} Let $\phi$ be a semiflow on $X$ and let $\hat{X}$
 be a proper compactification with inclusion map $j : X \to
\hat{X}$ and associated algebra $\mathcal{A} =
j^*(\mathcal{B}(\hat{X})) \subset \mathcal{B}(X)$. We call
$\hat{X}$ a \emph{dynamic compactification for $\phi$}
\index{dynamic compactification for a semiflow}%
\index{compactification!dynamic compactification for a semiflow}%
\index{semiflow!dynamic compactification} when
\begin{itemize}
\item $\mathcal{A} \subset \mathcal{B}_{\phi}(X)$.
\item $\mathcal{A}$ is $\phi^* \ +$invariant.
\item $\hat{X}$ is a dynamic compactification for the closed
relation $\phi^I$.
\end{itemize}
\end{df} \vspace{.5cm}

\begin{theo}\label{th4.12} Let $\hat{X} \supset X$ be a dynamic compactification
for a semiflow $\phi$ on $X$ with inclusion map $j : X \to \X$. There is a unique semiflow
$\hat{\phi} : {\mathbb R}_+ \times \hat{X} \to \hat{X}$ such that
$\hat{\phi}^t \circ j = j \circ \phi^t$ for all $t \in {\mathbb
R}_+$. Furthermore,
\begin{equation}\label{eq4.40}
\begin{split}
(X \times X) \cap \mathcal{G} \hat{\phi} \quad = \quad
\mathcal{G}\phi \hspace{2cm}\\
(X \times X) \cap \mathcal{G} (\hat{\phi}^J) \quad = \quad
\mathcal{G}(\phi^J).\hspace{1cm}
\end{split}
\end{equation}

 If $\hat{E} \subset |\mathcal{G} (\hat{\phi}^J)| $
  is an $\mathcal{G} (\hat{\phi}^J) \cap
\mathcal{G} (\hat{\phi}^J)^{-1}$  equivalence class with $E =
\hat{E} \cap X$ then exactly one of the following three
possibilities holds:
\begin{enumerate}
\item[(i)]  $\hat{E} \subset \hat{X} \setminus
X$  and $E = \emptyset$.
\item[(ii)] $E$ is a noncompact $\mathcal{G}(\phi^J) \cap
\mathcal{G}(\phi^J)^{-1}$ equivalence class $E \subset X$ whose
closure meets $\hat{X} \setminus X$ and is contained in $\hat{E}$.
\item[(iii)] $\hat{E}$ is contained in $ X$, i.e. $\hat{E} = E$, and it is a compact
$\mathcal{G}(\phi^J) \cap \mathcal{G}(\phi^J)^{-1}$ equivalence
class.
\end{enumerate}
\end{theo}

{\bfseries Proof:} The semiflow $\hat{\phi}$ is defined by Lemma
\ref{lem4.10}.

The functions $\hat{\phi}^t$ extend $\phi^t$  and for any compact
subset $K$ of ${\mathbb R}_+ \ \hat{\phi}^K$ is the closure in
$\hat{X} \times \hat{X}$ of $ \phi^K$. So (\ref{eq4.40}) for the reflexive
relation $\mathcal{G}\phi = \mathcal{G}(\phi^I)$ follows because we have assumed the
compactification is dynamic for $\phi^I$. The result for $\mathcal{G}(\phi^J)$ follows from (\ref{eq4.10}).
That is, for $x,y \in X , \ (x,y) \in \mathcal{G}(\hat{\phi}^J)$
iff $(f(x),y) \in \mathcal{G}\hat{\phi}$ and so iff $(f(x),y) \in
\mathcal{G}\phi$ iff $(x,y) \in \mathcal{G}(\phi^J)$.

By (\ref{eq4.40}) the compactification is dynamic for $\phi^J$. Hence, the
  remaining results  follow from Corollary  \ref{cor1.10new} applied to the +proper relation $\phi^J$.

$\Box$ \vspace{.5cm}

We will call $L : X \to {\mathbb R}$ a Lyapunov function  \index{semiflow!Lyapunov function for a semiflow}%
\index{Lyapunov function!for a semiflow} for a
semiflow  $\phi$ when it is a Lyapunov function for the closed
relation $\phi^I$ associated with  $\phi$ and so it is a
$\mathcal{G}(\phi)$ Lyapunov function as well. Equivalently, $L$
is a Lyapunov function for every map $\phi^t$ with $t \in {\mathbb
R}_+$ and this condition holds iff $\Delta_t L \ \geq \ 0 $ for every
$t \in {\mathbb R}_+$. It then follows that for all $t \in
{\mathbb R}_+ \ L \circ \phi^t$ is a Lyapunov function for $\phi$.
If $\phi$ is reversible this is true for all $t \in {\mathbb R}$.

To construct a Lyapunov compactification for the semiflow we need
a uniform version of Corollary \ref{cor1.11}.

\begin{theo}\label{th4.13} Let $\phi$ be a semiflow on a space $X$, $A$ be a closed subset of $X$
 and $x \in X \setminus A$. If $A = \mathcal{G}\phi(A)$
 then there exists a $\phi$ Lipschitz Lyapunov function $L$ such that
$L(x)= 0$ and $L(y) = 1$ for all $y \in A$.
\end{theo}
\vspace{.5cm}

{\bfseries Proof:}  Since $A$ is closed, there exists $\epsilon
\in (0,1]$ such that $\tilde{x} = \phi(\epsilon,x) \not\in A$. Let
$B = \mathcal{G}\phi^{-1}(\tilde{x})$. Since $A$ is $\mathcal{G}
\phi $ invariant, the closed set $B$ is disjoint from $A$.
Corollary \ref{cor1.9} implies that there exists a continuous
$\mathcal{G}\phi$ Lyapunov function $\tilde{L}: X \to [0,1]$ with
$\tilde{L}$ equal zero on $B$ and equal to one on $A$. The
function we want is given by
\begin{equation}\label{eq4.41}
L(z) \quad =_{def} \quad \frac{1}{\epsilon} \ \int_{0}^{\epsilon}
\ \tilde{L}(\phi(s,z)) \ ds.
\end{equation}

%Since $\tilde{L}$ is a $\phi^I$ Lyapunov function we have for
%every $s \in [0,\epsilon]$:
%%4.37(33)
%\begin{equation}
% \tilde{L}(\phi(s,x)) \geq  \tilde{L}(x) >  \tilde{L}(\phi(\epsilon,y)) \geq  \tilde{L}(\phi(s,y))
% \end{equation}
%from which it follows that $L(x) > L(y)$.

Since $\phi^{[0,\epsilon]}(x) \subset B$ and
$\phi^{[0,\epsilon]}(y) \subset A$ for all $y \in A$, we have that
$L(x) = 0$ and $L$ equals one on $A$.

From a little change of variables we see that
\begin{equation}\label{eq4.42}
L(\phi(t,z)) \quad = \quad \frac{1}{\epsilon}
\int_{t}^{t+\epsilon} \ \tilde{L}(\phi(s,z)) \ ds.
\end{equation}
This is differentiable in $t$ with derivative
\begin{equation}\label{eq4.43}
\frac{d}{dt} L(\phi(t,z)) \quad = \quad ( \frac{1}{\epsilon}
\Delta_{\epsilon} \tilde{L})(\phi(t,z)).
\end{equation}
Since this is nonnegative, $L$ is a Lyapunov function for $\phi$.
Since this is bounded by $\frac{2}{\epsilon}, \ L$ is $\phi$
Lipschitz with constant $\frac{2}{\epsilon}$ by the Mean Value
Theorem.

$\Box$ \vspace{.5cm}

We denote by $\mathcal{L}_{\phi}$ the set of bounded $\phi$
uniform Lyapunov functions for $\phi$. We call $\mathcal{L}$  a
\emph{sufficient set of Lyapunov functions for $\phi$}
\index{sufficient set of Lyapunov functions for a semiflow}%
\index{semiflow!sufficient set of Lyapunov functions} when it is a
subset of $\mathcal{L}_{\phi}$ which is a sufficient set of
Lyapunov functions for $\phi^I$, that is,
\begin{equation}\label{eq4.44}
 \bigcap_{\mathcal{L}}\{ \leq_{L} \} \quad = \quad \mathcal{G}\phi.
\end{equation}

It follows from Theorem \ref{th4.13} that the set $\mathcal{L}_{\phi}$
itself is a sufficient set for $\phi$. Once again, if $X$ is
metrizable  we can choose a countable sufficient set $\mathcal{L}
\subset \mathcal{L}_{\phi}$.

If $\mathcal{L}_0$ is a sufficient set of Lyapunov functions for
$\phi$  then $\mathcal{L} = \{ L\circ \phi^t: t \in {\mathbb R}_+
\}$ is a $\phi^* \ +$invariant sufficient set for $\phi$. If
$\phi$ is reversible we can let $t$ vary over ${\mathbb R}$ and
obtain $\mathcal{L}$ which is $\phi^*$ invariant.

Now assume that $\phi$ is a proper semiflow and that $\mathcal{L}$
is a $\phi^* \ +$invariant sufficient set of Lyapunov functions
for $\phi$.  Let $\mathcal{A}$ be the  closed subalgebra
of $\mathcal{B}(X)$ which is generated by $\mathcal{L}$ together with the functions
of compact support. By Lemma \ref{lem4.10} this is a $\phi^* \
+$invariant subalgebra of $\mathcal{B}_{\phi}(X)$ and the semiflow
$\phi$ extends to a semiflow $\hat{\phi}$ on the associated
compactification $\hat{X}$. As usual we will regard the embedding
$j : X \to \hat{X}$ as an inclusion. We will call such a
$\hat{\phi}$ a \emph{Lyapunov function compactification} for the proper
semiflow $\phi$.\index{Lyapunov function compactification for a semiflow}%
\index{semiflow!Lyapunov function compactification}%
\index{compactification!Lyapunov function compactification for a semiflow}

We will use $|\phi|$ \index{$|\phi|$}to denote the set of fixed points of the
semiflow $\phi$.  That is,
\begin{equation}\label{eq4.45}
|\phi| \quad =_{def} \quad \{ x \in X : \phi^t(x) = x \quad
\mbox{for all} \ t \in {\mathbb R}_+ \}.
\end{equation}

\begin{theo}\label{th4.14} Let $\phi$ be a proper semiflow on  $X$ and let $\hat{\phi}$ on $\hat{X}$ be a
Lyapunov function compactification for the semiflow $\phi$.

\begin{itemize}
\item[(a)] $\hat{X}$ is a dynamic compactification for $\phi^I$ with
\begin{equation}\label{eq4.46}
\begin{split}
(X \times X) \cap \mathcal{G} \hat{\phi} \quad = \quad
\mathcal{G}\phi \hspace{2cm}\\
(X \times X) \cap \mathcal{G} (\hat{\phi}^J) \quad = \quad
\mathcal{G}(\phi^J).\hspace{1cm}
\end{split}
\end{equation}
\item[(b)]  The compact set $\hat{X} \setminus X$ is $\phi \ +$invariant and
every generalized recurrent point of $\hat{\phi}$ which lies
in $\hat{X} \setminus X$ is a fixed point for $\hat{\phi}$.  That
is,
\begin{equation}\label{eq4.47}
(\hat{X} \setminus X) \cap |\mathcal{G}\hat{\phi}^J| \qquad \subset
\qquad |\hat{\phi}|
\end{equation}
\item[(c)] If $\hat{E} \subset |\mathcal{G} (\hat{\phi}^J)| $
  is an $\mathcal{G} (\hat{\phi}^J) \cap
\mathcal{G} (\hat{\phi}^J)^{-1}$  equivalence class with $E =
\hat{E} \cap X$ then exactly one of the following three
possibilities holds:
\begin{enumerate}
\item[(i)]  $\hat{E}$ consists of a single point of $\hat{X} \setminus
X$ which is a fixed point for $\hat{\phi}$ and $E = \emptyset$.
\item[(ii)] $\hat{E}$ is the one point compactification of a noncompact
$\mathcal{G}(\phi^J) \cap \mathcal{G}(\phi^J)^{-1}$ equivalence
class $E$. That is, there is a noncompact $\mathcal{G}(\phi^J)
\cap \mathcal{G}(\phi^J)^{-1}$ equivalence class $E \subset X$
whose closure is $\hat{E}$ and $\hat{E} \setminus E$ is a
singleton which is a fixed point of $\hat{\phi}$.
\item[(iii)] $\hat{E}$ is contained in $ X$, i.e. $\hat{E} = E$, and it is a  compact
$\mathcal{G}(\phi^J) \cap \mathcal{G}(\phi^J)^{-1}$ equivalence
class.
\end{enumerate}
\item[(d)]  For $x \in X$
the $\mathcal{G} \hat{\phi} \cap \mathcal{G} \hat{\phi}^{-1}$
equivalence class of $x$ is the closure in $\hat{X}$ of its
$\mathcal{G}\phi \cap \mathcal{G}\phi^{-1}$ equivalence class.

\item[(e)] If $X$ is metrizable and/or $\phi$ is reversible then there exist Lyapunov function
compactifications for $\phi$ with the same properties.
\end{itemize}
\end{theo}

{\bfseries Proof:} Because $\phi$ is proper a Lyapunov
compactification for $\phi$ satisfies the conditions of Definition
\ref{df4.11} and so we can apply Theorem \ref{th4.12}.

By Corollary \ref{cor4.5} (d) every $\mathcal{G} (\hat{\phi}^J) \cap
\mathcal{G} (\hat{\phi}^J)^{-1}$ class  $\hat{E}$ is $\hat{\phi}$
invariant and by Theorem \ref{th1.18} it intersects $\hat{X} \setminus X$
in at most one point. So if $z \in \hat{E} \setminus X$ then for
every positive $t$ there exists $z_t \in \hat{E} \setminus X $
with $\phi^t(z_t) = z$. As $\hat{E} \setminus X$ is a singleton it
follows that  $z_t = z$ for all $t$ and so $z$ is a fixed point.

(d) is clear from the cases in (c).

Finally, if $X$ is metrizable we can begin with a countable
sufficient set $\mathcal{L}_0$. The $\phi^* \ +$invariant (or
$\phi^*$ invariant) sufficient set $\mathcal{L}$ obtained by
closing up under the action of $\phi$ yields a countably generated
algebra because by Lemma \ref{lem4.10} we need only let $t$ vary over
rational values.

$\Box$ \vspace{.5cm}

In Akin (1993) Proposition 6.3(b) it is shown that if $\phi$ is a
semiflow on a compact metric space $X$ then there exists a
residual subset $T$ of $(0,\infty)$ such that:
\begin{equation}\label{eq4.49}
\Omega \phi \quad = \quad \Omega \phi^t \qquad \mbox{for all}
\quad t \in T,
\end{equation}
with $\Omega \phi$ defined by (\ref{eq4.6}).

We conclude this section by extending this result from $\Omega$ to
$\mathcal{G}$.

\begin{lem}\label{lem4.15}  If $F $ be a closed relation on a compact  space
$X$ then
\begin{equation}\label{eq4.50}
\Omega \mathcal{G} F \circ \omega F \ = \  \mathcal{G} F \circ
\omega F \ = \  \mathcal{G} F \circ \Omega  F \ = \ \mathcal{G} F
\circ \Omega \mathcal{G} F \ = \ \Omega \mathcal{G} F.
\end{equation}
\end{lem}

{\bfseries Proof:}  The first four expressions form and increasing
sequence of relations and the fourth equals the fifth by Akin
(1993) Proposition 2.4(b).  So it suffices to show that if $(x,y)
\in \Omega \mathcal{G} F$ then $(x,y) \in \Omega \mathcal{G} F
\circ \omega F$.

From Akin (1993) Proposition 2.4(c) and induction  it follows that
$\Omega \mathcal{G} F = \Omega \mathcal{G} F \circ F^n$ for $n =
1,2,...$ and so there exists $z_n \in F^n(x)$ such that $(z_n,y)
\in \Omega \mathcal{G} F$.  If $z$ is a limit point of the
sequence $\{ z_n \}$ then $(x,z) \in \omega F$ and $(z,y) \in
\Omega \mathcal{G} F$.

$\Box$ \vspace{.5cm}

\begin{cor}\label{cor4.16} If $\phi$ is a semiflow on a compact space $X$
then for any $t \in (0,\infty)$
\begin{equation}\label{eq4.51}
\Omega \phi \ = \ \Omega (\phi^t) \qquad \Longrightarrow \qquad
\Omega \mathcal{G} (\phi^J) \ = \ \Omega \mathcal{G} (\phi^t).
\end{equation}

In particular, if $X$ is metrizable then $\Omega \mathcal{G}
(\phi^J) \ = \ \Omega \mathcal{G} (\phi^t)$ for $t$ in a residual
subset of $(0,\infty)$.
\end{cor}

{\bfseries Proof:} First observe that for $K$ any nonempty compact
subset of $[0,\infty)$
\begin{equation}\label{eq4.52}
\phi^K \circ \Omega \phi \quad = \quad \Omega \phi \quad = \Omega
(\phi^J).
\end{equation}
with $\Omega \phi$ defined by (\ref{eq4.6}).

By (\ref{eq4.10})  applied to the semiflow $(s,x) \to \phi(ts,x)$ we have
for all $t \in (0,\infty)$ that
\begin{equation}\label{eq4.53}
\mathcal{G} (\phi^{[t,2t]}) \quad = \quad \mathcal{G} (\phi^t)
\circ \phi^{[0,t]}.
\end{equation}

By (\ref{eq4.9}) applied to this semiflow as well as the original one, we
have
\begin{equation}\label{eq4.54}
\phi^I \cup \mathcal{G}(\phi^J) \quad = \quad \mathcal{G}\phi
\quad = \quad \phi^{[0,t]} \cup \mathcal{G}(\phi^{[t,2t]}).
\end{equation}

Now we put all these together.

First:
\begin{equation}\label{eq4.55}
\Omega \mathcal{G}(\phi^J) \ = \ \mathcal{G}(\phi^J) \circ
\Omega(\phi^J) \ = (\phi^I \cup \mathcal{G}(\phi^J))\circ
\Omega(\phi^J) \ = \ \mathcal{G}\phi \circ \Omega \phi.
\end{equation}
The first equation from (\ref{eq4.50}) with $F = \phi^J$ and the third
from (\ref{eq4.54}) and (\ref{eq4.52}). The second follows from (\ref{eq4.52}) and the
inclusion of $\Omega (\phi^J)$ in $\Omega \mathcal{G}(\phi^J)$.

Continuing:
\begin{equation}\label{eq4.56}
\begin{split}
\mathcal{G}\phi \circ \Omega \phi \ = \ (\phi^{[0,t]} \cup
\mathcal{G}(\phi^{[t,2t]}))\circ \Omega \phi  \ = \
\mathcal{G}(\phi^{[t,2t]}) \circ \Omega \phi \\ = \ \Omega
\mathcal{G}(\phi^t) \circ \phi^{[t,2t]} \circ \Omega \phi \ = \
\mathcal{G}(\phi^t) \circ \Omega \phi.\hspace{1cm}
\end{split}
\end{equation}
Here we first use the same argument on the new semiflow and then
apply (\ref{eq4.53}) and (\ref{eq4.52}) again.

When $\Omega \phi = \Omega (\phi^t)$ then this last equals $\Omega
\mathcal{G}(\phi^t)$ by (\ref{eq4.50}) with $F = \phi^t$.

$\Box$\vspace{.5cm}

\section{Chain Compactifications}\label{secchain}

The chain relations are uniform space notions.  This is not so
apparent in the compact case because a compact space has a unique
uniformity consisting of all neighborhoods of the diagonal. Throughout
this section the  spaces $X$, still locally compact and $\sigma$-compact,
are assumed to be  uniform spaces
with uniformity denoted $\mathcal{U}_X$. \index{$\mathcal{U}_X$} We let  $\B_{\U}(X)$ \index{$\B_{\U}(X)$} denote
the closed subalgebra of uniformly continuous functions
in $\B(X)$. We use the development of uniform spaces in Kelley (1955) Chapter 6.

When we speak of different uniformities on a space $X$ we refer only to uniformities compatible with
the given topology.

The \emph{gage} \index{gage} of $\U_X$ is the set of all pseudo-metrics on $X$ which
are uniformly continuous on $X$.  For example, if $u \in \B_{\U}$ then $d_u$ defined
by $d_u(x,y) = |u(x) - u(y)|$ \index{$d_u$} is in the gage. The gage generates the uniformity in
the sense that for every $V \in \U_X$ there exists $d$ in the gage such that $V \supset V^d$
where $V^d = \{ (x,y) : d(x,y) < 1 \}$, \index{$V^d$} see Kelley (1955) Theorem 6.19.

A uniform space is compact iff it is complete and totally bounded, see Kelley (1955) Theorem 6.32.
If $X$ is compact then the unique uniformity on $X$ is the set of all neighborhoods of the diagonal,
the gage is the set of all pseudo-metrics continuous on $X \times X$ and of course $\B_{\U}(X) = \B(X)$.
In the compact case the pseudo-metrics of the form $d_u$ with $u \in \B_{\U}(X)$ generate the uniformity
but in general this is not true.

\begin{prop}\label{prop7.1new} For a uniform space $X$ assume that $\A$ is
a closed subalgebra of $\B_{\U}(X)$ which distinguishes points and
closed sets. Let $\U(\A)$ \index{$\mathcal{U}_A$}denote the uniformity generated by $\A$, i.e. the smallest uniformity which
contains $V^{d_u}$ for all $u \in \A$. Let $\X \supset X$ be proper compactification of $X$ associated
with the subalgebra $\A$.
\begin{enumerate}
\item[(a)] The uniformity $\U(\A)$ induces the original topology on $X$ and is coarser than
the original uniformity, i.e. $\U(\A) \subset \U_X$.
\item[(b)] The uniformity induced on $X$ from $\U_{\X}$ is $\U(\A)$.
The inclusion map $j: X \to \X$  is uniformly continuous with
respect to $\U_X$.
\item[(c)] The uniformity $\U(\A)$ is totally bounded and  $\X$ is the $\U(\A)$ completion of $X$.
\item[(d)] $\U(\B_{\U}(X)) = \U_X$ iff $\U_X$ is totally bounded.
\end{enumerate}
\end{prop}

{\bfseries Proof:} (a): To say that $\U(\A)$ is generated by $\A$ means that
 $V \in \U(\A)$ iff there is a finite subset $F \subset \A$
such that $V \supset \bigcap_{u \in F} \ V^{d_u}$.  Because $\A$ distinguishes points and closed sets
$\U(\A)$ is a uniformity on $X$ with the original topology. Since $\A \subset \B_{\U}(X)$,
 $\U(\A) \subset \U_X$.

 (b):Because $\X$ is the compactification associated with $\A$ we have $j^*\B(\X) = \A$ or, equivalently,
 $\B(\X) = \{ \hat u : u \in \A \}$. By uniqueness of the uniformity on $\X$, we have $\U(\B(\X)) = \U_{\X}$.
 Clearly, $V^{d_{\hat u}} \cap (X \times X) = V^{d_u}$.
 Thus, $\U(\A)$ is the uniformity on $X$ which is induced from $\U_{\X}$.
 Since $\U(\A)$ is coarser than $\U_X$ it follows that $j$ is uniformly continuous with respect to $\U_X$.

 (c): Since $\X$ is compact the induced uniformity $\U(\A)$ on the subset $X$ is totally bounded. Since
 $X$ is dense in the complete space $\X$ the latter is the completion of the former.

 (d): By (c) $\U(\B_{\U}(X))$ is totally bounded. So if $\U_X = \U(\B_{\U}(X))$ then $\U_X$ is totally bounded.
 On the other hand, let $\tilde X$ be the $\U_X$ completion of $X$. Because $j$ is uniformly continuous and
 $\X$ is complete, $j$  extends to a uniformly continuous map $\tilde j : \tilde X \to \X$.

 If $x \not= y \in \tilde X$
 then since $\tilde X$ is Hausdorff there exists an element of $\B_{\U}(\tilde X)$ which is 0 on a neighborhood of $x$
 and $1$ on a neighborhood of $y$.  Let $u \in \B_{\U}(X)$ be the restriction to $X$. If a net $\{ z_i \}$ in $X$
 converges to $x$ in $\tilde X$ then $\{ u(z_i) \}$ is eventually $0$ and so $\hat u( \tilde j(x)) = 0$.  Similarly,
 $\hat u(\tilde j(y)) = 1$. It follows that $\tilde j$ is injective.

 If $\U_X$ is totally bounded then $\tilde X$ is compact and so $\tilde j : \tilde X \to \X $ is a homeomorphism
 between two compact spaces and so is a uniform isomorphism. The uniformities $\U_{\tilde X}$ and $\U_{\X}$
 are equal and restrict on $X$ to $\U_X$ and $\U(\B_{\U}(X))$ respectively.

$\Box$ \vspace{.5cm}

For $f$ a closed relation on  a uniform space $X$  we define:
\index{$\mathcal{C} f$}
\begin{equation}\label{eq7.1}
\mathcal{C}f \quad =_{def} \quad \bigcap_{V \in \mathcal{U}_X} \
\mathcal{O} (V \circ f \circ V).
\end{equation}
When we need to indicate the uniformity $\U_X$ explicitly, we will write $\mathcal{C}_{\U_X} f$ .

Notice that $W \circ W \subset V$ implies that the closure of $W
\circ f \circ W$ is contained in $V \circ f \circ V$.  Hence,

\begin{equation}\label{eq7.2}
\mathcal{C}f \quad = \quad \bigcap_{V \in \mathcal{U}_X} \
\mathcal{N} (V \circ f \circ V).
\end{equation}
It follows that $\mathcal{C}f$ is a closed, transitive relation
which contains $f$.  Hence,

\begin{equation}\label{eq7.3}
\mathcal{G}f \subset \mathcal{C}f.
\end{equation}
 Clearly, we also have
have

\begin{equation}\label{eq7.4}
\mathcal{C}(f^{-1}) \quad = \quad (\mathcal{C}f)^{-1},
\end{equation}
and so we can omit the parentheses. Furthermore,
$\mathcal{C}\mathcal{C}f = \mathcal{C}f$ is an easy exercise using
the definition (\ref{eq7.1}).

\begin{prop}\label{prop7.1}  If  $f$ is a + proper relation on a uniform space $X$ then

\begin{equation}\label{eq7.5}
\mathcal{C}f \quad = \quad \bigcap_{V \in \mathcal{U}_X} \
\mathcal{O} (V \circ f )\quad = \quad \bigcap_{V \in \mathcal{U}_X} \
\mathcal{N} (V \circ f ).
\end{equation}
Furthermore,

\begin{equation}\label{eq7.6}
\mathcal{C}f \quad = \quad f \cup (\mathcal{C}f) \circ f,
\end{equation}
and if, in addition, $f$ is proper then

\begin{equation}\label{eq7.7}
\mathcal{C}f \quad = \quad f \cup f \circ \mathcal{C}f.
\end{equation}
\end{prop}

{\bfseries Proof:}  Fix $x \in X$ and let $W$ vary over closed
members of the uniformity such that $W(x)$ is compact.  By
Proposition \ref{prop1.2}(f)  the intersection of the $f(W(x))$'s is $f(x)$.
For $V \in \mathcal{U}$ choose $V_1 \in \mathcal{U}$ such that
$V_1\circ V_1 \circ V_1 \subset V$. By Proposition \ref{prop1.2}(g)(applied to $f^{-1}$) there exists a $W$
such that $ f(W(x)) \subset V_1(f(x))$.  We can assume that $W
\subset V_1$ and so $ (W \circ W)(f(W(x))) \subset
V(f(x))$.

It follows that

\begin{equation}\label{eq7.8}
\begin{split}
\bigcap_W \ \mathcal{O}(W \circ f \circ W)(x) \ \subset \ \bigcap_V \ \mathcal{O}( V
\circ f )(x)\\
\subset \ \bigcap_V \ \mathcal{N}( V
\circ f )(x) \ \subset \  \bigcap_V \  \mathcal{O}( V
\circ f \circ V)(x).
\end{split}
\end{equation}
proving (\ref{eq7.5}).  We then have

\begin{equation}\label{eq7.9}
\mathcal{C}f(x) \ \subset \ \bigcap_V \ (V(f(x)) \cup \mathcal{N}( V
\circ f  )(f(x)),
\end{equation}
Since $f(x)$ is compact, Proposition \ref{prop1.2}(f) implies that
$\mathcal{C}f $ is contained in $f \cup (\mathcal{C}f)\circ f$.
The reverse inclusion follows from transitivity of $\mathcal{C}f
\supset f$.

For (\ref{eq7.7}) apply (\ref{eq7.6}) to $f^{-1}$ and invert.

$\Box$ \vspace{.5cm}

{\bfseries Remark:} In particular, this theorem applies to any closed relation $f$ when the
space $X$ is compact (with its unique uniformity). Thus, the compact space definition
for $\mathcal{C} f$ given in Section \ref{seccomp} is consistent with the
more general one used here.
\vspace{.5cm}

\begin{df}\label{df7.8new} Let $f$ be a closed relation on a uniform space $X$ and let $(\X, \f)$ be
a proper compactification of the dynamical system $(X,f)$. We call $(\X, \f)$ a \emph{chain dynamic
compactification} \index{chain dynamic compactification}%
\index{compactification!chain dynamic compactification}
of $(X,f)$ when it satisfies
\begin{itemize}
\item The inclusion $j : X \to \X$ is uniformly continuous, or, equivalently, the uniformity
$\U_X$ is finer than the restriction to $X$ of the uniformity $\U_{\X}$.
\item $(X \times X) \cap \mathcal{C} \hat{f} \quad = \quad \mathcal{C} f$.
\end{itemize}
\end{df}
\vspace{.5cm}

\begin{theo}\label{th7.8} Let $\hat{X}$ be a uniform space with $X$ a dense
subset equipped with the induced uniformity. Assume that $f$ is a closed
relation on $X$.   If  $\hat{f}$ is the closed relation on $\hat{X}$
which is  the $\hat{X} \times \hat{X}$ closure of $f$,
then

\begin{equation}\label{eq7.20}
(X \times X) \cap \mathcal{C} \hat{f} \quad = \quad \mathcal{C} f.
\hspace{2cm}
\end{equation}
\end{theo}

{\bfseries Proof:}  For $x,y \in X$ and $V$ open in
$ \mathcal{U}_{\hat{X}}$ a $V$ chain for $\hat{f}$ from $x$ to $y$
is a sequence $a_1, b_1, c_1,....,a_n, b_n, c_n,a_{n+1}$ in
$\hat{X}$ with $a_1 = x, a_{n+1} = y$, and $(a_i,b_i),
(c_i,a_{i+1}) \in V, (b_i, c_i) \in \hat{f}$ for $i = 1,...,n$.
Assume, inductively, that the terms up to and including $a_i$ lie
in $X$. First, perturb $(b_{i}, c_i)$ to a point in $f$ close
enough that $b_i \in V(a_i)$ and $c_i \in V^{-1}(a_{i+1})$.  If $i
= n$  then we have a $V$ chain in $X$ for $f$ from $x$ to $y$. If
$i < n$ then perturb $a_{i+1}$ to a point in $X \cap V(c_i) \cap
V^{-1}(b_{i+1})$.

$\Box$ \vspace{.5cm}

\begin{cor}\label{cor7.9new} Let $f$ be a closed relation on a uniform space
$X$ with $\U_X$ totally bounded. If $\X$ is the completion of $X$ and $\f$
is the $\X \times \X$ closure of $f$, then $(\X, \f)$ is a chain dynamic compactification
of $(X,f)$.
\end{cor}

{\bfseries Proof:}  Since $X$ is totally bounded, the completion $\X$ is compact. In any case
the uniformity on the completion restricts to the original uniformity on $X$.  Hence, the
inclusion is uniformly continuous and so the compactification is chain dynamic by Theorem \ref{th7.8}.

$\Box$ \vspace{.5cm}

For $A, B \subset X$ we write
\index{$\subset \subset_u$}
\begin{equation}\label{eq7.10}
A \ \subset \subset_u \ B \qquad \Longleftrightarrow \qquad V(A) \
\subset \ B \quad \mbox{for some} \ V \in \mathcal{U}_X,
\end{equation}
i.e. $B$ is a uniform neighborhood of $A$. Clearly, $A \subset
\subset_u B$ implies $A \subset \subset B$ and $\overline{A}
\subset \subset_u B^{\circ }$.  If $A \subset \subset_u B$ then $X
\setminus B \subset \subset_u X \setminus A$.

\begin{prop}\label{prop7.2} For $A, B$ subsets of a uniform space $X$, the
uniform inclusion $A \subset \subset_u B$ holds iff there exists a
uniformly continuous function $L : X \to [0,1]$ such that

\begin{equation}\label{eq7.11}
A \ \subset \ L^{-1}(1) \ \subset \ L^{-1}((0,1]) \subset B.
\end{equation}
We will say that such a function $L$ \emph{establishes the uniform
inclusion} $A \subset \subset_u B$.

If, in addition, $C \subset \subset_u B \setminus A$ then $L$ can be
chosen with $C \subset L^{-1}(0,1)$.
\end{prop}

{\bfseries Proof:} The existence of such a function implies the
uniform inclusion by uniform continuity.  For the converse we can
choose a pseudo-metric $d$ in the gage of the uniformity  so that
$V^{d}(A) \subset B$.  Define

\begin{equation}\label{eq7.12}
L(x) \quad =_{def} \quad d(x,X \setminus B)/(d(x,A) + d(x,X
\setminus B)).
\end{equation}
Since the denominator is at least 1, it follows that $L$ is
uniformly continuous.

If $C \subset \subset_u B \setminus A$ then we can choose $d$ so that,
in addition, $V^{d}(C) \subset B \setminus A$. Thus, for $x \in C$,
$d(x,A), d(x,X \setminus B) \geq 1$ and so $L(x) \in (0,1)$.

$\Box$ \vspace{.5cm}

For a closed relation $f$ on $X$ we call a closed set $U$
\emph{uniformly inward} \index{uniformly inward set} for $f$ if

\begin{equation}\label{eq7.13}
f(U) \subset \subset_u U.
\end{equation}
Clearly, a uniformly inward set is $\mathcal{C}f$ + invariant.
Also, $X \setminus U^{\circ} = \overline{X \setminus U} $ is
uniformly inward for $f^{-1}$.

If $L : X \to [0,1]$ establishes this inclusion then

\begin{equation}\label{eq7.14}
(x,y) \in f \qquad \Longrightarrow \qquad L(x) = 0 \quad \mbox{or}
\quad L(y) = 1.
\end{equation}
We will call  a uniformly continuous  function which
satisfies (\ref{eq7.14}) an \emph{elementary uniform Lyapunov function}\index{elementary uniform Lyapunov function}%
\index{Lyapunov function!elementary uniform Lyapunov function} for
$f$.  For an elementary uniform Lyapunov function $L$ the closed
sets $L^{-1}([\epsilon,1])$ are uniformly inward for every $\epsilon \in
(0,1)$. A fortiori these sets are $\mathcal{C} f$ +invariant and so $L$ is a $\mathcal{C}f$ Lyapunov function.
In fact,  $L$ is an elementary uniform Lyapunov function for $\mathcal{C} f$.

For if $(x,y) \in \mathcal{C} f$ then (\ref{eq7.6}) implies that $(x,y) \in f$ or there exists $z \in X$
such that $(x,z) \in f$ and $(z,y) \in \mathcal{C}f$. Hence, $L(x) = 0$ or $L(y) \geq L(z) = 1$ and so
$L(y) = 1$. 

For a function $L : X \to [0,1]$ define $\eta L : X \times X \to [0,1] $ by $\eta L(x,y) = L(x)(1 - L(y))$.
Thus, $L$ is an elementary uniform Lyapunov function when it is a uniformly continuous function such that
\begin{equation}\label{eq7.14newb}
f \quad \subset \quad \eta L^{-1}(0). \hspace{2cm}
\end{equation}
Note that for any $L : X \to [0,1]$ we have
\begin{equation}\label{eq7.14newbb}
 \eta L^{-1}(0) \quad \subset \quad \leq_L. \hspace{2cm}
\end{equation}

\begin{theo}\label{th7.3}  Let $f$ be a closed relation on a uniform space $X$
and let $A, B$ be compact subsets of $X$.
If  $B \cap \mathcal{C} f(A) = \emptyset$, then there exists an elementary
uniform Lyapunov function $L$ such that

\begin{equation}\label{eq7.15}
\begin{split}
x \in B \cup \mathcal{C}f^{-1}(B) \qquad \Longrightarrow \qquad L(x) = 0; \hspace{.7cm}\\
x \in  \mathcal{C}f(A) \qquad \Longrightarrow \qquad L(x) =
1. \hspace{1cm}
\end{split}
\end{equation}

If, in addition, $B \cap A = \emptyset$ then $L$ can be chosen so that $L(x) = 1$ for all $x \in A$.

If, in addition, $B, A, \mathcal{C} f(A)$ are all pairwise disjoint then $L$ can instead be chosen so that
$0 < L(x) < 1$ for all $x \in A$.
\end{theo}

{\bfseries Proof:}  By (\ref{eq7.2}) and Proposition \ref{prop1.2}(h)  there exists
an open element $V$ of the uniformity such that
from $\mathcal{N}(V \circ f \circ V)(A) \subset X \setminus B$.  Choose $W \in \U_X$
 with $W \circ W \subset V$
and let $U =_{def}  \mathcal{O}(V \circ f \circ V)(A) $ .
Clearly, $W \circ \mathcal{O}(W \circ f \circ W) \circ W(U) \subset U$. Hence,
$\overline{U}$ is uniformly inward. In fact,
\begin{equation}\label{eq7.15newb1}
\mathcal{C} f (A) \cup f(\overline{U}) \quad \subset \quad
\mathcal{O}(W \circ f \circ W)(\overline{U}) \quad \ \subset \subset_u \ U .
\end{equation}
Let $L : X \to [0,1]$ establish the latter inclusion. $L$ is an
 elementary uniform Lyapunov function with $L = 1$ on $\mathcal{C} f(A)$.  Since \\ 
 $U \subset X \setminus B, \ L = 0$ on $B$.

 Since  $L$ is a $\mathcal{C}f$ Lyapunov function, it
follows that $\mathcal{C}f^{-1}(B) \subset L^{-1}(0)$.

If $B$ is disjoint from $A$  then by compactness we can choose $V$ so that $B$ is disjoint
from $V(A) \cup \mathcal{N}(V \circ f \circ V)(A)$.  With $W $
such that $W \circ W \circ W \subset V$
let $U =_{def} (W \cup \mathcal{O}(V \circ f \circ V)(A) $ .
Again $W \circ \mathcal{O}(W \circ f \circ W)\circ W(U) \subset U$. Proceed
as before. This time we have
\begin{equation}\label{eq7.15newb2}
A \cup \mathcal{C} f (A) \cup f(\overline{U}) \ \subset \
A \cup \mathcal{O}(W \circ f \circ W)(\overline{U}) \ \subset \subset_u \ U \ \subset \ X \setminus B.
\end{equation}

If, on the other hand, $B, A$ and $\mathcal{C} f(A)$ are pairwise disjoint then  we can choose
$V_1$ a closed element of the uniformity so that $B, V_1(A)$ are disjoint compacta, both disjoint from the closed set
$\mathcal{C} f(A)$. Then use Proposition \ref{prop1.2}(h) again to
choose $V \subset V_1 $ so that $\mathcal{N}(V \circ f \circ V)(A) \subset
X \setminus (B \cup V_1(A))$. Use (\ref{eq7.15newb1})
again and apply Proposition \ref{prop7.2} with $C = A$ to get $L(x) \in (0,1)$ for $x \in A$.

$\Box$ \vspace{.5cm}

\begin{cor}\label{cor7.4} If $f$ is a closed relation on a uniform space $X$ then

\begin{equation}\label{eq7.16}
1_X \cup \mathcal{C}f \quad = \quad \bigcap_{L} \ \leq_L,
\hspace{2cm}
\end{equation}
and
\begin{equation}\label{eq7.16newb}
|\mathcal{C} f| \quad = \quad \bigcap_{L} \  L^{-1}\{0,1\}, \hspace{2cm}
\end{equation}
where $L$ varies over all elementary uniform Lyapunov functions
for $f$.

If $X$ is second countable then there is a countable
collection of elementary uniform Lyapunov functions with the same
intersections.

\end{cor}

{\bfseries Proof:}  For $x,y$ distinct points of $X$ with $(x,y)
\not\in \mathcal{C}f$, we can apply Theorem \ref{th7.3} with $A = \{ x \}$
and $B = \{ y \}$ to obtain $L$ with $L(x)=1 $ and $ L(y)= 0$. Let $L_{x,y} =
2 (max(min(L,\frac{3}{4}),\frac{1}{4}) - \frac{1}{4})$.  Thus,
$L_{x,y}$ is an elementary uniform Lyapunov function with $U_{x,y}
= L_{x,y}^{-1}(1) \times L_{x,y}^{-1}(0)$ a neighborhood of
$(x,y)$. If $X$ is second countable then we can choose a countable
number of such functions so that the $U_{x,y}$ cover the
Lindel\"{o}f complement of $1_X \cup \mathcal{C}f $ in $X \times
X$.

If $x \not\in |\mathcal{C} f|$ then $B, A$ and $\mathcal{C}f(A)$ are pairwise
disjoint and so by the last case of Theorem \ref{th7.3} we can choose $L$ so
that $L(y) = 0 < L(x) < 1$ with $L = 1$ on $\mathcal{C} f(x)$.

If $L$ is an elementary uniform Lyapunov function and $L(x) \in (0,1)$ then
$L = 1$ on $\mathcal{C} f(x)$ and $L = 0$ on $\mathcal{C} f^{-1}(x)$. Hence,
$x \not\in |\mathcal{C} f|$.  Again when $X$ is second countable we can
cover $X \setminus |\mathcal{C} f|$ with a countable collection of the open
sets $L^{-1}(0,1)$.

$\Box$ \vspace{.5cm}

%
%\begin{cor}\label{cor7.5} If $f$  is a + proper closed relation on a
%locally compact uniform space $X$ and $x \not\in |\mathcal{C} f|$ then
%there exists a compact neighborhood $U$ of $x$ and an elementary uniform Lyapunov
%function $L$ for $f$ such that $L = 0$ on $U \cup \mathcal{C} f^{-1}(U)$ and $L = 1$ on $\mathcal{C} f(U)$.
%In particular, $L$ is a splitting Lyapunov function for $f$.
%\end{cor}
%
%{\bfseries Proof:}  By Lemma 1.12 there exists a compact neighborhood $U$ of $x$ such that $U, \mathcal{C}f(U)$
%and $\mathcal{C}f^{-1}(U)$ are pairwise disjoint. Let $A = f(U)$, which is compact because $f$ is
%+ proper, and let $B = U$.  By Theorem 7.3 there is an elementary uniform Lyapunov function $L$ which is $0$
%on $U \cup \mathcal{C}f^{-1}(U)$ and is $1$ on $A \cup \mathcal{C}f(A)$ which equals $\mathcal{C}f(U)$ by (\ref{eq7.6}).
%
%$\Box$ \vspace{.5cm}

\begin{df}\label{df7.6} Let $X$ be a  uniform space and $f$ be a
 closed relation on $X$. A collection $\L$ of elementary
uniform Lyapunov functions is called $\mathcal{L}$ \emph{a
sufficient set of elementary uniform Lyapunov functions} \index{sufficient set of elementary uniform Lyapunov functions}
when

\begin{equation}\label{eq7.17}
1_X \cup \mathcal{C}f \quad = \quad \bigcap_{L \in \mathcal{L}} \
\leq_L. \hspace{2cm}
\end{equation}
\end{df}
\vspace{.5cm}

From Corollary \ref{cor7.4} it follows that the collection of all elementary uniform
Laypunov functions for $f$ is a sufficient set and if $X$ is second countable then a countable
sufficient set exists. Notice that $X$ second countable is equivalent to the assumption that
the underlying space is metrizable.  We do not need that the uniformity $\U_X$ is metrizable.
For example, if $X$ is a noncompact metric space then the uniformity of all neighborhoods of the
diagonal does not have a countable base and so is not metrizable.

If $\mathcal{U}_1, \mathcal{U}_2$ are uniformities on $X$ then we
denote by $\mathcal{C}_{\U_1}$ and $\mathcal{C}_{\U_2}$ the corresponding
chain operators.  If $\mathcal{U}_1 \supset \mathcal{U}_2$ then
for any relation $f$ on $X$ $\mathcal{C}_{\U_1} f \subset \mathcal{C}_{\U_2}
f$.  That is, the coarser uniformity has the larger chain
relation.

\begin{theo}\label{th7.7}  Assume that $X$ is a  uniform space
with uniformity $\mathcal{U}_X$ and that $f$ is a +proper
relation on $X$. Let $\mathcal{A}$ be a closed
subalgebra of $\B_{\U}(X)$ which distinguishes points and closed sets
and let $(\X, \f)$ be the associated proper compactification of $(X,f)$.
Assume that $\A$ contains
 a sufficient set of elementary
uniform Lyapunov functions for $f$ then

\begin{enumerate}
\item[(a)]$(\X, \f)$ is a chain dynamic compactification, i.e. the inclusion
$j : X \to \X$ is uniformly continuous and
\begin{equation}\label{eq7.18}
(X \times X) \cap \mathcal{C} \hat{f} \quad = \quad \mathcal{C} f. \hspace{2cm}
\end{equation}
\item[(b)] If $C$ is a compact $\mathcal{C} f$ unrevisited subset of $X$
then $C$ is a $\mathcal{C} \f$ unrevisited subset of $\X$.
That is, $ \mathcal{C} \f (C) \cap \mathcal{C} \f^{-1}(C) \subset C$.
\item[(c)] If $C$ is a compact $\mathcal{C} f$ +invariant subset of $X$
then $C$ is a $\mathcal{C} \f$ +invariant subset of $\X$.
That is, $ \mathcal{C} \f (C)  \subset C$.
\item[(d)] If $\hat{E} \subset |\mathcal{C} \hat{f}| $ is
 a $\mathcal{C} \hat{f} \cap
\mathcal{C} \hat{f}^{-1}$ equivalence class with $E = \hat{E} \cap
X$, then exactly one of the following four possibilities holds:
\begin{itemize}
\item[(i)]  $\hat{E} \subset \hat{X} \setminus
X$ and $E = \emptyset$.
\item[(ii)] $E$ is contained in $ |\mathcal{C}f|$  and is a noncompact
$\mathcal{C} f \cap \mathcal{C}f^{-1}$ equivalence class with
$\hat{E} \setminus E \not= \emptyset$.
\item[(iii)] $\hat{E} = E$ is contained in $|\mathcal{C}f|$ and is a compact
$\mathcal{C} f \cap \mathcal{C} f^{-1} $ equivalence class.
\end{itemize}

\item[(e)] If $x,y \in |\mathcal{C}f|$  lie in distinct $\mathcal{C}f \cap
\mathcal{C}f^{-1}$ equivalence classes then their equivalence
classes have disjoint closures in $\hat{X}$.

\item[(f)] If $L \in \mathcal{A}$ then $L$ is an elementary uniform Lyapunov function for $f$ iff
$\hat L$ is an elementary uniform Lyapunov function for $\f$.
\end{enumerate}
\end{theo}

{\bfseries Proof:}  (a): By definition
\begin{equation}\label{eq7.18new}
\mathcal{C} f = \mathcal{C}_{\U_X} f \qquad \mbox{ and } \qquad
\mathcal{C} \f = \mathcal{C}_{\U_{\X}} \f.
\end{equation}

Since $\mathcal{U}_{\U(\A)}$ is coarser than
$\mathcal{U}_X$ it follows that $j$ is uniformly continuous and $\mathcal{C}_{\U_X} f \subset
\mathcal{C}_{\U(\A)} f$.  Because the functions of $\mathcal{L}$ are, by definition,
$\mathcal{U}(\A)$ uniformly continuous, they are elementary Lyapunov functions for $f$ which
are uniform with respect $\mathcal{U}(\A)$.  Since $\L$ is a $ \mathcal{C}_{\U_X} f$ sufficient set we have

\begin{equation}\label{eq7.19}
1_X \cup \mathcal{C}_{\U(\A)} f \quad   \subset \quad \bigcap_{L \in \L} \ \leq_L \quad = \quad
1_X \cup \mathcal{C}_{\U_X} f.
\end{equation}

We must show
that $(x,x) \in \mathcal{C}_{\U(\A)} f $ implies $(x,x) \in \mathcal{C}_{\U_X}
f$.

If $(x,x) \in f$ then $f \subset \mathcal{C}_{\U_X} f$ implies $(x,x)
\in \mathcal{C}_{\U_X} f$.

Assume $(x,x) \not\in f$.  Since $(x,x) \in \mathcal{C}_{\U(\A)} f,$
(\ref{eq7.6}) implies there exists $y$ such that $(x,y) \in f$ and $(y,x) \in
\mathcal{C}_{\U(\A)} f$.  Since $(x,x) \not\in f$ it follows that $x
\not= y$ and so by (\ref{eq7.19}) $(x,y), (y,x) \in \mathcal{C}_{\U_X} f$.  By
transitivity $(x,x) \in \mathcal{C}_{\U_X} f$.

Thus we have
\begin{equation}\label{eq7.20new}
 \mathcal{C} f \quad = \quad  \mathcal{C}_{\U_X} f \quad = \quad  \mathcal{C}_{\U(\A)} f. \hspace{2cm}
\end{equation}

By Proposition \ref{prop7.1new}
$\U_{\X}$ induces the uniformity $\U(\A)$ on $X$. So
by Theorem \ref{th7.8}  we obtain
\begin{equation}\label{eq7.21new}
(X \times X) \cap \mathcal{C} \f \quad = \quad (X \times X) \cap
\mathcal{C}_{\U_{\X}} \f \quad = \quad \mathcal{C}_{\U(\A)} f.
\end{equation}

(b): Let $\hat C = C \cup (\mathcal{C}\f(C) \cap \mathcal{C}\f^{-1}(C)$. Proceed as in the proof of
Theorem \ref{th1.9new}(a). Instead of using Theorem \ref{th2.1} apply Akin (1993)Theorem 4.5 which
implies that $(\mathcal{C} \f)_{\hat C} = \mathcal{C}(\f_{\hat C})$.  That is, if $x,y \in \hat C$
and $y \in \mathcal{C} \f(x)$ then for every $V \in \U_{\X}$ there is a $V$ chain in $\hat C$ from
$x$ to $y$.

(c): Let $\hat C = C \cup \mathcal{C} \f(C)$. Proceed as in Theorem \ref{th1.9new} (b), using Akin
(1993) Theorem 4.5 in place of Theorem \ref{th2.1}.

(d): Proceed as in the proof of Theorem \ref{th1.6} (a). Apply (b) to eliminate the analogue of
case (iv) as in the proof of Corollary \ref{cor1.10new}.

(e): Follow the proof of Theorem \ref{th1.6} (c).

(f):  $L \in \A$ maps to $[0,1]$ iff $\hat L$ does. Clearly, $X \times X \cap \eta \hat L^{-1}\{ 0,1 \} =
\eta L^{-1} \{ 0,1 \}$. Hence, $f \subset \eta L^{-1}\{ 0,1 \}$ iff $\f \subset \eta \hat L^{-1} \{ 0,1 \}$.

$\Box$ \vspace{.5cm}

{\bfseries Remark:} If $f$ is closed but not + proper we still get
\begin{equation}\label{eq7.18newa}
1_X \cup \mathcal{C}_{\U_X} f \quad = \quad 1_X \cup \mathcal{C}_{\U(\A)} f. \hspace{2cm}
\end{equation}
\vspace{.5cm}

\begin{cor}\label{cor7.9}  Let  $f$ be a + proper, closed
relation on a uniform space $X$. Let $\mathcal{L} \subset \mathcal{B}_{\U}(X)$ be a sufficient
set of elementary uniform Lyapunov functions for $f$ and
 $\mathcal{A}$ be the closed subalgebra generated by
$\mathcal{L}$ and the continuous functions with compact support.
If $(\X, \f)$ is the  $\L$ compactification of the dynamical system
$(X,f)$ then $(\X, \f)$ is a chain dynamic compactification of $(X,f)$.

Furthermore, if $\hat{E} \subset |\mathcal{C} \hat{f}| $ is
 a $\mathcal{C} \hat{f} \cap
\mathcal{C} \hat{f}^{-1}$ equivalence class with $E = \hat{E} \cap
X$, then exactly one of the following four possibilities holds:
\begin{itemize}
\item[(i)]  $\hat{E}$ consists of a single point of $\hat{X} \setminus
X$.
\item[(ii)] $E$ is contained in $ |\mathcal{C}f|$  and is a noncompact
$\mathcal{C} f \cap \mathcal{C}f^{-1}$ equivalence class with
$\hat{E} \ $its one point compactification. That is, there is a
noncompact equivalence class $E \subset |\mathcal{C} f|$ whose
closure in $\hat{X}$ is $\hat{E}$ and $\hat{E} \setminus E$ is a
singleton.
\item[(iii)] $\hat{E} = E$ is contained in $|\mathcal{C}f|$ and is a compact
$\mathcal{C} f \cap \mathcal{C} f^{-1} $ equivalence class.
\end{itemize}

\end{cor}

{\bfseries Proof:}  The functions with compact support are uniformly continuous.
Hence, $\A \subset \B_{\U}(X)$. Hence, the compactification is chain dynamic by Theorem \ref{th7.7} (a).

 As with Lemma \ref{lem1.17} the elementary uniform Lyapunov functions in
$\L$ distinguish the points of $\X \setminus X$ and
so $\hat E \cap (\X \setminus X)$ contains at most one point. From (\ref{eq7.18}) we see that $E = \hat E \cap X$ is
either empty or is a $\mathcal{C} f \cap \mathcal{C} f^{-1}$ equivalence class. So the remaining results follows
from Theorem \ref{th7.7} (d).

$\Box$ \vspace{.5cm}

If $(X,f)$ is a cascade with $f$ uniformly continuous
and $\L_0$ is a sufficient set of elementary uniform Lyapunov functions for $f$ then
$\L = \{ L \circ f^n : L \in \L_0 $ and $ n \in \Z_+ \}$ is an $f^*$ + invariant sufficient set of
elementary Lyapunov functions.  If $f$ is a uniform isomorphism then we can let $n$ vary over $\Z$ to get
an $f^*$ invariant set.

\begin{cor}\label{cor7.10} Let $f$ be a proper uniformly continuous map on a uniform space $X$. For $\L$  an $f^*$
+invariant sufficient set of elementary uniform Lyapunov functions for $f$ the $\L$ compactification
$(\X,\f)$ is a chain dynamic cascade compactification of $(X,f)$ with
$|\mathcal{C} \f| \cap (\X \setminus X) \subset | \f |$. If $f$ is a homeomorphism and
$\L$ is $f^*$ invariant then $\f$ is a homeomorphism.

If $X$ is second countable then there exists a countable $f^*$ + invariant sufficient set of elementary
uniform Lyapunov functions and if $f$ is a uniform isomorphism it can be chosen $f^*$ invariant.
In that case, the space $\X$
is metrizable.
\end{cor}

{\bfseries Proof:} Since $f$ is proper the functions of compact support form an $f^*$ +invariant set
($f^*$ invariant when $f$ is a homeomorphism).
Hence, the algebra $\A$ generated by $\L$ and these functions is $f^*$ +invariant ($f^*$ invariant when
$f$ is a homeomorphism and $\L$ is $f^*$ invariant).  The $\L$ compactification is a cascade compactification
by Theorem \ref{th1.21new} (c) and is chain dynamic by
Corollary \ref{cor7.9}. As in Theorem \ref{th1.21new} (a) the
$\f$ invariance of the $\mathcal{C} \f \cap \mathcal{C} \f^{-1}$ equivalence class $\hat E$ implies
that the singleton $\hat E \cap (\X \setminus X)$ is a fixed point of $\f$.

We have seen above that $\L_0$ can be chosen countable when $X$ is second countable and so in those
cases $\A$ is countably generated and $\X$ is metrizable.

$\Box$ \vspace{.5cm}

If $U$ is a uniformly inward set for a closed relation $f$ then \\ $\bigcap_{n \in \Z_+} \ f^n(U)$
is the associated \emph{attractor}. \index{attractor} If $f$ is a homeomorphism then $\bigcap_{n \in \Z_+} \ f^n(U)
= \bigcap_{n \in \Z_+} \ f^n(U^{\circ})$ implies that an attractor is a $G_{\delta}$.

\begin{ex}\label{ex7.10} Even with $(X,f)$ a  cascade and $X$ compact, an attractor need not be
a $G_{\delta}$.\end{ex}
Let $e$ be a point in a compact space $X_0$ which does not have a countable base of neighborhoods and
so $\{ e \}$ is not $G_{\delta}$. Let $X = X_0 \times \Z_+  \cup \{ \infty \}$ be the one point compactification
of $X_0 \times \Z_+$ and define $f$ on $X$ by
\begin{equation}\label{eq7.20newb}
f(x) \quad = \quad \begin{cases} e \qquad \qquad \ \mbox{if} \quad x = (a,0), \\
(a,n-1) \quad \mbox{if} \ x = (a,n) \ \mbox{and} \ n > 0, \\
\infty \qquad \qquad \mbox{if} \quad x = \infty.
\end{cases}
\end{equation}
$U = X_0 \times \{0\}$ is inward with attractor $\{ e \}$.

$\Box$ \vspace{.5cm}

If $X$ is a compact metric space there are only countably many attractors,
see, e.g. Akin (1993) Proposition 3.8.

\begin{ex}\label{ex7.11} Even with $(X,f)$ a cascade on a metric space with $f$ a uniform isomorphism
 there may be uncountably many attractors. \end{ex}
  If $f$ is the time-one map for the flow
associated with the gradient of the function $t \mapsto cos(2 \pi
t)$ then every subset of ${\mathbb Z}$ is an attractor for the
homeomorphism $f$ and so $f$ has uncountably many attractors.

$\Box$ \vspace{.5cm}

Because a locally compact and
$\sigma$-compact space is paracompact,  the set of all neighborhoods of the
diagonal is a uniformity on $X$ with the appropriate topology, it
is clearly the finest such uniformity.  If $A \subset \subset B$
then with respect to this uniformity $A \subset \subset_u B$.  To
see this choose for each $x \in \overline{A}$ an open set $U_x$
containing $x$ and contained in $B^{\circ}$.  Together with $X \setminus
\overline{A}$ these form an open cover and by paracompactness
there exists a neighborhood $V$ of the diagonal such that $ \{
V(x) : x \in X \}$ refines this cover.  Hence, $V(A) \subset B$.
The chain relation for a closed relation $f$ associated with this
uniformity is the smallest one.
\vspace{.5cm}

\section{Stopping at Infinity}\label{secstop}

In this section we apply an idea from Beck (1958).

The usual way of obtaining a flow on a smooth manifold $X$ is by
integrating a smooth vectorfield $\xi$ on $X$.  Some boundedness
condition is necessary to avoid reaching infinity in finite time.
It suffices that $\xi$ be bounded in norm with respect to a
complete Riemannian metric. The flow $\phi$ is obtained by
integrating the differential equation:

\begin{equation}\label{eq5.1}
\frac{dx}{dt} \quad = \quad \xi (x), \hspace{2cm}
\end{equation}

If $V : X \to {\mathbb R}_+$ is a smooth function then the orbits
of the flow $\psi$ associated with the vectorfield $V \cdot \xi$
agree as oriented sets with those of the original flow except
where they are interrupted by the new fixed points introduced by
the zeroes of $V$. This is best seen by regarding the new flow as
obtained by a time-change with the new time $\tau$ related to the
original time $t$ via the time change:

\begin{equation}\label{eq5.2}
\frac{dt}{d\tau} \quad = \quad V(x)
\end{equation}
which combines with (5.1) to yield:

\begin{equation}\label{eq5.3}
\frac{dx}{d\tau} \quad = \quad V(x) \cdot \xi(x).
\end{equation}

To be a bit more precise the flows $\phi$ and $\psi$ satisfy

\begin{equation}\label{eq5.4}
\begin{split}
\psi(\tau,x) \quad = \quad \phi(t,x) \qquad \mbox{where} \\
\frac{dt}{d\tau} \quad = \quad V(\phi(t,x)).\hspace{1.5cm}
\end{split}
\end{equation}

The equations (\ref{eq5.1}) and (\ref{eq5.3}) require differentiating a path in
$X$ which only makes sense for a manifold.  But Beck observed that
the time change equation \index{time change} in (\ref{eq5.4}) uses $x$ in $X$ as a parameter
and makes perfect sense for any space $X$.  Furthermore, the whole
time change procedure works fine in this general context.

Let $\phi$ be a semiflow on  $X$ and let $V : X \to {\mathbb R}_+$
be continuous. Let $X_0 = V^{-1}(0)$, the closed zero-set of $V$.
For any $x \in X$ let

\begin{equation}\label{eq5.5}
t^*(x) \quad =_{def} \quad sup \{ t \in [0,\infty] : V(\phi(s,x))
> 0 \quad \mbox{ for all} \ s \in [0,t) \}.
\end{equation}
Thus, $t^*(x) = 0$ iff $x \in X_0$ and $t^*(x) = \infty$ iff the
solution curve $\phi(t,x)$ never hits $X_0$. The set

\begin{equation}\label{eq5.6}
\begin{split}
D_{V} \quad =_{def} \quad \{ (s,x) \in {\mathbb R}_+ \times X : s
< t^*(x) \}\quad =  \hspace{.5cm} \\
\{ (s,x) \in {\mathbb R}_+ \times X : V(\phi(s_1,x))
> 0 \quad \mbox{ for all} \ s_1 \in [0,s] \}
\end{split}
\end{equation}
is  open in ${\mathbb R}_+ \times X$ and the second
coordinate projection maps $D_V$ onto $X \setminus X_0$.

Notice that

\begin{equation}\label{eq5.7}
s < t^*(x) \qquad \Longrightarrow \qquad t^*(\phi(s,x)) = t^*(x) -
s.
\end{equation}

We now introduce the version we will need of the boundedness
condition which was required by $\xi$ above.

\begin{df}\label{df5.1}Let $\phi$ be a semiflow on  $X$.
 We call $V : X \to {\mathbb R}_+$ a \emph{$\phi$ regular function}
 \index{$\phi$ regular function} when
\begin{itemize}
\item  For every $x \in X \setminus X_0$
\begin{equation}\label{eq5.8}
\int_0^{t^*(x)} \ \frac{ds}{V(\phi(s,x))} \quad = \quad \infty
\end{equation}
\item  For every $x \in X_0$ such that $\phi(s,x) \in X \setminus
X_0$ for all $s$ in some interval $(0,\epsilon)$ with $\epsilon
>0$
\begin{equation}\label{eq5.9}
\int_{0}^{\epsilon} \ \frac{ds}{V(\phi(s,x))} \quad = \quad \infty
\end{equation}
\end{itemize}
\end{df}
\vspace{.5cm}

{\bfseries Remark:}  If $\phi$ is a reversible semiflow then
condition (\ref{eq5.9}) for $\phi$ is equivalent to (5.8) for the reverse
flow $\phi^{-1}$.  In particular, if $V$ is regular for $\phi$
then it is regular for $\phi^{-1}$ when $\phi$ is reversible.
\vspace{.5cm}

When $V$ is $\phi$ regular, the map $\bar{\tau} : D_V \to {\mathbb
R}$ given by
\begin{equation}\label{eq5.10}
\bar{\tau}(t,x) \quad =_{def} \quad \int_0^{t} \
\frac{ds}{V(\phi(s,x))}
\end{equation}
maps $[0,t^*(x))\times \{ x \}$ homeomorphically onto ${\mathbb
R}$ with $\bar{\tau}(0,x) = 0$ for every $x \in X \setminus X_0$.
We extend the inverse by defining $\bar{t} : {\mathbb R}_+\times X
\to {\mathbb R}_+$ so that for every $\tau \in {\mathbb R}_+$
\begin{equation}\label{eq5.11}
\begin{split}
\bar{t}(\tau,x) \quad = \quad 0 \qquad  \qquad \mbox{for} \ x \in X_0 \hspace{1.5cm}\\
 \tau
\quad = \quad \int_0^{\bar{t}(\tau,x)} \ \frac{ds}{V(\phi(s,x))}
\qquad \mbox{for} \ x \in X \setminus X_0.
\end{split}
\end{equation}

Thus, for each $x \in X \setminus X_0 \ \tau \mapsto
\bar{t}(\tau,x)$ takes ${\mathbb R}_+$ onto $[0,t^*(x))$.
Furthermore since $\bar{t}(\tau,x) < t^*(x)$,
$\phi(\bar{t}(\tau,x),x)$ remains in $X \setminus X_0$ for all
$\tau \in {\mathbb R}_+$.

We will call $\bar{t}$ the \emph{time change map} for $\phi$
associated with $V$.

It is easy to see that $\tau \times \pi_2 : D_V  \to
{\mathbb R}_+ \times (X \setminus X_0)$ and the restriction
$\bar{t} \times \pi_2 : {\mathbb R}_+ \times (X \setminus X_0) \to
D_V$ are inverse homeomorphisms. Continuity of $\bar{t}$ at the
points of ${\mathbb R}_+ \times X_0$ requires a more delicate
argument.

\begin{lem}\label{lem5.2} Assume that $V$ is regular for a semiflow $\phi$ on
$X$ and that $X_0 = V^{-1}(0)$. For every $ \epsilon, M >0$ and $x
\in X_0$ there exists $\delta > 0$ such that
\begin{equation}\label{eq5.12}
\int_0^{\epsilon} \ \frac{ds}{max(\delta,V(\phi(s,x)))} \quad >
\quad M.
\end{equation}
\end{lem}

{\bfseries Proof:} Case i : There exists a positive $\epsilon_1 \leq
\epsilon$ such that $V(\phi(s,x)) = 0$ for all $s \in
[0,\epsilon_1]$: For any $\delta > 0$
\begin{equation}\label{eq5.13}
\int_0^{\epsilon} \ \frac{ds}{max(\delta,V(\phi(s,x)))} \geq
\int_0^{\epsilon_1} \ \frac{ds}{max(\delta,V(\phi(s,x)))} =
\epsilon_1/\delta.
\end{equation}
Choose $\delta < \epsilon_1/M$.

Case ii : There exists $\epsilon_1 \leq \epsilon$ such that
$V(\phi(s,x)) > 0$ for all $s \in (0,\epsilon_1]$: By (5.9) there
exists $q \in (0,\epsilon_1)$ such that:
\begin{equation}\label{eq5.14}
\int_q^{\epsilon_1} \ \frac{ds}{V(\phi(s,x))} > M.
\end{equation}
Choose $\delta \leq inf \{ V(\phi(s,x)) : s \in [q,\epsilon_1]
\}$.

Case iii : There exist $0 < s_1 < s_2 < \epsilon$ such that
$V(\phi(s_1,x)) > 0$ and $V(\phi(s_2,x)) = 0$: Let $x_1 =
\phi(s_1,x) \in X \setminus X_0$. Since $\phi(s_2 - s_1,x_1) \in
X_0$ it follows that $t^*(x_1) \leq s_2 - s_1$. By (5.8) there
exists $q \in (s_1,t^*(x_1) + s_1) \subset (s_1,s_2)$ such that
\begin{equation}\label{eq5.15}
\int_{s_1}^{q} \ \frac{ds}{V(\phi(s,x))} = \int_0^{q - s_1} \
\frac{ds}{V(\phi(s,x_1))} > M.
\end{equation}
Choose $\delta \leq inf \{ V(\phi(s,x)) : s \in [s_1,q] \}$.

$\Box$ \vspace{.5cm}

\begin{theo}\label{th5.3} Let $\phi$ be a semiflow on  $X$ and $V : X \to {\mathbb R}_+$ be a $\phi$ regular function
with $X_0 = V^{-1}(0)$. The associated time change map  $\bar{t}
:{\mathbb R}_+ \times X \to {\mathbb R}_+$ is continuous.

Furthermore, for all $x \in X $ and $\tau_1, \tau_2 \in {\mathbb
R}_+$ the following cocycle condition \index{cocycle condition} holds:
\begin{equation}\label{eq5.16}
 \bar{t}(\tau_2,\phi(\bar{t}(\tau_1,x),x)) + \bar{t}(\tau_1,x) \quad
= \quad \bar{t}(\tau_2 + \tau_1,x).
\end{equation}
\end{theo}

{\bfseries Proof:} Let $\{(\tau_i,x_i) \}$ be a sequence in
${\mathbb R}_+ \times X$ converging to $(\tau,x) \in {\mathbb R}_+
\times X_0$. We must prove that $\{ \bar{t}(\tau_i,x_i) \}$
converges to $0$. If this is not true then by going to a
subsequence we can find a positive $\epsilon$ such that
$\bar{t}(\tau_i,x_i) > \epsilon$ for all $i$ and we can assume as
well that $\tau_i < 2 \tau + 1$ for all $i$ as well. From this we
will derive a contradiction.

By Lemma \ref{lem5.2} we can choose $\delta > 0$ so that with $z = x$
\begin{equation}\label{eq5.17}
\int_0^{\epsilon} \ \frac{ds}{max(\delta,V(\phi(s,z)))} \quad >
\quad 2 \tau + 1.
\end{equation}
By continuity of the integral as a function of $z$ the inequality
holds for all $z$ in some neighborhood of $x$. Hence, it holds for
$z = x_i$ when $i$ is large enough.  But
\begin{equation}\label{eq5.18}
2 \tau + 1 > \tau_i =  \int_0^{\bar{t}(\tau_i,x_i)} \
\frac{ds}{V(\phi(s,x_i))} \geq \int_0^{\epsilon} \
\frac{ds}{max(\delta,V(\phi(s,x_i)))}.
\end{equation}
This contradiction proves continuity of $\bar{t}$ at $(x,\tau)$.

Finally, (\ref{eq5.16}) says $0 + 0 = 0$  when $x \in X_0$. When $x \in
X\setminus X_0$ let  $t_1 = \bar{t}(\tau_1,x)$ and observe
\begin{equation}\label{eq5.19}
\begin{split}
\tau_2 = \int_0^{\bar{t}(\tau_2,\phi(t_1,x))} \
\frac{ds}{V(\phi(s,\phi(t_1,x)))} \\ =
\int_0^{\bar{t}(\tau_2,\phi(t_1,x))} \ \frac{ds}{V(\phi(s +
t_1,x))}
\\ = \int_{t_1}^{\bar{t}(\tau_2,\phi(t_1,x)) + t_1} \
\frac{ds}{V(\phi(s,x))}
\end{split}
\end{equation}
from which (5.16) follows because
\begin{equation}\label{eq5.20}
\tau_1 \quad = \quad \int_0^{t_1} \
\frac{ds}{V(\phi(s,x))}.\hspace{1cm}
\end{equation}

$\Box$ \vspace{.5cm}

\begin{cor}\label{cor5.4} Let $\phi$ be a semiflow on  $X$ and $V : X \to {\mathbb R}_+$ be a $\phi$ regular function
with $X_0 = V^{-1}(0)$ and $\bar{t} : {\mathbb R}_+ \times X \to
{\mathbb R}_+$ the associated time change map. If $\psi : {\mathbb
R}_+ \times X \to {\mathbb R}_+$ is defined by

\begin{equation}\label{eq5.21}
\psi(\tau,x) \quad =_{def} \quad \phi(\bar{t}(\tau,x),x).
\end{equation}
then $\psi$ is a semiflow on $X$ with $|\psi| = |\phi| \cup X_0$.
Furthermore, if $\phi$ is reversible then $\psi$ is reversible.
\end{cor}

{\bfseries Proof:} Continuity of $\psi$ follows from Theorem \ref{th5.3}
and the semigroup equation $\psi^{\tau_2} \circ \psi^{\tau_1} =
\psi^{\tau_2 + \tau_1}$ follows from the cocycle equation (5.16).

Clearly, the points of $|\phi|$ and $X_0$ are fixed by $\psi$ and
any point of $X \setminus X_0$ which is not fixed by $\phi$ is not
fixed by $\psi$.

If $\phi$ is reversible then, as was noted above, $V$ is regular
for the reverse flow $\phi^{-1}$. Extend $\phi$ to the flow $\phi
: {\mathbb R} \times X \to X$ and let
\begin{equation}\label{eq5.22}
\begin{split}
t^{**}(x) \ =_{def} \ sup \{ t \in [0,\infty] :
V(\phi(-s,x))
> 0 \quad \mbox{ for all} \ s \in [0,t) \}. \\
D^{\pm}_V \quad =_{def} \quad \{ (x,t) \in {\mathbb R} \times X :
-t^{**}(x) < t < t^*(x) \}. \hspace{2cm}
\end{split}
\end{equation}
The formula (\ref{eq5.10}) extends to define $\bar{\tau} : D^{\pm}_V \to
{\mathbb R}$ and (\ref{eq5.11}) extends to $\bar{t} : {\mathbb R} \times X
\to {\mathbb R}$. All of the maps with negative $t$ are just the
corresponding maps for the reverse flow with the signs changed. In
particular using the extended definition of $\bar{t}$ in (\ref{eq5.21})
yields the flow extension of the semiflow $\psi$.

$\Box$ \vspace{.5cm}

In order to apply these results  we must construct $\phi$
regular functions. First, some preliminary results.

For $A,B$ subsets of $X$ with disjoint closures we define the
\emph{$\phi$ distance} between them:
\begin{equation}\label{eq5.23}
\delta_{\phi}(A,B) \quad =_{def} \quad sup \{ s \in [0,1] :
\phi([0,s] \times A) \cap B = \emptyset = A \cap \phi([0,s] \times
B) \}.
\end{equation}

\begin{lem}\label{lem5.5} Let $\phi$ be a semiflow on $X$  and let  $A,B$ be subsets of $X$ with disjoint closures.
\begin{enumerate}
\item[(a)]   If $A$ is bounded then $\delta_{\phi}(A,B) > 0$.
\item[(b)]   Let $V : X \to {\mathbb R}_+$ be a continuous function
with $X_0 = V^{-1}(0) \subset B$. Assume that $V$ is bounded on $X
\setminus B$ with $V \leq k$ on $X \setminus (A \cup B)$.

If $x \in A$  then

\begin{equation}\label{eq5.24}
\int_0^{t^*(x)} \ \frac{ds}{V(\phi(s,x))} \quad \geq \quad
\delta_{\phi}(A,B)/k.
\end{equation}

If $x \in B$ with $V(\phi(x,s)) > 0$ for all $s \in (0,\epsilon]$
and $\phi(\epsilon,x) \in A$ then

\begin{equation}\label{eq5.25}
\int_0^{\epsilon} \ \frac{ds}{V(\phi(s,x))} \quad \geq \quad
\delta_{\phi}(A,B)/k.
\end{equation}
\end{enumerate}
\end{lem}

{\bfseries Proof:} (a):  Replacing the sets by their closures if
necessary, we can assume that $A$ is compact and $B$ is a disjoint
closed set. By compactness there exists a positive $s_1$ such that
$\phi([0,s_1] \times A)$ is disjoint from $B$. Now let $A_1$ be a
compact neighborhood of $A$ which is disjoint from $B$ and let
$B_1$ be the topological boundary of $A_1$, i.e. $B_1 = A_1
\setminus Int A_1$. Again there exists a positive $s_2$ such that
$\phi([0,s_2] \times B_1)$ is disjoint from $A$. Since $B_1$
separates $A$ and $B$, any path which begins in $B$ must cross
$B_1$ before it reaches $A$. Hence, $\phi([0,s_2] \times B)$ is
disjoint from $A$.  Thus, $\delta_{\phi}(A,B) \geq min(s_1,s_2) >
0$.

(b): For $x \in A$ let $t^+ $ be the first entrance time to $B$.
That is,
\begin{equation}\label{eq5.26}
t^+  \quad = \quad sup \{ t : \phi(s,x) \in X \setminus B \quad
\mbox{ for all} \ 0 \leq s < t\}.
\end{equation}
Since $X_0 \subset B$, $t^+ \leq t^*(x)$. If $t^+ = \infty$ then
the integral in (\ref{eq5.24}) in infinite because $V$ is bounded on the
complement of $B$. Hence, (\ref{eq5.24}) is true when $t^+ =
\infty.$

When $t^+$ is finite, $\phi(t^+,x) \in B$ and we let $t^-$ be the
last exit time from $A$ before $t^+$. That is,
\begin{equation}\label{eq5.27}
t^- \quad = \quad sup \{ t < t^+ : \phi(t,x) \in A \}.\hspace{2cm}
\end{equation}
By definition of $\delta_{\phi}$, $t^+ - t^- \geq
\delta_{\phi}(A,B)$.  For all $s \in (t^-,t^+)  \ \phi(s,x) \in X
\setminus (A \cup B)$ on which $V$ is bounded by $k$.  So (\ref{eq5.24})
holds in this case as well.

%(b):  Any $\phi$ path which begins in $A$ remains in $X \setminus
%B$ for all $t \in [0,\delta_{\phi}(A,B))$. Since $X_0 \subset B$
%it follows that $t^*(x) \geq \delta_{\phi}(A,B)$. These
%inequalities imply (5.24) since $V \leq k$ on $X \setminus B$.

Similarly, for (\ref{eq5.25})with $x \in B$ let $t^+$ be the entrance time
to $A$

\begin{equation}\label{eq5.28}
t^+ \quad = \quad sup \{ t : \phi(s,x) \in X \setminus A \quad
\mbox{ for all} \ 0 \leq s < t \}.
\end{equation}
By hypothesis $\epsilon \geq t^+ $ and we let $t^-$ be the exit
time from $B$

\begin{equation}\label{eq5.29}
t^- \quad = \quad sup \{ t \in [0,t^+] : \phi(t,x) \in B \}.
\hspace{2cm}
\end{equation}
As in the previous case $t^+ - t^- \geq \delta_{\phi}(A,B)$ and
for all $s \in (t^-,t^+)  \ \phi(s,x) \in X \setminus (A \cup B)$
on which $V$ is bounded by $k$.  Thus, (\ref{eq5.25}) follows.

$\Box$\vspace{.5cm}

\begin{lem}\label{lem5.6} Let $\{ A_n \}$ be a sequence of closed subsets of $X$ with
$A_0 = \emptyset$, $A_n \subset \subset A_{n+1}$ for $n = 0,1,...$
and $\bigcup_n A_n = X$.  Let $\{ a_n \}$ be a sequence in
${\mathbb R}_+$ with $a_0 = 1$ and $a_n < a_{n+1}$ for $n =
0,1,...$. There exists a continuous function $u : X \to
[1,\infty)$ with

\begin{equation}\label{eq5.30}
x \in A_{n+1} \setminus A_n \qquad \Longrightarrow \qquad a_n \
\leq \ u(x) \ \leq \ a_{n+1}.
\end{equation}
\end{lem}

{\bfseries Proof:}  Choose continuous functions $u_n : X \to I$ such that $u_n = 0$ on
$A_n$ and $= 1$ on $X \setminus A_{n+1}$. Let $u(x) = 1 +
\Sigma_{i=0}^{\infty} (a_{i+1} - a_i)\cdot u_i(x)$ which equals
$a_n + (a_{n+1} - a_n) \cdot u_{n}(x)$ if $x \in A_{n+1} \setminus
A_n$. For any $x \in X,$ there exists $n$ such that $ A_n$ is a neighborhood of $x$. On $A_n  \ u_i = 0$
for all $i > n$. Hence, $u$ is continuous.

$\Box$ \vspace{.5cm}

\begin{theo}\label{th5.7} Let $\phi$ be a semiflow on $X$ and $X_0$ be a
closed $G_{\delta}$ subset of $X$.  There exists $V : X \to [0,1] $ continuous
and $\phi$ regular with $X_0 = V^{-1}(0)$. If $X_0$ is compact
then there exists $V : X \to {\mathbb R}_+$ continuous
and $\phi$ regular with $X_0 = V^{-1}(0)$
and with $Lim_{x \to \infty} V(x) = \infty$, i.e. for every
positive real number $M$, the set $V^{-1}([0,M])$ is compact.
\end{theo}

{\bfseries Proof:} Choose a sequence $\{ A_n \}$ of closed subsets
of $X$ such that $A_0 = \emptyset$, $A_{n-1} \subset \subset A_n$
for $n = 1,2,...$ and $\bigcup_n A_n = X \setminus X_0$. If $X_0$
is compact then choose $A_1$ so that $X \setminus A_1$ is bounded
and hence $X \setminus A_n$ is bounded for every $n$. If $X_0$ is
unbounded then choose the $A_n$'s to be compact.

Let $\epsilon_0 = 1$ and for $n = 1,2,...$ inductively choose
positive

\begin{equation}\label{eq5.31}
\epsilon_n \quad \leq \quad min(\epsilon_{n-1}/2,
\delta_{\phi}(A_{n},X \setminus A_{n+1})).
\end{equation}
By Lemma \ref{lem5.5}  each $\delta_{\phi}(A_{n},X \setminus A_{n+1})$ is
positive.  Apply Lemma \ref{lem5.6} and take the reciprocal to get a
continuous $V : X \setminus X_0 \to (0,1]$ such that for

\begin{equation}\label{eq5.32}
x \in A_{n+1} \setminus A_{n} \qquad \Longrightarrow \qquad
\epsilon_{n+1}^2 \ \leq \ V(x) \ \leq \ \epsilon_n^2 .
\end{equation}
Define $V(x) = 0$ for $x \in X_0$ so that $V : X \to I$ is
continuous and $X_0 = V^{-1}(0)$.

If $x \in X \setminus X_0$ then for all $n$ large enough, $x \in
A_n$ and (\ref{eq5.24}) implies that

\begin{equation}\label{eq5.33}
\int_0^{t^*(x)} \ \frac{ds}{V(\phi(s,x))} \quad \geq \quad
\frac{1}{\epsilon_n} .
\end{equation}
As this is true for all large $n$, the integral is infinite,
proving (\ref{eq5.8}) for $V$. Similarly, for $x \in X_0$ we obtain (\ref{eq5.9})
from (\ref{eq5.25}).

This establishes the existence of a bounded $\phi$ regular
function $V$.

Now assume that $X_0$ is compact.  Letting $B_1 $ be the
closure of the complement of $A_1$ we have assumed that $X_0 \subset
\subset B_1$ with $B_1$ compact.  The $\phi$ regular function that
we constructed above, now renamed
$U : X \to [0,1]$, satisfies $X_0 = U^{-1}(0)$ and $U = 1$ on $X \setminus B_1$.

Inductively choose $B_{n+1}$ compact so that
$\phi([0,(n+1)^2],B_n) \subset \subset B_{n+1}$ and $\bigcup_n B_n
= X$.

Use Lemma \ref{lem5.6} to construct a continuous
function $V_0$ on $X$ such that $V_0 = 1$ on $B_1$ and

\begin{equation}\label{eq5.34}
x \in B_{n+1} \setminus B_n \qquad \Longrightarrow \qquad n  \
\leq \ V_0(x) \ \leq \ n+1.
\end{equation}
Since $U$ and $V_0$ are both constantly $1$ on the boundary of $B_1$ we can define the continuous
function $V$ to be $U$ on $B_1$ and $V_0$ on $X \setminus B_1$.

 Clearly, $V$ tends to infinity as
$x$ does, i.e. $V : X \to {\mathbb R}_+$ is a proper map.

Because $V$ agrees with $U$ on a neighborhood of $X_0$ condition
(\ref{eq5.9}) for $V$ follows because it holds for $U$ and the same is
true for (\ref{eq5.8}) if $t^*(x) < \infty$.

Finally, let $x \in X \setminus X_0$ with $t^*(x) = \infty$. There
exists a positive integer $N(x)$ so that  $x \in B_n$ when $n \geq
N(x)$.

For $n \geq N(x)$, let

\begin{equation}\label{eq5.35}
t^{+}_{n} \quad = \quad sup \{ t : \phi(s,x) \in X \setminus
B_{n+1} \quad \mbox{for all} \ s \in [0,t) \}.
\end{equation}
If for some $n \geq N(x)$ \ $t^{+}_{n} = \infty$ then
$V(\phi(s,x)) \leq n+1$ for all $s$ and so

\begin{equation}\label{eq5.36}
\int_0^{t^*(x)} \ \frac{ds}{V(\phi(s,x))} \quad = \quad
\int_0^{\infty} \ \frac{ds}{V(\phi(s,x))}
\end{equation}
is infinite.

Otherwise,we can define for $n \geq N(x)$

\begin{equation}\label{eq5.37}
 t^{-}_n \quad = \quad sup \{ t \in
[0,t^{+}_{n}] : \phi(t,x) \in B_n \}.
\end{equation}
By construction, $t^{+}_n - t^{-}_n \geq (n+1)^2$ and
$V(\phi(s,x)) \leq n+1$ for $s \in (t^{-}_n,t^{+}_{n})$. Hence,

\begin{equation}\label{eq5.38}
\int_0^{t^*(x)} \ \frac{ds}{V(\phi(s,x))} \quad \geq \quad
\int_{t^{-}_n}^{t^{+}_{n}} \ \frac{ds}{V(\phi(s,x))} \quad \geq
\quad  n+1 .
\end{equation}
Since this is true for all large $n$, $\int_0^{t^*(x)} \
\frac{ds}{V(\phi(s,x))}$ is infinite here too. Thus, (5.8) holds
for $V$ and so $V$ is $\phi$ regular.

$\Box$ \vspace{.5cm}

Our application of the Beck results is the following:

\begin{theo}\label{th5.8} If $\phi$ is a reversible semiflow on $X$, then it admits a
reversible Lyapunov function compactification  $\hat{\phi}$ on
$\hat{X}$ such that every point of $\hat{X} \setminus X$ is a
fixed point. That is,

\begin{equation}\label{eq5.39}
\hat{X} \setminus X \quad \subset \quad |\hat{\phi}|.
\end{equation}
If $X$ is metrizable  then the
compactification  can be chosen metrizable.
\end{theo}

{\bfseries Proof:} Use Theorem \ref{th5.7} to choose a proper $\phi$
regular function $V : X \to (0,\infty)$.  That is, for every
$\delta > 0$ there is a compact $X_{\delta} \subset X$ such that
on $X \setminus X_{\delta} \ V > 1/\delta$. Let $\psi : {\mathbb
R} \times X \to X$ be the flow associated with $\phi$ via the
time-change function $\bar{t}$. Since $X_0 = \emptyset \ \bar{t}
\times \pi_2$ is a homeomorphism on ${\mathbb R} \times X$. For
any $t \in {\mathbb R}_+$ let $\delta_t = \frac{\delta}{|t|+1}$
and $ X_{t,\delta} = (\phi^{[0,t]})^{-1}( X_{\delta_{t}})$ which
is a compact subset of $X$ because $\phi$ is proper.  If $x \in X
\setminus X_{t,\delta}$ then

\begin{equation}\label{eq5.40}
|\bar{\tau}(t,x))| \ = \ | \int_{0}^{t}  \ \frac{ds}{V(\phi(s,x))}
| \ < \  \delta.
\end{equation}

For each $x \in X$ the positive orbits of $\phi$ and $\psi$ are
the same set which says

\begin{equation}\label{eq5.41}
\mathcal{O}\phi \quad = \quad \mathcal{O}\psi \qquad \mbox{and so}
\qquad \mathcal{G}\phi \quad = \quad \mathcal{G}\psi.
\end{equation}
In particular, $L$ is a Lyapunov function for $\phi$ iff it is a
Lyapunov function for $\psi$.

Now let  $\mathcal{L}_0 \subset \mathcal{B}_{\psi}(X)$ be a
sufficient set of Lyapunov functions for $\psi$ . In particular, each $L$ is
$\psi$ uniform.  Hence, for every $L \in
\mathcal{L}_0$ and $\epsilon
> 0$ there exists a $\delta > 0$ such that $0 < \tau_1 - \tau_2  <
\delta$ imply $0 < L(\psi(\tau_1,x)) - L(\psi(\tau_2,x)) <
\epsilon $ for all $x \in X$. Since $\phi(t,x) =
\psi(\bar{\tau}(t,x),x)$ it follows from (\ref{eq5.40}) that

\begin{equation}\label{eq5.42}
x \in X \setminus X_{t,\delta} \qquad \Longrightarrow \qquad
|L(\phi(t,x)) - L(x)| \ \leq \ \epsilon.
\end{equation}
It easily follows that $\mathcal{L}_0 \subset \mathcal{B}_{\phi}(X)$ as well.
Hence, the compactification $\X$ obtained by using the $L \circ
\phi^t$'s is a Lyapunov compactification for the flow $\phi$ and  for every $z \in \hat{X} \setminus X$
and every $t \in {\mathbb R} \ \hat{L}(\hat{\phi}(t,z)) =
\hat{L}(z)$. Since these Lyapunov functions distinguish the points
at infinity, it follows that $\hat{\phi}(t,z) = z$ for all $(t,z)
\in {\mathbb R} \times (\hat{X} \setminus X)$ proving (\ref{eq5.39}).

As usual if $X$ is metrizable we can choose $\mathcal{L}_0$
countable and obtain a metrizable compactification.

$\Box$ \vspace{.5cm}

We prove the analogous result for a cascade by using the flow
result.  This requires the \emph{suspension} \index{suspension} construction which
builds a semiflow $\phi : {\mathbb R}_+ \times Y \to Y$ from a
cascade $f$ on $X$.  Begin with the trivial translation flow
$\phi_0 : {\mathbb R}_+ \times Y_0 \to Y_0$ on the product $Y_0 =
{\mathbb R}_+ \times X $ defining $\phi_0(t,(s,x)) = (t + s,x)$.
On $Y_0$ take the equivalence relation such that $(s + n,x) $ is
identified with $(s,f^n(x))$ for every positive integer $n$. The
quotient space is the same as the one obtained from $I \times X$
by identifying $(0,x)$ with $(1,f(x))$. On the quotient we obtain
the semiflow $\phi$. We regard embedding $X \to Y$ given by $x
\mapsto (0,x)$ as an identification so that $X$ is a $+$invariant
set for the time-one map $\phi^1$ with $f = \phi^1|X$. We use $f$
to stand for the time-one map on all of $Y$. We will write the
points of $Y$ as $[s,x]$ with $s \in I$

Notice that the projection $\pi_1 : Y_0 \to {\mathbb R}_+$ factors
through the equivalence relation to define the map $\pi : Y \to
{\mathbb R}/{\mathbb Z}$ which maps $\phi$ on $Y$ to the
translation flow on the circle. In particular, $\pi$ maps $f$ on
$Y$ to the identity on the circle and so we obtain:

\begin{equation}\label{eq5.43}
\mathcal{G}f \quad \subset \quad (\pi \times \pi)^{-1}(1_{{\mathbb
R}/{\mathbb Z}}).
\end{equation}
From this we obtain for $(t,x), (s,y) \in [0,1] \times X$

\begin{equation}\label{eq5.44}
\begin{split}
[s,y] \in [1_Y \cup \mathcal{G}f]([t,x]) \qquad
\Longleftrightarrow \hspace{2cm} \\ [s,y] \in
\mathcal{G}\phi([t,x]) \quad \mbox{and} \quad s=t \ \mbox{or} \
\{s,t\} = \{0,1\}
\end{split}
\end{equation}
because by (\ref{eq4.9}) $\mathcal{G}\phi = \phi^I \cup
\mathcal{G}(\phi^J)$ and by (\ref{eq4.10}) $\mathcal{G}(\phi^J) =
\mathcal{G}f \circ \phi^I$.

Notice that if $f$ is a homeomorphism on $X$ then the semiflow
$\phi$ on $Y$ is reversible. It is easy to check that if $f$ is a
proper continuous map then the suspension semiflow is proper.

\begin{theo}\label{th5.9}  If $f$ is a homeomorphism on $X$ there exists a
Lyapunov function compactification $(\hat{X},\hat{f})$ of the
cascade $(X,f)$ such that $\hat{f}$ is a homeomorphism on
$\hat{X}$ and every point of $\hat{X} \setminus X$ is a fixed
point of $\hat{f}$. That is,

\begin{equation}\label{eq5.45}
\hat{X} \setminus X \quad \subset \quad |\hat{f}|. \hspace{2cm}
\end{equation}

If $X$ is metrizable  then the
compactification $(\hat{X},\hat{f})$ can be chosen metrizable.
\end{theo}

{\bfseries Proof:} Apply the proof of Theorem \ref{th5.8} to the
suspension flow $\phi$ on $Y$ and its associate $\psi$. We obtain
a sufficient set $\mathcal{L}_0$ of $\psi$ uniform Lyapunov
functions  which we regard as $f$ Lyapunov functions by
restricting to $X$. The time change map may destroy the
factorization over the circle flow but this does not matter since
the $L$'s are $\mathcal{G}\phi$ Lyapunov functions.

By (\ref{eq5.44}) if $(x,y) \not\in \mathcal{G}f$ then $([0,x],[0,y])
\not\in \mathcal{G}\phi$ and so there exists an $L \in
\mathcal{L}_0$ such that $L(x) > L(y)$.  Hence, we can use the set
$\mathcal{L} = \{ L\circ f^m : L \in \mathcal{L}_0 \}$ to define a
Lyapunov function compactification for $f$. Since $f$ is the
time-one map for $\phi$ we have that for every $L \in \mathcal{L}$
and $ m \in \Z \ \  L(f^m(x)) - L(x)$ tends to zero as $x$ tends to infinity.
Hence, just as in Theorem \ref{th5.8} all of the points of $\hat{X}
\setminus X$ are fixed points.

$\Box$ \vspace{.5cm}

The suspension construction allows us to compare
 the Lyapunov function compactification of a semiflow
$\phi$ and the Lyapunov compactification of its time-one map
$f$.

\begin{ex}\label{ex5.10}The Lyapunov compactification of a semiflow and
that of its time-one map can be different.\end{ex}
 Begin with a homeomorphism $f$ on a
locally compact space $X$ which has a noncompact $\mathcal{G}f
\cap \mathcal{G}f^{-1}$ equivalence class $E_0$. For example we
can begin with a topologically transitive homeomorphism on a
compact space which has a fixed point. Removing the fixed point we
get a topologically transitive homeomorphism on a non compact
space $X$ and the entire space is a single $\mathcal{G}f \cap
\mathcal{G}f^{-1}$ equivalence class.  Let $\phi$ be the flow on
$Y$ which is the suspension of $f$ on $X$ extend $f$ to denote the
time-one map on $Y$. With $\pi$ the projection from $Y$ to the
circle $ {\mathbb R}/{\mathbb Z}$ (5.43) says that $\mathcal{G}f
\subset (\pi \times \pi)^{-1}(1_{{\mathbb R}/{\mathbb Z}})$.
Hence, for each $t \in [0,1)$ the image $E_t = \phi^t(E)$ is a
separate $\mathcal{G}f \cap \mathcal{G}f^{-1}$ equivalence class
and so in any Lyapunov function compactification of $f$ the
$E_t$'s all have pairwise disjoint closures by Theorem \ref{th1.6} (c). But
their union $E$ is a single $\mathcal{G}\phi \cap
\mathcal{G}\phi^{-1}$ equivalence class and so the closure $E$ in
any Lyapunov compactification of $\phi$ is the one point
compactification of $E$. The difficulty comes from the fact that
if $x \in E_t$ and $y \in E_s$ with $t \not= s$ then $x$ and $y$
can be distinguished by some $f$ Lyapunov function but not by any
$\phi$ Lyapunov function.

$\Box$ \vspace{.5cm}

\section{Parallelizable Systems}\label{secpara}

In this section we describe the results of Antosiewicz and
Dugundji (1961), see also Markus (1969)

Let $h :X_1 \to X_2$ be a continuous map.  We say that $h$ is an
\emph{action map} from a cascade $(X_1,f_1)$ to a cascade
$(X_2,f_2)$, or just that $h$ maps $f_1$ to $f_2$, if

\begin{equation}\label{eq6.1}
h \circ f_1 \quad = \quad f_2 \circ h,
\end{equation}
or, equivalently, if $h \circ f_{1}^n =  f_{2}^n \circ h$ for
every $n \in {\mathbb Z}_+$. If both $f_1$ and $f_2$ are
homeomorphisms then this equation holds for all $n \in {\mathbb
Z}$. If $h$ is an action map homeomorphism then $h^{-1}$ is an
action map from $(X_2,f_2)$ to $(X_1,f_1)$. In that case, we will
call $h$ an \emph{isomorphism} \index{isomorphism} between the cascades.

If $\phi_1$ and $\phi_2$ are semiflows on $X_1$ and $X_2$,
respectively, then we say that $h$ is an \emph{action map} from
$\phi_1$ to $\phi_2$, or just that $h$ maps $\phi_1$ to $\phi_2$,
if $h$ maps $\phi_{1}^t$ to $\phi_{2}^t$ for every $t \in {\mathbb
R}_+$, or, equivalently, if

\begin{equation}\label{eq6.2}
h(\phi_1(t,x)) \quad = \quad \phi_2(t,h(x))
\end{equation}
for every $(t,x) \in {\mathbb R}_+ \times X_1$. If $\phi_1$ and
$\phi_2$ are reversible and we extend to the associated flows then
these results hold for all $t \in {\mathbb R}$.  If $h$ is an
action map homeomorphism then $h^{-1}$ is an action map from
$\phi_2$ to $\phi_1$ and we call $h$ an isomorphism between the
semiflows or between the flows.

For a space $Y$ recall that the translation reversible semiflow,
and associated flow,  $\tau$ on ${\mathbb R} \times Y$ are defined
by

\begin{equation}\label{eq6.3}
\tau(t,(s,y)) \quad = \quad (t+s,y).
\end{equation}
with $t \in {\mathbb R}_+$ or ${\mathbb R}$. Clearly,

\begin{equation}\label{eq6.4}
\begin{split}
\mathcal{O} \tau \quad = \quad \{ ((s_1,y_1),(s_2,y_2)) : s_1 \ \leq
\ s_2 \ \mbox{and} \ y_1 \ = \ y_2 \}. \hspace{1cm}\\
\mathcal{O} (\tau^J) \quad = \quad \{ ((s_1,y_1),(s_2,y_2)) : s_1 +
1 \  \leq \ s_2 \ \mbox{and} \ y_1 \ = \ y_2 \}.
\end{split}
\end{equation}
Since these relations are closed as well as transitive, we have

\begin{equation}\label{eq6.5}
\begin{split}
\mathcal{O} \tau \  = \ \mathcal{G} \tau \quad \mbox{and}\quad
\mathcal{O}( \tau^J) \ = \ \mathcal{G} (\tau^J) \\
|\mathcal{G} (\tau^J)| \quad = \quad \emptyset. \hspace{2cm}
\end{split}
\end{equation}

Recall from (\ref{eq4.6})

\begin{equation}\label{eq6.6}
\begin{split}
\mathcal{O} \phi  \quad =_{def}  \quad \mathcal{O} (\phi^I) \quad
= \quad
\bigcup \{ \phi^t : t \geq 0 \}.  \\
\mathcal{R}\phi   \quad =_{def}  \quad \mathcal{R} (\phi^I). \hspace{3cm} \\
\mathcal{N} \phi  \quad =_{def}  \quad \mathcal{N} (\phi^I) \quad
= \quad \overline{\mathcal{O} \phi} \quad = \quad \mathcal{O} \phi \ \cup \ \Omega \phi, \\
\mbox{where} \quad \Omega \phi \quad =_{def} \quad Limsup_{t \to
\infty} \{ \phi^t \}. \\
\mathcal{N} (\phi^J) \quad = \quad \overline{\mathcal{O} (\phi^J)}
\quad = \quad \mathcal{O} (\phi^J) \ \cup \ \Omega \phi.
\end{split}
\end{equation}
It easily follows that

\begin{equation}\label{eq6.7}
|\Omega \phi | \quad = \quad | \mathcal{O}(\phi^J)|. \hspace{2cm}
\end{equation}
We call a point \emph{non-wandering} for $\phi$ when it lies in
this set.

For the translation semiflow we have

\begin{equation}\label{eq6.8}
\Omega \tau \quad = \quad \emptyset. \hspace{2cm}
\end{equation}

For a reversible semiflow $\phi$ on $X$, we will call $\phi$,or
the associated flow, \emph{parallelizable} \index{parallelizable flow} if there is a
homeomorphism from $X$ to ${\mathbb R} \times Y$ for some space
$Y$ which is an isomorphism from $\phi$ to the translation $\tau$.
The space  $Y$ is then called a \emph{section} for $\phi$.

If $A$ is a closed $\phi \ +$invariant subset of $X$ then we let
$\phi_A$ denote the semiflow on $A$ obtained by restricting
$\phi$.  If $A$ is $\phi$ invariant and $\phi$ is reversible then
$\phi_A$ is reversible. All this is true if $A$ is open as well
except that here a slight quibble arises. An open subset of a
space $X$ is locally compact but it is also $\sigma$-compact iff
it is an $F_{\sigma}$ set, or, equivalently, the complement of the
zero-set of a real-valued continuous function on $X$. For a
general space $X$ this is not true of every open set but the open
$F_{\sigma}$ sets do form a basis. For convenience, we will call
$U$ a $\phi$ \emph{open} set when it is a $\phi$ invariant open
$F_{\sigma}$ set.  Then we can define the restriction $\phi_U$
which is reversible when $\phi$ is.

\begin{theo}\label{th6.1} Let $\phi$ be a reversible semiflow on $X$ with time-one
homeomorphism $f$ and let $x \in X$. The following conditions are
equivalent. \begin{enumerate}
\item[(i)] The point $x$ is a wandering point, i.e. $x \not\in
|\mathcal{N}(\phi^J)|$.
\item[(ii)] There exists a $\phi$ open set $U$ with $x \in U$ such that is not
a generalized recurrent point for the restriction  $\phi_U$, i.e.
$x \not\in |\mathcal{G}((\phi_U)^J)|$.
\item[(iii)] There exists a $\phi$ open set $U$ with $x \in U$ and
$L : U \to [-1,1]$ a Lyapunov function for $\phi_U$ such that

\begin{equation}\label{eq6.9}
 L(x) \ = \ 0 \quad \mbox{and} \quad L(f(x)) \ = \ 1.
\end{equation}
\item[(iv)] There exists a $\phi$ open set $U$ with $x \in U$ such that
the restriction $\phi_U$ is parallelizable.
\end{enumerate}
\end{theo}

{\bfseries Proof:} (ii) $\Rightarrow$ (iii): By Corollary \ref{cor4.5} (a) $x
\not\in |\mathcal{G}((\phi_U)^J)|$ implies $(f(x),x) \not\in
\mathcal{G} (\phi_U)$. By Theorem \ref{th4.13} there  is a $\phi_U$
Lyapunov function $L : U \to I$ such that $L(x)  = 0$ and $L(f(x))  = 1$.

(iii) $\Rightarrow$ (iv):  For $z \in U$ define $K(z) = \int_{0}^1
\ L(\phi(u,z)) \ du$. As in the proof of Theorem \ref{th4.13} this is a
Lyapunov function for $\phi_U$ and

\begin{equation}\label{eq6.10}
\dot{K}(z) \quad =_{def}\quad \frac{d}{ds} K(\phi(s,z)) |_{s=0}
\quad = \quad L(f(z)) - L(z)
\end{equation}
is nonnegative on $U$ and positive at $x$. Choose $V_0$ an open
subset of $U$ with $x \in V_0$ such that

\begin{equation}\label{eq6.11}
z \in V_0 \quad \Longrightarrow \quad \dot{K}(z) > 0.
\end{equation}
There exist positive $\epsilon$ and $\delta$ such that
$\phi([-\epsilon,\epsilon] \times \{ x \} \subset V_0,$
$K(\phi(-\epsilon,x)) < - \delta $ and $K(\phi(\epsilon,x)) >
\delta$. Choose  $W_0$ an open $F_{\sigma}$ with $x \in W_0$ and
with compact closure $\overline{W_0} \subset V_0$ such that

\begin{equation}\label{eq6.12}
\begin{split}
\phi([-\epsilon,\epsilon] \times \overline{W_0} ) \quad \subset
\quad V_0, \hspace{3cm}\\
z \in \overline{W_0} \quad \Longrightarrow \quad
K(\phi(-\epsilon,z)) \ < \ - \delta, \quad K(\phi(\epsilon,z)) \ >
\ \delta.
\end{split}
\end{equation}

%6.9
%\begin{equation}
%V_0 \quad = \quad L^{-1}((-\frac{1}{2},\frac{1}{2})) \cap f(
%L^{-1}([-1,\frac{1}{2}))) \cap f^{-1}(L^{-1}((\frac{1}{2},1]))
%\end{equation}
%$V_0$   an open $F_{\sigma}$  subset  with $x \in V_0$, compact
%closure $ \overline{V_0}$ contained in $ U$ and such that $y \in
%\overline{V_0}$ implies
%%6.9
%\begin{equation}
%L(f^{-1}(y)) \quad < \quad - \ \frac{1}{2} \quad < \quad L(y)
%\quad < \quad \frac{1}{2} \quad < \quad L(f(y))).
%\end{equation}
Now let

\begin{equation}\label{eq6.13}
\begin{split}
W \quad =_{def} \quad \bigcup_{t \in {\mathbb R}} \phi^t(W_0)
\quad = \quad \bigcup_{k \in {\mathbb Z}_+} \phi^{[-k,k]}(W_0) \\
\tilde{W} \quad =_{def} \bigcup_{k \in {\mathbb Z}_+}
\phi^{[-k,k]}(\overline{W_0}). \hspace{3cm}
\end{split}
\end{equation}
$W$ is a $\phi$ open set containing $x$, $\tilde{W}$ is a $\phi$
invariant set containing $W$ and each
$\phi^{[-k,k]}(\overline{W_0})$ is a compact subset of
$\tilde{W}$.
%Furthermore, (6.9) implies
%%6.11
%\begin{equation}
%\tilde{V} \ \cap \ |\mathcal{G}((\phi_U)^J)| \quad = \quad
%\emptyset,
%\end{equation}
% i.e. no point of $\tilde{V}$ is generalized recurrent for the
%restriction $\phi_U$.

By the Intermediate Value Theorem there exists for every $z \in
\overline{W_0}$ a time $t \in (-\epsilon,\epsilon)$ such that
$K(\phi(t,z)) = 0$. By (\ref{eq6.11}) the Lyapunov function $K$ is
strictly increasing in $V_0$. Hence, for each $z \in \tilde{W}$
the time $t \in {\mathbb R}$ such that $K(\phi(t,z)) = 0$ is
unique.

So if we define $Y = K^{-1}(0) \cap W$ its closure $\overline{Y}$
is contained in $\phi^{[-\epsilon,\epsilon]}(\overline{W_0})$ and
so is compact. From (\ref{eq6.12}) it then follows that for all $y \in
\overline{Y}$

\begin{equation}\label{eq6.14}
\begin{split}
\delta \quad \leq \quad K(\phi(s,y)) \qquad \mbox{for all} \
s \geq 2\epsilon, \\
- \delta \quad \geq \quad K(\phi(s,y)) \qquad \mbox{for all} \ s
\leq -2\epsilon.
\end{split}
\end{equation}

 If the map $\tilde{h} :
{\mathbb R} \times \overline{Y} \to U$ is given by

\begin{equation}\label{eq6.15}
\tilde{h}(s,y) \quad =_{def} \quad \phi(s,y)
\end{equation}
then $\tilde{h}$ is injective and continuous. The restriction, $h:
{\mathbb R} \times Y \to W$ is a continuous bijection which maps
the translation flow $\tau$ on ${\mathbb R} \times Y$ to the
restriction $\phi_W$. So $h$ is a homeomorphism on compacta. To
show that $h$ is a homeomorphism we need only prove it is proper.

Now suppose instead that $\{ (s_i,y_i) \}$ is an unbounded net in
${\mathbb R} \times Y$  with $ \{ \phi(s_i,y_i) \}$ tending to $z
\in W$.  Since $\tilde{h}$ is a continuous injection on each
compact set $\phi^{[-k,k]}(\overline{W_0})$ it cannot happen that
$\{ s_i \}$ remains bounded in ${\mathbb R}$. We can assume that
$\{ s_i \}$ tends to $+ \infty$. The limit point $z$ is in
$\phi^t(Y)$ for some $t$ and so by replacing $s_i$ by $s_i - t$ we
can assume that $z \in Y$ and so $K(z) = 0$. However, each $y_i
\in Y$ and so $K(\phi(s_i,y_i)) \geq \delta$ once $s_i \geq 2 \epsilon$
by (\ref{eq6.13}). This contradiction completes the proof of (iv).

(iv) $\Rightarrow$ (i):  If $x$ is a nonwandering point for $\phi$
then it is a nonwandering point for any restriction $\phi_U$ when
$U$ is a $\phi$ open set containing $x$. If $\phi_U$ is
parallelizable then every point is wandering by (\ref{eq6.5}).

(i) $\Rightarrow$ (ii): Since $x \not\in \Omega(\phi)(x)$ there
exists an open $F_{\sigma}$ subset $U_0$ with $x \in U_0$ and
compact closure $\overline{U_0}$ such that $\overline{U_0}$ is
disjoint from the closed invariant set
$\Omega(\phi)(\overline{U_0})$. As above define $U = \bigcup_t \{
\phi^t(U_0) \}$ to get a $\phi$ open set which contains $x$ and is
disjoint from  $\Omega(\phi)(\overline{U_0})$. The closure of
$\mathcal{O}(\phi_U)$ in $U \times U$ is contained in the
intersection of $\mathcal{N}\phi = \mathcal{O}\phi \cup
\Omega(\phi)$ with $U \times U$.  But $\Omega(\phi)(U) \cap U =
\emptyset$. Hence, the intersection is $\mathcal{O}(\phi_U)$.

This means that $\mathcal{O}(\phi_U)$ is closed as well as
transitive in $U \times U$. Hence, $\mathcal{O}(\phi_U) =
\mathcal{G}(\phi_U)$ and since $(f(x),x) \not\in
\mathcal{O}(\phi_U)$ (e.g. because $x$ is wandering) it follows
from Corollary \ref{cor4.5}(a) that $x$ is not generalized recurrent for
$\phi_U$.

$\Box$ \vspace{.5cm}

\begin{prop}\label{prop6.2} Let $E$ be a closed equivalence relation on a space
$X$. Let $X/E$ denote the set of equivalence classes with $\pi : X
\to X/E$ the canonical projection. With the quotient topology
$X/E$ is a $\sigma$-compact Hausdorff space (not necessarily
locally compact).
\end{prop}

{\bfseries Proof:} If $A, B$ are disjoint closed subsets of $X$
with $E(A) = A$ and $E(B) = B$ then by Corollary \ref{cor1.10} there is a
Lyapunov function $L : X \to I$ for $E$ which is zero on $A$ and
one on $B$.  As $L$ is constant on every $E$ equivalence class $L$
factors through $\pi$ to define a map $L : X/E \to I$ which is
continuous by definition of the quotient topology.  Thus, the
real-valued continuous functions on $X/E$ distinguish closed sets
and so $X/E$ is a normal space.  Since individual classes are
closed (i.e. each $E(x)$ is a closed subset of $X$) it follows
that $X/E$ is Hausdorff. The continuous image of a
$\sigma$-compact space is $\sigma$-compact and so $X/E$ is
$\sigma$-compact. However, the continuous image of a locally
compact space need not be locally compact, e.g. move an exterior
point to the boundary of an open square in ${\mathbb R}^2$.
With $X = (0,1) \times [0,1)$ and  $E = 1_X \cup ((0,1)\times \{ 0 \})^2$
the quotient space, obtained by smashing $(0,1) \times \{ 0 \}$ to a point,
 is not locally compact. To see this, map $X$ to $Y = \{(0,0) \} \cup \{ (x,y) \in (0,1) \times (0,1) : y < x \}$
 by $(x,y) \mapsto (x,xy)$. The map factors through $E$ to obtain a homeomorphism of
 $X/E$ onto $Y$.  

$\Box$ \vspace{.5cm}

\begin{theo}\label{th6.3} Let $\phi$ be a reversible semiflow on $X$.
\begin{enumerate}
\item[(a)] The following conditions are equivalent:
\begin{itemize}
\item[(i)] The reflexive, transitive relation $\mathcal{O}\phi$ is
closed.
\item[(ii)] The transitive relation $\mathcal{O}(\phi^J)$ is closed.
\item[(iii)] $\mathcal{O}\phi = \mathcal{N}\phi$.
\item[(iv)] $\mathcal{O}\phi = \mathcal{G}\phi$.
\item[(v)] $\mathcal{O}(\phi^J) = \mathcal{G}(\phi^J)$.
\end{itemize}
The above conditions imply that the equivalence relation
$\mathcal{O}(\phi \cup \phi^{-1})$ is closed, or, equivalently,
 $\mathcal{O}(\phi \cup \phi^{-1}) = \mathcal{G}(\phi \cup
\phi^{-1})$.

\item[(b)]  The following are equivalent:
\begin{itemize}
\item[(i)] $\mathcal{O}\phi$ is closed and there are no periodic points, i.e. $|\mathcal{O}(\phi^J)|
= \emptyset$.
\item[(ii)] $\mathcal{O}\phi$ is closed and all points are wandering, i.e. $|\mathcal{N}(\phi^J)|
= \emptyset$.
\item[(iii)] $\mathcal{O}\phi$ is closed and there are no generalized recurrent points, i.e. $|\mathcal{G}(\phi^J)|
= \emptyset$.
\item[(iv)] $\Omega \phi = \emptyset$.
\item[(v)] $\mathcal{O}(\phi \cup \phi^{-1})$ is closed and all points are wandering, i.e. $|\mathcal{N}(\phi^J)|
= \emptyset$.
\item[(vi)] $\phi$ is parallelizable.
\end{itemize}
\end{enumerate}
\end{theo}

{\bfseries Proof:} (a): (ii) $\Rightarrow$ (i): $\mathcal{O}\phi =
\phi^I \cup \mathcal{O}(\phi^J)$ and $\phi^I$ is closed.

(i) $\Rightarrow$ (ii): $\mathcal{O}(\phi^J) = \mathcal{O}\phi
\circ \phi^J$. Apply Lemma \ref{lem4.1}.

(i) $\Leftrightarrow$ (iii): $\mathcal{N}\phi$ is the closure of
$\mathcal{O}\phi$.

(i) $\Leftrightarrow$ (iv), (ii) $\Leftrightarrow$ (v): For any
closed relation $f$, the transitive relation $\mathcal{O}f$ is
closed iff $\mathcal{O}f = \mathcal{G}f$. This also shows that
$\mathcal{O}(\phi \cup \phi^{-1})$ is closed iff
 $\mathcal{O}(\phi \cup \phi^{-1}) = \mathcal{G}(\phi \cup
\phi^{-1})$.

Finally, since $\phi$ is reversible, $\mathcal{O}(\phi \cup
\phi^{-1}) = (\mathcal{O}\phi) \cup (\mathcal{O}\phi)^{-1}$.
Hence, $\mathcal{O}(\phi \cup \phi^{-1})$ is closed if
$\mathcal{O}\phi $ is closed.

(b): (i) $\Leftrightarrow$ (ii) $\Leftrightarrow$ (iii): By part
(a) $\mathcal{O}\phi$ is closed iff $\mathcal{O}\phi =
\mathcal{N}\phi$ and iff $\mathcal{O}\phi = \mathcal{G}\phi$.
Furthermore, these conditions each imply $\mathcal{O}(\phi^J) =
\mathcal{N}(\phi^J) = \mathcal{G}(\phi^J)$ and so the latter all
have the same cyclic set.

(vi) $\Rightarrow$ (iv): Follows from (\ref{eq6.8}).

(iv) $\Rightarrow$ (ii):  $\mathcal{N}\phi = \mathcal{O}\phi \cup
\Omega \phi$ and if $x$ is nonwandering then $x \in |\Omega
\phi|$.

(i) $\Rightarrow$ (v): Again this follows from part (a).

(v) $\Rightarrow$ (vi): Let $E = \mathcal{O}(\phi \cup
\phi^{-1})$. By (v) this is a closed equivalence relation and so
by Proposition \ref{prop6.2} the quotient space $X/E$ is Hausdorff. Let $\pi
: X \to X/E$ be the quotient map. For any set $A \subset X$

\begin{equation}\label{eq6.16}
\pi^{-1}(\pi(A)) \quad = \quad \bigcup_{t \in {\mathbb R}} \{
\phi^t(A) \}.
\end{equation}
With $A$ open we see that $\pi$ is an open map.  Hence, $X/E$ is
locally compact as well as $\sigma$-compact and Hausdorff.

Any $x \in X$ is wandering and so by Theorem \ref{th6.1} there exists a
$\phi$ open set $U$ with $x \in U$ and an isomorphism $h_U$ from a
translation flow $\tau$ on $ {\mathbb R} \times Y$ to the
restriction $\phi_U$.  Since $U$ is $\phi$ invariant it equals
$\pi^{-1}(\pi(U))$ and $h_U$ induces a homeomorphism between
$\pi(U)$ and $Y$ identified with the space of $\mathcal{O}(\tau)$
equivalence classes. Using this homeomorphism we can assume that
$Y = \pi(U)$ and that $h_U$ is an isomorphism from $\tau_U$ on
${\mathbb R} \times \pi(U)$ to $\phi_U$ inducing the identity on
$\pi(U)$. If $h_V$ is similarly an isomorphism from $\tau_V$ on
${\mathbb R} \times \pi(V)$ to $\phi_V$ then we obtain an
automorphism $h_{U,V} = (h_V)^{-1} \circ h_U$ of $\tau_{U \cap V}$
on ${\mathbb R} \times \pi(U \cap V)$. Furthermore, this map is of
the form

\begin{equation}\label{eq6.17}
h_{U,V}(s,a) \quad = \quad (H_{U,V}(a) + s,a)
\end{equation}
where the \emph{transition map} $H_{U,V} : \pi(U \cap V) \to
{\mathbb R} $ is the projection to the ${\mathbb R}$ coordinate of
$h_{U,V}(0,a)$. We will call these open sets $\pi(U)$ of $\pi(X)$ the
\emph{trivializing open sets}.

Thus, the isomorphisms $h_U$ give the map $\pi : X \to X/E$ the
structure of a principal ${\mathbb R}$ bundle. The result (vi)
follows from the fact that a bundle with a contractible fiber has
a section and a principal bundle which admits a global section is
trivial.

In detail, let $\{A_1,A_2,... \}$ be a sequence of compact subsets
whose interiors cover $\pi(X)$ and such that $A_i$ is contained in
the trivializing open set $\pi(U_i)$. Let $B_n = A_1 \cup ... \cup
A_n$ for $n = 1,...$

Assume inductively that $h_{B_n}$ is an isomorphism from
$\tau_{B_{n}}$ on ${\mathbb R} \times B_n$ to the restriction of
$\phi$ to $\pi^{-1}(B_n)$. Using the trivializing neighborhood
$U_{n+1}$ we have an isomorphism $h_{A_{n+1}}$ from
$\tau_{A_{n+1}}$ to $\pi^{-1}(A_{n+1})$. Use the Tietze Extension
Theorem to extend the transition map $H_{B_n,A_{n+1}} : B_n \cap
A_{n+1} \to {\mathbb R}$ to a continuous map $H_{n+1} : A_{n+1}
\to {\mathbb R}$. Now define $h_{B_{n+1}}$ by

\begin{equation}\label{eq6.18}
h_{B_{n+1}}(s,a) \quad = \quad \begin{cases} h_{B_n}(s,a) \hspace{3cm}\mbox{for} \ a \in B_n \\
h_{A_{n+1}}(H_{n+1}(a) + s,a) \qquad \mbox{for} \ a \in A_{n+1}.
\end{cases}
\end{equation}
This extends $h_{B_n}$.  Since every point of $X/E = \pi(X)$ is
eventually in the interior of some $B_n$, these isomorphisms fit
together to define a parallelism of $\phi$.

$\Box$ \vspace{.5cm}

\begin{ex}\label{ex6.4}The requirement in Theorem \ref{th6.3} (b) that $\mathcal{O}\phi$ be closed
is necessary. \end{ex}
 We recall the  example of Akin (1993)
Problem 4.22. On $ \tilde{X} = [0,1] \times [-1,1]$ identify
$(0,y)$ with $(1,-y)$ to obtain a M\"{o}bius strip \index{M\"{o}bius strip} $X_0$. Let $g$ be a
smooth nonnegative function on $\tilde{Y}$ with $g(0,y) = g(1,-y)$
and with $g^{-1}(0) = 0 \times [-1,0] \cup 1 \times [0,1]$. Let
$\phi_0$ be the reversible semiflow which is the solution of the
differential equations $\frac{dx}{dt} = g(x,y), \frac{dy}{dt} =
0$. By removing the set $g^{-1}(0)$ of fixed points from $X_0$ we
obtain a $\phi_0$ open set $X$. Let $\phi$ be the restriction of
$\phi_0$ to $X$. It is easy to check that
%6.19
\begin{equation}\label{eq6.19}
\mathcal{O}(\phi \cup \phi^{-1}) \quad = \quad \{
((x_1,y_1),(x_2,y_2)) \in X \times X : y_1 = \pm y_2 \}
\end{equation}
and so is a closed relation. Clearly there are no periodic points,
i.e. $|\mathcal{O}(\phi^J)| = \emptyset$. On the other hand,
 the points on the central
circle $(0,1) \times 0$ are nonwandering from which it follows
that $\mathcal{O}\phi$ is not closed.

$\Box$ \vspace{.5cm}

%Let $X$ denote the unit disc in ${\mathbb R}^2$ with the origin
%and $(1,0)$ removed. Consider the flow which is the solution of
%the system -in polar coordinates- $\frac{dr}{dt} = 0, \frac{d
%\theta}{dt} = (1-r) + sin^2 (\frac{ \theta}{2})$. It is clear that
%$\mathcal{O}(\phi \cup \phi^{-1})$ is closed and easy to check
%that $\mathcal{O}\phi$ is not.

\begin{theo}\label{th6.5}  Let $f$ be a homeomorphism on $X$. For the cascade
$(X,f)$ the following conditions are equivalent:
\begin{itemize}
\item[(i)] The transitive relation $\mathcal{O}f$ is closed and
there are no periodic points, i.e. $|\mathcal{O}f| = \emptyset$.
\item[(ii)]$\mathcal{O}f = \mathcal{N}f $ and
all points are wandering, i.e. $|\mathcal{N}f| = \emptyset$.
\item[(iii)] $\mathcal{O}f = \mathcal{G}f $ is closed and
there are no generalized recurrent points, i.e. $|\mathcal{N}f| =
\emptyset$.
\item[(iv)] $\Omega f = \emptyset$.
\item[(v)] The equivalence relation $\mathcal{O}(f \cup f^{-1})$
is closed and all points are wandering,i.e. $|\mathcal{N}f| =
\emptyset $.
\item[(vi)] There exists a translation flow $\tau$ on ${\mathbb R}
\times Y$ for some space $Y$  and a closed $\tau^1$ invariant
subset $Z$ of ${\mathbb R} \times Y$ such that the cascade $(X,f)$
is isomorphic to $(Z,g)$ where $g$ is the restriction to $Z$ of
the time-one map $\tau^1$.
\end{itemize}
\end{theo}

{\bfseries Proof:}  The proofs that (i) $\Leftrightarrow$ (ii)
$\Leftrightarrow$ (iii) and that  (vi) $\Rightarrow$ (iv)
$\Rightarrow$ (ii) $\Rightarrow$ (v) proceed just as in Theorem
\ref{th6.3}.

(v) $\Rightarrow$ (vi): Use the suspension construction as in the
proof of Theorem \ref{th5.9}.  We obtain a reversible flow $\phi$ on a
space $\tilde{X}$ and a closed $\phi^1$ invariant subset
$\tilde{Z}$ such that $(X,f)$ is isomorphic to
$(\tilde{Z},\phi^1)$. The product structure on $\tilde{X}$ shows
that $\mathcal{O}(\phi \cup \phi^{-1})$ is an extension of
$\mathcal{O}(f \cup f^{-1})$ which is closed when the latter is.
It is easy to check that every point of $\tilde{X}$ is wandering
for $\phi$ because every point of $X$ is wandering for $f$.

It follows from Theorem \ref{th6.3} that $\phi$ is parallelizable,
isomorphic via $h$ to a translation flow $\tau$. The required
subset $Z$ is the image of $\tilde{Z}$ under $h$.

$\Box$ \vspace{.5cm}

\section{Appendix: Limit Prolongation Relations}\label{secprolong}

Let $f$ be a closed relation on a space $X$.

The \emph{omega limit point set} of the orbit of $x \in X$ is
$\omega f(x) = lim sup \{ f^i(x) \}$ which defines the relation
$\omega f$. Recall that for a sequence $\{ A_i \} \ lim sup \{ A_i
\} = \bigcap_i \overline{\bigcup_{j \geq i } \{ A_j \}}$ so  that
$\overline{ \bigcup_i \{A_i\}} = \bigcup_i \{A_i \} \cup lim sup
\{ A_i \}$ when all of the sets $A_i$ are closed.

If $f$ is a + proper relation, e.g. a continuous map, then each
iterate $f^n$ is closed. So when $f$ is + proper, $\mathcal{R}f(x)
= \mathcal{O}f(x) \cup \omega f(x)$ is the closure of the orbit
$\mathcal{O}f(x)$ for each point $x \in X$. but $\mathcal{R}f$ may
be a proper subset of $\mathcal{N}f$, the closure of
$\mathcal{O}f$ in $X \times X$. If we define $\Omega f = lim sup
\{ f^i \}$ then $\mathcal{N} f = \mathcal{O}f \cup \Omega f $ when
that $f$ is + proper.

In Example \ref{ex1.1new} the discontinuous map $g$ is a closed relation such that
$g^2$ is not closed.

\begin{prop}\label{prop9.1}  Let $F$ be a  relation on $X$ and $f : X \to
X$ be a continuous map.
\begin{enumerate}
\item[(a)]  If $f \circ F \subset F$ then $f^n \circ F \subset F$
for $n = 1, 2, ...$ and

\begin{equation}\label{eq9.1}
\begin{split}
f \circ F^n \ \subset \ F^n  \quad \mbox{for} \ n = 1, 2, ...
\hspace{1cm} \\
f \circ \Omega F \ \subset \ \Omega F, \qquad f \circ \mathcal{N}F
\ \subset \ \mathcal{N}F,\\
f \circ \mathcal{G} F \ \subset \ \mathcal{G} F.\hspace{2cm}
\end{split}
\end{equation}

\item[(b)] If $f \circ F = F$ and $f$ is a proper map then $f^n \circ F = F$
for $n = 1, 2, ...$ and

\begin{equation}\label{eq9.2}
\begin{split}
f \circ F^n \ = \ F^n  \quad \mbox{for} \ n = 1, 2, ...
\hspace{1cm} \\
f \circ \Omega F \ = \ \Omega F, \qquad f \circ \mathcal{N}F
\ = \ \mathcal{N}F,\\
f \circ \mathcal{G} F \ = \ \mathcal{G} F.\hspace{2cm}
\end{split}
\end{equation}

\item[(c)]  If $ F \circ f^{-1} \subset F$ then $F \circ f^{-n} \subset F$
for $n = 1, 2, ...$ and

\begin{equation}\label{eq9.3}
\begin{split}
  F^n \circ f^{-1} \ \subset \ F^n  \quad \mbox{for} \ n = 1, 2, ...
\hspace{1cm} \\
\Omega F \circ f^{-1} \ \subset \ \Omega F, \qquad \mathcal{N}F
\circ f^{-1}
\ \subset \ \mathcal{N}F,\\
\mathcal{G} F \circ f^{-1} \ \subset \ \mathcal{G} F.\hspace{2cm}
\end{split}
\end{equation}

\item[(d)]  If $F \subset F \circ f$, $F$ is closed and $X$ is compact then

\begin{equation}\label{eq9.4}
F \ \subset \ F \circ f \ \subset \ F \circ f^2 \ \subset \ F
\circ f^3 ... \ \subset \ F \circ \omega f \ \subset \ F \circ
\Omega f.
\end{equation}
\end{enumerate}
\end{prop}

{\bfseries Proof:}  (a), (c):  Observe that $f \circ F \subset F$
iff the subset $F$ is + invariant for the map $1_X \times f$ on $X
\times X$ and $F \circ f^{-1} \subset F$ iff $F$ is + invariant
for the map $f \times 1_X$.  With equality in each case equivalent
to invariance. The class of + invariant sets for a map is closed
under taking unions, intersections and closures.

Clearly, $f \circ F \subset F$ implies that the sequence $\{ f^n
\circ F \}$ is decreasing and that $f \circ F \circ G \subset F
\circ G$ for any relation $G$, it follows that $\mathcal{N} F$ and
$\Omega F$ are $1_X \times f$ + invariant.  By transfinite
induction $\mathcal{N}_{\alpha} F$ is + invariant for every
ordinal $\alpha$ and so $\mathcal{G} F$ is + invariant as well.
This proves (a) and the proof for (c) is similar.

(b):  For a proper (and hence closed) map the class of invariant
sets is closed under unions and closure. For any filter
$\mathcal{A}$ of closed invariant subsets, the intersection is
closed and invariant when the map is proper.  Now proceed as in
(a) above.

(d): That $\{ F \circ f^n \}$ is an increasing sequence is clear.
If $(x,y) \in F \circ f^n$, i.e. $(f^n(x),y) \in F$ then it
follows that for $i = 1,2,...$ $(f^{n+i}(x),y) \in F$.  If $z \in
\omega f(x)$ then $(z,y)$  is a limit point of this sequence and
so $(z,y) \in F$ because $F$ is closed.  $\omega f(x)$ is nonempty
by compactness. The final inclusion is obvious.

$\Box$ \vspace{.5cm}

\begin{theo}\label{th9.2}  Let $f$ be a continuous map on a space $X$. If
either $f$ is a homeomorphism or $X$ is compact then

\begin{equation}\label{eq9.5}
\mathcal{G} \Omega f \quad = \quad \Omega \mathcal{G} f \quad =
\quad \bigcap_{n=1}^{\infty} \ (\mathcal{G} f)^n. \hspace{1cm}
\end{equation}
\end{theo}

{\bfseries Proof:} For any continuous map $f$ the iterates are
closed relations and so  $\mathcal{N}f$, the closure of the orbit
relation $\mathcal{O}f$, is equal to $\mathcal{O}f \cup \Omega f$.
Since $\mathcal{G} \Omega f$ contains $\Omega f$ it follows that

\begin{equation}\label{eq9.6}
\mathcal{O}f \cup \mathcal{G}\Omega f \quad = \quad \mathcal{N}f
\cup \mathcal{G}\Omega f \hspace{1cm}
\end{equation}
is a closed relation.

Since $\Omega f \subset \mathcal{N} f \subset \mathcal{G} f$ it
follows that

\begin{equation}\label{eq9.7}
\mathcal{G} \Omega f \ \subset \ \mathcal{G} \mathcal{G} f \ = \
\mathcal{G} f. \hspace{1cm}
\end{equation}

Next define for $n = 1,2,...$

\begin{equation}\label{eq9.8}
\mathcal{O}_n f \quad =_{def} \quad \bigcup_{i=n}^{\infty} \
f^i.\hspace{2cm}
\end{equation}

  Since $f \circ \mathcal{O}_n f =
\mathcal{O}_{n+1} f$ it follows by taking closures and
intersections that $f \circ \Omega f = \Omega f$.  This only
requires that that $f$ is a proper map and so holds in both cases.
Hence, from Proposition \ref{prop9.1} (b) it follows that $f \circ
\mathcal{G}\Omega f  = \mathcal{G} \Omega f$ and so

\begin{equation}\label{eq9.9}
\mathcal{O} f \circ \mathcal{G} \Omega f \quad = \quad \mathcal{G}
\Omega f. \hspace{2cm}
\end{equation}

 At this point the proofs for the two cases diverge.
\vspace{.25cm}

{\bfseries Homeomorphism Case:} Apply (\ref{eq9.9}) to the inverse
homeomorphism $f^{-1}$ and then invert the equation to get

\begin{equation}\label{eq9.10}
\mathcal{G} \Omega f  \circ \mathcal{O} f \quad = \quad
\mathcal{G} \Omega f. \hspace{2cm}
\end{equation}

From (\ref{eq9.9}) and (\ref{eq9.10}) it follows that

\begin{equation}\label{eq9.11}
(\mathcal{O} f \cup \mathcal{G} \Omega f)^n \quad = \quad
\mathcal{O}_n f \cup \mathcal{G} \Omega f. \hspace{1cm}
\end{equation}
In particular, the closed relation $\mathcal{O} f \cup \mathcal{G}
\Omega f$ is transitive and so contains $\mathcal{G} f$. From
(\ref{eq9.7}) we obtain

\begin{equation}\label{eq9.12}
\mathcal{O} f \cup \mathcal{G} \Omega f \quad = \quad \mathcal{G}
f. \hspace{2cm}
\end{equation}

Now (\ref{eq9.11}) says that

\begin{equation}\label{eq9.13}
(\mathcal{G}f)^n \quad = \quad \mathcal{O}_n f \cup \mathcal{G}
\Omega f. \hspace{2cm}.
\end{equation}

A priori the relations $(\mathcal{G}f)^n $ need not be closed, but
as in (\ref{eq9.6}), (\ref{eq9.13}) implies that these iterates are closed and so
their intersection is $ Lim sup \{ (\mathcal{G}f)^n  \} = \Omega
\mathcal{G} f$. Also (\ref{eq9.13}) implies that the intersection is
$\mathcal{G} \Omega f$ because the intersection of the
$\mathcal{O}_n f$'s is contained in $\Omega f$. \vspace{.25cm}

{\bfseries Compact Space Case:}  It is an easy exercise to prove
that $\Omega f \subset \Omega f \circ f$ (see Prop. 1.12(a) of Akin (1993)).  That is, if $(x,y) \in
\Omega f$ then $(f(x),y ) \in \Omega f$. Hence,

\begin{equation}\label{eq9.14}
\mathcal{G} \Omega f \ \subset \ \mathcal{G}(\Omega f \circ f).
\hspace{2cm}
\end{equation}

%
%From Proposition 9.1(d) it follows that $\Omega f \circ f \subset
%(\Omega f)^2 \subset \mathcal{G} \Omega f$. Hence,
%\begin{equation}
%\mathcal{G} (\Omega f \circ f) \ \subset \ \mathcal{G} \Omega f.
%\end{equation}

By compactness $(\mathcal{G} \Omega f) \circ f $ is a closed
relation and  $ f \circ \mathcal{G} \Omega f = \mathcal{G} \Omega
f$ implies that the relation is transitive as well.  Hence,

\begin{equation}\label{eq9.15}
\mathcal{G} \Omega f \ \subset \ \mathcal{G}(\Omega f \circ f)  \
\subset \ (\mathcal{G}\Omega f) \circ f.
\end{equation}

Now Proposition \ref{prop9.1} (d) with $F =   \mathcal{G}\Omega f $ implies

\begin{equation}\label{eq9.16}
 (\mathcal{G}\Omega f) \circ f  \ \subset \ (\mathcal{G}\Omega f) \circ \Omega
 f  \ \subset \ \mathcal{G}\Omega f .
 \end{equation}

 From (\ref{eq9.15}) and (\ref{eq9.16})  we obtain

 \begin{equation}\label{eq9.17}
 (\mathcal{G}\Omega f) \circ f \quad = \quad    \mathcal{G}\Omega
 f.
 \end{equation}

Thus, (\ref{eq9.10}) holds in the compact case as well.  The proof is
completed as in the homeomorphism case.

$\Box$ \vspace{.5cm}

\begin{ex}\label{ex9.3}When the space is not compact then the results need not hold even
for proper maps which are close to being homeomorphisms.\end{ex}
Let $f_0$ be the time-one map for the
flow on the square $[0,1]^2$ given by

\begin{equation}\label{eq9.18}
\begin{split}
\frac{dx}{dt} \quad = \quad (x(1-x))^2, \hspace{2cm}\\
\frac{dy}{dt} \quad = \quad x(1-x)y(1-y). \hspace{1.7cm}
\end{split}
\end{equation}
For any point $(x,y) \in (0,1) \times (0,1] ,  \Omega f_0 (x,y)
=\omega f_0 (x,y) = \{ (1,1) \}$, for $x \in [0,1], \Omega f_0
(x,0) = \{ 1 \} \times [0,1]$ and for $y \in [0,1], \Omega f_0
(1,y) = \{ (1,y_1) : y \leq y_1 \leq 1 \}$.

Now let $X = (0,1) \times [0,1) \cup \{ 1 \} \times (1/2,1)  \cup
\{ (1/2, -1) \}$. Let $f(1/2,-1) = (1/2,0)$ and otherwise let
$f(x,y) = f_0(x,y)$. This is a proper map on the locally compact
space $X$.  Clearly, $\Omega f(x,y) = \emptyset$ except when $y =
0$ in which case it is the set $ \{ 1 \} \times (1/2,1)$ and when
$x = 1, y \in (1/2,1)$ in which case it it the set $\{ 1 \} \times
[y,1)$.  Since $\Omega f$ is transitive it equals $\mathcal{G}
\Omega f$. On the other hand, $\mathcal{N} f = \mathcal{O} f \cup
\mathcal{G} \Omega f$ is not transitive because $\mathcal{N}f
\circ f (1/2,-1) = \{ 1 \} \times (1/2,1)$. Hence, $\Omega
\mathcal{G}f$ equals $\mathcal{G} \Omega f  \cup \{(1/2,-1) \}
\times \{ (1,y) : 1/2 < y < 1 \}$.

$\Box$ \vspace{1cm}

\section{Appendix: Paracompactness}\label{secparacompact}

While our dynamics results have been stated for $\sigma$-compact, locally
compact Hausdorff spaces they are actually true a bit more generally
because of the following Theorem based on results in Kelley (1955). \index{paracompact space}

\begin{theo}\label{th8.1} A locally compact Hausdorff space is paracompact iff
it admits a partition $\mathcal{Q}$ by clopen, $\sigma$-compact
subsets.

If  $f$ is a proper closed relation on a locally
compact, paracompact space $X$ then $X$ admits a partition
$\mathcal{Q}$ by clopen, $\sigma$-compact subsets each of which is
+ invariant for $f$ and $f^{-1}$.
\end{theo}

{\bfseries Proof:} It is clear from Theorem 5.28 of Kelley (1955)
that any Lindel\"{o}f Hausdorff space, and a fortiori any
$\sigma$-compact Hausdorff space, is paracompact and that if $X$
admits a partition by clopen paracompact subsets then $X$ is
paracompact.

Conversely, if $X$ is a locally compact, Hausdorff space then we
can choose a cover by bounded open sets. If $X$ is also
paracompact then by Kelley Theorem 5.28 again it is even and so
admits a refinement $ \{ V(x) : x \in X \}$ where $V$ is a
neighborhood of the diagonal. By Kelley Lemma 5.30 the
neighborhoods of the diagonal form a uniformity and so we can
choose a closed, symmetric neighborhood $W$ such that  $W \circ W \subset V$.

I claim that $W$ is proper.  By symmetry it suffices to show that $A$ compact implies
$W(A)$ is compact.  From Proposition \ref{prop1.2}(a) it follows that $W(A)$ is closed.  By compactness there
exists  a finite subset $F$ of $A$ such that $A \subset W(F)$.  Then $W(A) \subset W(W(F)) \subset V(F)$.
Since $V(F)$ is bounded, $W(A)$ is compact.

Thus, if
$A_1 \subset X$ is compact then, inductively $A_{n+1} = W(A_n)$ is
compact and $\bigcup_n \{ A_n \} = (\mathcal{O}W)(A_1)$ is
$\sigma$-compact. Thus, $\mathcal{O}W$ is an equivalence
relation on $X$ with $\sigma$-compact equivalence classes.

An equivalence relation $E$ which
contains a neighborhood of the diagonal has open -and hence clopen- equivalence classes. Hence,
$\{ E(x) \times E(y) : x,y \in X \}$ is a clopen partition of $X \times X$ and
so $E = \bigcup \{ E(x) \times E(x) : x \in X \}$ is itself clopen. The equivalence classes form
the required partition of $X$.

See also Kelley Exercises 6L and 6T.

If $f$ is a proper closed relation on $X$ then we let
\begin{equation}\label{eq8.1}
W_f \quad =_{def} \quad W \circ (f \cup 1_X \cup f^{-1}) \circ W.
\end{equation}
By Proposition \ref{prop1.2} (d) the reflexive, symmetric relation $W_f$ is  closed
 and proper. Clearly, $W \subset W_f$.

 As above,   $\mathcal{O} W_f$
 is a  clopen equivalence relation and each equivalence class is
 $\sigma$-compact as well as clopen.  Since $f \cup f^{-1} \subset
 W_f$ it follows that each equivalence class is + invariant for $f
 \cup f^{-1}$.

$\Box$ \vspace{.5cm}

 A. H. Stone's theorem, see Kelley
Corollary 5.35, says that every metric space is paracompact and so
the above applies to every locally compact metric space whether
separable or not.  Of course only a separable metric space admits
a metrizable compactification.
\vspace{.5cm}

\section*{References}

E. Akin (1993) {\bfseries The general topology of dynamical
systems} Amer. Math. Soc., Providence.
\vspace{.5cm}

E. Akin (1997) {\bfseries Recurrence in topological dynamics:
Furstenberg families and Ellis actions} Plenum Publishing Co., New
York.\vspace{.5cm}

 A. Antosiewicz and J. Dugundji (1961)
\emph{Parallelizable flows and Lyapunov's second method} Ann.
Math. {\bfseries 73:}543-555.
\vspace{.5cm}

J. Auslander (1964) \emph{Generalized recurrence in
dynamical systems} Contr. Differential Equations {\bfseries 3:}
65-74.
\vspace{.5cm}

A. Beck  (1958) \emph{On invariant sets} Annals of Math
{\bfseries 67:}99-103.
\vspace{.5cm}

C. Conley  (1978) {\bfseries Isolated invariants sets and the Morse
index} Amer. Math. Soc., Providence.
\vspace{.5cm}

C. Conley (1988)  \emph{The gradient structure of a flow:
I} Ergod. Th. \& Dynam. Sys. {\bfseries 8:}11-26.
\vspace{.5cm}

J. Kelley (1955) {\bfseries General topology} D. Van
Nostrand Co., Princeton.
\vspace{.5cm}

L. Markus (1969) \emph{Parallel dynamical systems}
Topology  {\bfseries 8:}47-57.
\vspace{.5cm}

L. Nachbin (1965) {\bfseries Topology and order} D. Van
Nostrand Co., Princeton.
\vspace{.5cm}

\printindex

\end{document}